\documentclass[12pt]{article}

\usepackage{amsmath,enumerate,amsbsy,amsfonts,amssymb,mathabx,amscd,graphicx,multirow,xcolor,subcaption, pdflscape}

\usepackage{afterpage}
\usepackage{caption}
\usepackage[round,authoryear]{natbib}
\usepackage{authblk}
\usepackage{mathrsfs}
\usepackage[flushleft]{threeparttable}
\usepackage{algorithmic,algorithm}
\RequirePackage[colorlinks,citecolor=blue,urlcolor=blue]{hyperref}
\usepackage{epsfig}

\newcommand{\blind}{1}

\newtheorem{theorem}{Theorem}

\newtheorem{lemma}{Lemma}

\newtheorem{proposition}{Proposition}
\newtheorem{condition}{Condition}
\newtheorem{inequality}{Inequality}

\newtheorem{remark}{Remark}
\usepackage{multirow}
\newcommand{\beps}{\boldsymbol \epsilon}
\newcommand{\bLambda}{\boldsymbol \Lambda}

\newcommand{\bSigma}{\boldsymbol \Sigma}

\newcommand{\bDelta}{\boldsymbol \Delta}
\newcommand{\bXi}{\boldsymbol \Xi}

\newcommand{\bPsi}{\boldsymbol \Psi}
\newcommand{\bphi}{\boldsymbol \phi}
\newcommand{\bPhi}{\boldsymbol \Phi}
\newcommand{\bPi}{\boldsymbol \Pi}

\newcommand{\bxi}{\boldsymbol \xi}
\newcommand{\btheta}{\boldsymbol \theta}

\newcommand{\bbeta}{\boldsymbol \beta}

\newcommand{\bbE}{{\mathbb{E}}}
\newcommand{\bbP}{{\mathbb{P}}}

\newcommand{\md}{{\rm d}}

\newcommand{\be}{{\mathbf e}}

\newcommand{\br}{{\mathbf r}}

\newcommand{\bs}{{\mathbf s}}

\newcommand{\bD}{{\bf D}}
\newcommand{\bA}{{\bf A}}
\newcommand{\bB}{{\bf B}}

\newcommand{\bE}{{\bf E}}
\newcommand{\bI}{{\bf I}}

\newcommand{\bS}{{\bf S}}
\newcommand{\bX}{{\bf X}}
\newcommand{\bY}{{\bf Y}}

\newcommand{\bR}{{\bf R}}
\newcommand{\bU}{{\bf U}}
\newcommand{\bV}{{\bf V}}
\newcommand{\bQ}{{\bf Q}}

\newcommand{\bM}{{\bf M}}

\newcommand{\bcB}{\boldsymbol{\cal B}}

\newcommand{\cA}{{\cal A}}
\newcommand{\cB}{{\cal B}}

\newcommand{\cE}{{\cal E}}
\newcommand{\cF}{{\cal F}}
\newcommand{\cG}{{\cal G}}
\newcommand{\cL}{{\cal L}}

\newcommand{\cH}{{\cal H}}

\newcommand{\cU}{{\cal U}}

\newcommand{\cS}{{\cal S}}

\newcommand{\eR}{\mathbb{R}}

\newcommand{\bbN}{\mathbb{N}}

\newcommand{\tF}{\text{F}}

\newcommand{\cov}{\text{Cov}}

\newcommand{\low}{{\scriptstyle\textup{low}}}
\newcommand{\high}{{\scriptstyle\textup{high}}}

\providecommand{\keywords}[1]{\textbf{\textit{Keywords: }} #1}

\def\T{{ \mathrm{\scriptscriptstyle T} }}

\allowdisplaybreaks

\begin{document}

\def\spacingset#1{\renewcommand{\baselinestretch}%
	{#1}\small\normalsize} \spacingset{1}


\if1\blind
{
	\spacingset{1.25}
	\title{\bf \Large Large-scale multiple testing of cross-covariance functions with applications to functional network models}
	\author[1]{Qin Fang}
	\author[2]{Qing Jiang}
	\author[3]{Xinghao Qiao}
 \affil[1]{\it University of Sydney Business School, Sydney, Australia}
	\affil[2]{Center for Statistics and Data Science, Beijing Normal University, Zhuhai, China}
 \affil[3]{\it Faculty of Business and Economics, The University of Hong Kong, Hong Kong}
	\setcounter{Maxaffil}{0}
	
	\renewcommand\Affilfont{\itshape\small}
	\date{\vspace{-5ex}}
	\maketitle
} \fi

\if0\blind
{\spacingset{2}
	\bigskip
	\bigskip
	\bigskip
	\begin{center}
		{\Large\bf Large-scale multiple testing of cross-covariance functions with applications to functional network models}
	\end{center}
	\medskip
} \fi

\bigskip

\spacingset{1.5}
\begin{abstract}
The estimation of functional networks through functional covariance and graphical models have recently attracted increasing attention in settings with high dimensional functional data, where the number of functional variables $p$ is comparable to, and maybe larger than, the number of subjects. However, the existing methods all depend on regularization techniques, which make it unclear how the involved tuning parameters are related to the number of false edges. In this paper, we first reframe the functional covariance model estimation as a tuning-free problem of simultaneously testing $p(p-1)/2$ hypotheses for cross-covariance functions, and introduce a novel multiple testing procedure. We then explore the multiple testing procedure under a general error-contamination framework and establish that our procedure can control false discoveries asymptotically. Additionally, we demonstrate that our proposed methods for two concrete examples: the functional covariance model for discretely observed functional data and, importantly, the more challenging functional graphical model, can be seamlessly integrated into the general error-contamination framework, and, with verifiable conditions, achieve theoretical guarantees on effective false discovery control. Finally, we showcase the superiority of our proposals through extensive simulations and brain connectivity analysis of two neuroimaging datasets.

\end{abstract}

\noindent%
.

\spacingset{1.68}

\section{Introduction}
Recent advances in information technology have led to the growing prevalence of multivariate or even high-dimensional functional datasets across various applications. Examples include time-course gene expression data in genomics, and various types of brain imaging data in neuroscience, such as electrocorticography, magnetoencephalography and functional magnetic resonance imaging scans, among others. The brain signals are typically recorded in the form of a location-by-time matrix for each individual subject, where the rows and columns correspond to a set of brain regions and a number of observational time points spanning over minutes, respectively. To capture the non-stationary and dynamic patterns, recent proposals involve modeling brain signals as multivariate random functions, treating the time course data of each region as a random function \citep[e.g.,][]{zapata2019,lee2023,xue2023}. Such high-dimensional functional data can be represented as 
$\bX_i(\cdot)=\{X_{i1}(\cdot), \dots, X_{ip}(\cdot)\}^\T$ defined on a compact interval $\cU$ with marginal- and cross-covariance functions, which together form the covariance function matrix 
$$\bSigma(\cdot,\cdot) = \{\Sigma_{jk}(\cdot,\cdot)\}_{p \times p}, ~~~\Sigma_{jk}(u,v)=\cov\{X_{ij}(u), X_{ik}(v)\}~~\text{for}~~(u,v) \in \cU^2.$$
We observe $\bX_i(\cdot)$ for $i=1, \dots, n,$ where the dimension $p$ is large relative to, and maybe greater than the number of subjects $n.$

Our motivation lies in the brain connectivity analysis based on brain imaging data. The literature primarily focuses on two types of multivariate functional data analysis approaches for estimating brain connectivity networks, consisting of $p$ nodes. 
The first method, proposed by \cite{fang2023}, considers a functional covariance model depicting the marginal correlation information in $\bX_i(\cdot).$ This approach aims to identify the functional sparsity structure in $\bSigma(\cdot,\cdot),$ i.e., discovering edges $(j,k)$'s such that $\Sigma_{jk}(u,v)\neq0$ for some $(u,v) \in \cU^2.$ By applying adaptive functional thresholding to the entries of the sample version of $\bSigma(\cdot,\cdot),$ they achieve estimation and support recovery consistencies.
The second method frames the network estimation as a functional graphical modelling problem \cite[]{zhu2016,qiao2019}. This model characterizes the conditional dependence structure of $p$ Gaussian random functions, i.e., nodes $j$ and $k$ are connected by an edge if and only if $X_{ij}(\cdot)$ and $X_{ik}(\cdot)$ are dependent, conditional on the remaining $p-2$ functions. This line of work has witnessed numerous advancements, see, e.g., \cite{li2018,zapata2019,solea2022,lee2023} and \cite{tsai2023}.
Both types of models involve developing the regularized estimation of functional network structures that leverage the functional sparsity information.
However, the selection of practically reasonable regularization parameters poses a challenging task. While large regularization parameters result in sparse networks and may not be effective in discovering edges with small weights, small regularization parameters can produce excessive false edges, leading to high false discovery rates, as evidenced by simulation results in Section~\ref{sec:sim}.

In this paper, we reframe the functional covariance model estimation as a tuning-free problem of simultaneously testing $p(p-1)/2$ hypotheses for cross-covariance functions
\begin{equation}
\label{test_cov}
H_{0,jk}: \Sigma_{jk}(u,v)=0 \text{ for any } (u,v) \in \cU^2~~\text{vs}~~H_{1,jk}: \Sigma_{jk}(u,v) \neq 0 \text{ for some } (u,v) \in \cU^2,
\end{equation}
where $1 \leq j < k \leq p$ and we identify a significant edge between nodes $j$ and $k$ if and only if $H_{0,jk}$ is rejected. 
Compared to the multiple testing for the local covariance structures of high-dimensional vectors \cite[]{Cai2017}, the infinite-dimensional nature of functional data introduces additional methodological and theoretical complexities that arise from constructing the test statistic under a functional norm and addressing various technical obstacles within an abstract Hilbert space. Furthermore, when handling the practical scenario for brain imaging data with each trajectory $X_{ij}(\cdot)$ being observed over a densely sampled grid, nonparametric smoothing is frequently employed to obtain estimated curves $\widehat X_{ij}(\cdot),$ which requires considering functional error-contaminated versions of $X_{ij}(\cdot)$ satisfying: 
\begin{equation}
\label{error.X}
\widehat X_{ij}(\cdot) = X_{ij}(\cdot) + e_{ij}(\cdot),~~i=1, \dots, n,~j=1, \dots, p.
\end{equation}
Additionally, we present a novel procedure based on nodewise functional regressions that converts the functional graphical model estimation into the multiple testing task for the cross-covariance structures of functional regression errors.
As a result, we need to deal with estimated functional regression errors (i.e., functional residuals) instead of true ones, following a similar form as (\ref{error.X}). In both scenarios, accounting for the functional errors, such as $e_{ij}(\cdot)$'s in (\ref{error.X}), introduces an extra layer of complexity to our theoretical analysis.

Our paper aims to develop a general framework for large-scale multiple testing of cross-covariance functions from both methodological and theoretical perspectives. We begin by constructing a Hilbert--Schmidt-norm-based test statistic for each pair $j<k,$ which is shown to have asymptotically the same limiting null distribution as an infinite mixture of 
chi-squares. We then employ normal quantile transformations for all test statistics, upon which a multiple testing procedure is proposed to account for the multiplicity and dependence among the transformed test statistics. Theoretically, we establish that our procedure can control false discoveries asymptotically under both fully observed and error-contaminated functional scenarios, and, furthermore, be applied to two concrete examples: the functional covariance model with discrete observations and the more challenging functional graphical model.
Specifically, we demonstrate the seamless integration of our proposed methods for both examples into the general error-contamination framework, and, with verifiable conditions, achieve theoretical guarantees on false discovery control. Empirically, we conduct simulations to showcase the uniform superiority of our proposals over competitors in terms of false discovery control and power for both fully and discretely observed functional data within both functional covariance and graphical models. We also apply our method to identify network structures using two brain imaging datasets, and observe scientifically interpretable patterns. 

The main contributions of our paper are threefold.

\begin{itemize}
\item First, we make the first attempt in the literature of functional data analysis and multiple testing to develop a general procedure with theoretical guarantees for the simultaneous testing of a large collection of hypotheses for the functional sparse covariance structures. Our approach is fully functional in the sense it does not rely on dimension reduction techniques, thereby avoiding any incurred information loss.

\item Second, we extend our method and theory to the more general error-contamination setting (\ref{error.X}), encompassing the common practical scenario of discretely observed functional data as a special case. On the method front, our procedure remains valid by simply replacing each $X_{ij}(\cdot)$ with its estimated surrogate $\widehat X_{ij}(\cdot).$ Theoretically, we demonstrate that, under an additional condition, our proposal still ensures control over false discoveries. Such condition can be verified under mild circumstances for discretely observed functional data.

\item Third, we propose a novel method to formulate the functional graphical model estimation as a multiple testing problem for the cross-covariance structures between $p(p-1)/2$ pairs of functional regression errors formed by respectively regressing each pair $X_{ij}(\cdot)$ and $X_{ik}(\cdot)$ ($1 \leq j< k \leq p$) on the remaining $p-2$ functional variables. Employing the penalized functional regression technique to estimate functional coefficients, we obtain 
$p(p-1)/2$ pairs of functional residuals, which can be nicely integrated into a generalized version of our error-contamination framework. By deriving relevant convergence rates to validate the extra condition and applying our established theory, we show that such proposal combined with our multiple test procedure achieves effective false discovery control.
\end{itemize}

Our work lies in the intersection of high-dimensional functional data analysis and large-scale multiple testing.
In addition to aforementioned functional covariance and graphical models, various functional sparsity assumptions are commonly imposed on the model parameter space for high-dimensional functional data, and different regularized estimation procedures have been developed for the respective learning tasks including functional additive regressions \cite[]{fan2015,kong2016,luo2017,chang2023a}, functional linear discriminant analysis \cite[]{xue2023} and sparse functional principal components analysis \cite[]{hu2021}.
There is also a wealth of literature focused on the simultaneous testing of a large collection of hypotheses with false discovery control for the local covariance structures. Examples include, but are not limited to, studies on correlations \cite[]{cai2016}, conditional dependence between variables or subvectors in Gaussian graphical models \cite[]{Liu2013,Xia2018b}, regression coefficients over multiple responses \cite[]{xia2018a}, cross-covariance matrices between subvectors in multimodal integrative analysis \cite[]{Xia2020} and more complex dependence structures of multimodal imaging data \cite[]{Chang2023b}.

The rest of the paper is organized as follows. Section~\ref{sec:method} presents the test for a given pair of cross-covariance functions and a multiple testing procedure for all pairs, along with the established theoretical properties. Section~\ref{sec:error} explores the multiple testing under a general error-contamination framework and presents the supporting theory. Section~\ref{sec:app} illustrates the proposed methods with theoretical guarantees using two functional network models. We demonstrate the superior finite-sample performance of the proposed methods through extensive simulations in Section~\ref{sec:sim} and the analysis of two brain imaging datasets in Section~\ref{sec:real}. 
All technical proofs are relegated to the Supplementary Material.

We summarize some notation to be used throughout the paper. Denote by $I(\cdot)$ the indicator function. Let $\overset{d}{=}$ denote the  equality in distribution. For a positive integer $q,$ we write $[q]=\{1,\ldots,q\}$.
 For any $x,y\in\mathbb{R}$, let $x\vee y=\max(x,y)$ and $(x)_{+} = x \vee 0$. For any set $\cA$, we denote its cardinality by $|\cA|.$ For two sequences of positive numbers $\{a_n\}$ and $\{b_n\}$, we write $a_n\lesssim b_n$ or $b_n\gtrsim a_n$ if there exists some positive constant $c$ such that $\lim\sup_{n\to \infty}a_n/b_n \leq c$, and $a_n\ll b_n$ if $\lim\sup_{n\rightarrow\infty}a_n/b_n=0$. 
We write $a_n\asymp b_n$ if and only if $a_n\lesssim b_n$ and $b_n\lesssim a_n$ hold simultaneously. 
 For a vector ${\bf b}=(b_1,\ldots,b_p)^{\T}$, we denote its $\ell_2$ norm by $\|{\bf b}\|=(\sum_{j=1}^p|b_j|^2)^{1/2}$.  For a matrix $ \bB = (B_{ij})_{p \times q}$, we denote its Frobenius norm by $\|\bB\|_{\tF}=\big(\sum_{i,j}  B_{ij}^2\big)^{1/2}.$ 
Let $\cL_2(\cU)$ be a Hilbert space of square-integrable functions on a compact interval $\cU$ equipped with the inner-product $\langle f,g \rangle = \int_{\cU}f(u)g(u)\,\md u$ for $f,g\in\cL_2(\cU)$ and the induced norm $\|\cdot\|=\langle\cdot,\cdot\rangle^{1/2}$. 
For a Hilbert space $\mathbb{H} \subseteq L_2(\cU),$ we denote the $p$-fold Cartesian product by $\mathbb{H}^p = \mathbb{H} \times \cdots \times \mathbb{H}$.
For any $\cB\in\mathbb{S} \equiv L_2(\cU \times \cU)$,  we denote its Hilbert--Schmidt norm by $\|\cB\|_{\cS}=\{\iint\cB(u,v)^2\,{\rm d}u{\rm d}v\}^{1/2}$. For any $\bcB=(\cB_{1},\ldots,\cB_p)^{\T}$ with each $\cB_{j}\in\mathbb{S}$, we denote its functional versions of 
$\ell_1$ and $\ell_{\infty}$ norms by $\|\bcB\|_{\cS,1}=\sum_{j=1}^p\|\cB_{j}\|_{\cS}$ and $\|\bcB\|_{\cS,\max}=\max_{j\in[p]}\|\cB_{j}\|_{\cS},$ respectively.

\section{Test with fully observed functional data}
\label{sec:method}
This section considers the setting where each trajectory $X_{ij}(\cdot)$ for $i \in[n]$ and $j \in [p]$ is fully observed without any error. In Section~\ref{sec:test_full}, we propose a marginal test statistic for (\ref{test_cov}) with a given pair $1 \leq j \neq k \leq p,$ and study its theoretical properties. 
In Section~\ref{sec:mt}, we develop a procedure for simultaneously testing multiple hypotheses of 
(\ref{test_cov}) for $1 \leq j <k \leq p$ with theoretical guarantees on false discovery control.

\vspace{-0.3cm}
\subsection{Testing a cross-covariance function}
\label{sec:test_full}
Based on $n$ i.i.d. functional observations $\bX_{1}(\cdot), \dots, \bX_n(\cdot),$ we obtain the sample covariance function matrix $\widehat\bSigma(\cdot,\cdot)=\{\widehat\Sigma_{jk}(\cdot,\cdot)\}_{p \times p}$ with its $(j,k)$th entry given by
\begin{equation}\label{eq:sigmajk}
	\widehat{\Sigma}_{jk}(u,v)=\frac{1}{n-1}\sum_{i=1}^n\{X_{ij}(u)-\widebar{X}_j(u)\}\{X_{ik}(v)-\widebar{X}_{k}(v)\} ~\text{for}~(u,v) \in \cU^2~\text{and}~j,k \in [p]
\end{equation}
and $\widebar{X}_j(\cdot)=n^{-1}\sum_{i=1}^n X_{ij}(\cdot)$ for $j\in[p].$ 



We propose the following Hilbert--Schmidt-norm-based test statistic for testing $H_{0,jk}$ with a given pair $j \neq k,$ i.e.,
\begin{equation}
\label{eq:statH0}
	T_{n,jk}=n \|\widehat{\Sigma}_{jk} \|_{\cS}^2\,.
\end{equation}
\vspace{-1.0cm}
\begin{remark}
Note that the test statistic $T_{n,jk}$ measures the discrepancy between the sample cross-covariance function $\widehat\Sigma_{jk}(u,v)$ and the true function using an $L_2$-type norm, which tends to be small under the null hypothesis $H_{0,jk}.$ For functional data with local spikes, one may consider a supremum-norm-based test statistic, upon which our multiple testing procedure in Section~\ref{sec:mt} remains valid under the associated technical analysis. See more detailed discussions in Section~\ref{sec:discuss}. By comparison, our $L_2$ type-based test performs better when a (eventually not so strong) signal in the cross-covariance function is spread out over the entire domain $\cU^2,$ making it well-suited for accommodating the brain imaging data.
See Figures~\ref{HCP_cov.est} and \ref{EEG_fgm.est} of the Supplementary Material for the supporting empirical evidence.
\end{remark}

Denote the $j$th subject-effect function $Z_{ij}(\cdot)=X_{ij}(\cdot)-\mu_j(\cdot)$ with the mean function $\mu_j(\cdot)=\mathbb E\{X_{ij}(\cdot)\}$ for $j \in [p].$
Denote $\Gamma_{jk}(u_1,v_1,u_2,v_2)=\cov\{Z_{j}(u_1)Z_k(v_1),Z_j(u_2)Z_k(v_2)\}$ for $(u_1,v_1), (u_2, v_2) \in \cU^2$, which can be shown to be the asymptotic covariance function of $\sqrt{n}\widehat{\Sigma}_{jk}(u,v)$ for $(u,v) \in \cU^2.$ By performing an eigenanalysis of $\Gamma_{jk}(\cdot,\cdot,\cdot,\cdot)$, we obtain the sorted eigenvalues
$\lambda_{jk1} \geq \lambda_{jk2} \geq \cdots > 0$ and the associated orthonormal eigenfunctions $\varphi_{jk1}(\cdot,\cdot), \varphi_{jk2}(\cdot,\cdot), \dots$ satisfying
\begin{equation}
\label{eigenanalysis}
\iint\Gamma_{jk}(u_1,v_1,u_2,v_2)\varphi_{jkr}(u_2,v_2)\,{\rm d}u_2{\rm d}v_2 = \lambda_{jkr} \varphi_{jkr}(u_1,v_1) ~~\text{for}~~r\in \mathbb N_+.
\end{equation}
We compute the sample estimator of $\Gamma_{jk}$ by
\begin{align*}
    \widehat{\Gamma}_{jk}(u_1,v_1,u_2,v_2) = \frac{1}{n}\sum_{i=1}^n \{\widetilde{X}_{ij}(u_1)\widetilde{X}_{ik}(v_1)-\widehat{\Sigma}_{jk}(u_1,v_1)\}\{\widetilde{X}_{ij}(u_2)\widetilde{X}_{ik}(v_2)-\widehat{\Sigma}_{jk}(u_2,v_2)\} 
\end{align*}
with $\widetilde{X}_{ij}(\cdot)=X_{ij}(\cdot)-\widebar{X}_j(\cdot)$ and then carry out an eigenanalysis of $\widehat\Gamma_{jk}$ in the form of (\ref{eigenanalysis}), which leads to the estimated eigenvalue/eigenfunction pairs $\{\hat \lambda_{jkr}, \hat \varphi_{jkr}(\cdot,\cdot)\}_{r \in \mathbb N_+}.$

Before investigating the theoretical properties, we need some regularity conditions.

\begin{condition}\label{cond:covfun}
 {\rm (i)} $\max_{j\in[p]}\int\Sigma_{jj}(u,u)\,\md u=O(1);$
 {\rm (ii)}$\sup_{u\in\cU}\Sigma_{jj}(u,u)=O(1)$ for any $j\in[p]$.
\end{condition}
\begin{condition}\label{cond:eta}
$\bbE(\|Z_{1j}\|^4)=O(1)$ for any $j \in [p]$ and $\sup_{(u,v) \in \cU^2}\bbE\{Z_{1j}^2(u)Z_{1k}^2(v)\} = O(1)$ for any $j,k \in [p].$
\end{condition}

Condition~\ref{cond:covfun} (i) implies that $\int\Sigma_{jj}(u,u)\md u=O(1)$ for any $j \in [p],$ which together with Conditions~\ref{cond:covfun} (ii) and \ref{cond:eta} are standard in functional data analysis literature, see, e.g., Chapter 10 of \cite{Zhang2013}. 
We are now ready to present the theorem regarding the limiting null distribution and the asymptotic power.
\begin{theorem}\label{thm:nonH0} 
Suppose that Conditions {\rm \ref{cond:covfun} and \ref{cond:eta}} hold, then the following two assertions hold.\\
{\rm (i)} Under the null hypothesis $H_{0,jk}$, we have that
	$T_{n,jk}\overset{d}{\rightarrow}T_{0,jk}$ as $n\rightarrow\infty,$
	where
		$T_{0,jk}\overset{d}{=}\sum_{r=1}^{\infty} \hat\lambda_{jkr}A_{r}$ with $\{A_{r}\}_{r \in \mathbb N_+}$ being i.i.d. $\chi^2_1$ random variables. \\
 {\rm (ii)} Under the alternative hypothesis $H_{1,jk}$ with $\|\Sigma_{jk}\|_{\cS} >0$, we have that $\bbP\{T_{n,jk}\geq T_{0,jk}(\alpha)\}\rightarrow 1$ as $n \rightarrow \infty$, where $T_{0,jk}(\alpha)$ is the $100\alpha\%$ upper percentile of $T_{0,jk}$.
\end{theorem}


\begin{remark}
Theorem~\ref{thm:nonH0} indicates that $T_{n,jk}$ is asymptotically an infinite $\chi^2$-type mixture under the null hypothesis $H_{0,jk}.$ Its distribution can be approximated by a noncentral chi-squared distribution with parameters determined by $\sum_{r=1}^{\infty}\hat\lambda_{jkr},$ $\sum_{r=1}^{\infty}\hat\lambda_{jkr}^2,$ $\sum_{r=1}^{\infty}\hat\lambda_{jkr}^3$ and $\sum_{r=1}^{\infty}\hat \lambda_{jkr}^4$, see, e.g., \cite{Liu2009}. We will use this approximation in our numerical studies. It is noteworthy that, instead of calculating each individual eigenvalue $\hat\lambda_{jkr}$ explicitly, it is more convenient to compute the aforementioned infinite sums. For example, $\sum_{r=1}^{\infty}\hat\lambda_{jkr}$ and $\sum_{r=1}^{\infty}\hat\lambda_{jkr}^2$ can be obtained through the discretization of functional data and the respective simple Riemann sum approximations to the integrals $\iint \widehat\Gamma_{jk}(u,v,u,v)\,{\rm d}u {\rm d}v$ and $\iiiint\widehat\Gamma_{jk}(u_1,v_1,u_2,v_2)^2\,{\rm d}u_1 {\rm d}v_1 {\rm d}u_2 {\rm d}v_2.$ In this sense, our procedure is fully functional without relying on any dimension reduction techniques that may result in information loss.
\end{remark}



\subsection{Multiple testing of cross-covariance functions}
\label{sec:mt}

To develop a general multiple testing procedure, we let $\cH=\{(j,k):j,k\in[p],j< k\}$ and reformulate (\ref{test_cov}) as simultaneously testing $Q=p(p-1)/2$ hypotheses:
\begin{align}\label{eq:mt}
	H_{0}^{(q)}:\Sigma_{j_qk_q}(u,v)=0 \mbox{ for any } (u,v)\in\cU^2 \mbox{ and } (j_q,k_q)\in\cH
	~~\mbox{vs}~~
	H_{1}^{(q)}:H_{0}^{(q)} \mbox{ is not true}\,,
\end{align}
for $q\in[Q]$. The test statistic of the $q$th hypothesis test $H_0^{(q)}$ is constructed as 
\begin{equation}
\label{T.stat}
	T_n^{(q)} :=T_{n,j_qk_q}= n\|\widehat{\Sigma}_{j_qk_q} \|_{\cS}^2~,~~q \in [Q].
\end{equation}
Then we reject $H_{0}^{(q)}$ when $T_n^{(q)}$ takes some large values. 
Let $\cH_0=\{q\in [Q]:H_{0}^{(q)}~\textrm{is true}\}$ 
be the set of all true null hypotheses with $Q_0=|\cH_0|.$
By Theorem \ref{thm:nonH0}(i), we obtain that
\begin{equation}
\label{prob.conv}
\bbP\{T_n^{(q)}\leq x\} - \bbP\{T_0^{(q)}\leq x\} = o(1) 
\text{  for any  } q\in\cH_0,
\end{equation} 
where $T_0^{(q)}\overset{d}{=}\sum_{r=1}^{\infty}\hat\lambda_{r}^{(q)} A_r$ with 
$\hat{\lambda}_{r}^{(q)}:=\hat{\lambda}_{j_qk_qr}$ for notational simplicity. 
Then the p-value of $H_{0}^{(q)}$ is given by ${\rm pv}^{(q)}=\bbP\{T_0^{(q)}\geq T_n^{(q)}\}$, and its normal quantile transformation takes the form $V_q=\Phi^{-1}\{1-{\rm pv}^{(q)}\},$ where $\Phi(\cdot)$ is the cumulative distribution function of the standard normal distribution. 
For any threshold level $t \in {\mathbb R}$ such that $H_{0}^{(q)}$ is rejected if $V_q\geq t,$ denote the total number of false positives by $R_0(t)=\sum_{q\in\mathcal{H}_0}I(V_q \geq t)$ and the total number of rejections by $R(t)=\sum_{q=1}^Q I(V_q \geq t)$. Then, the false discovery proportion (FDP) and the false discovery rate (FDR) are respectively defined as
\begin{align*}
	{\rm FDP}(t)=\frac{R_0(t)}{R(t)\vee1}~~~\textrm{and}~~~{\rm FDR}(t)=\mathbb{E}\{{\rm FDP}(t)\}\,.
\end{align*}

Given a prespecified level $\alpha\in(0,1)$, the main target is to find the smallest $\hat t$ such that ${\rm FDP}(\hat t) \leq \alpha.$
To achieve this, we first consider ${\rm FDP}(t),$ the numerator of which needs to be estimated due to the unknown set $\cH_0.$ Using (\ref{prob.conv}), we can approximate $\bbP(V_q \geq t)$ as $1-\Phi(t)$ for any $q\in\mathcal{H}_0.$ An ideal estimate for ${\rm FDR}(t)$ can be obtained by 
$[Q_0\{1-\Phi(t)\}]/\{1\vee \sum_{q=1}^Q I(V_q\geq t)\}$, 
where $Q_0$ is unknown in practice. Thus we propose to estimate ${\rm FDP}(t)$ in a more conservative way by replacing $Q_0$ with $Q:$
\begin{equation}
\label{est.FDP}
\widehat{\rm FDP}(t) = \frac{Q\{1-\Phi(t)\}}{1\vee \sum_{q=1}^Q I(V_q\geq t)}. 
\end{equation}
We summarize the proposed multiple testing procedure in Algorithm~\ref{alg1}.

\begin{algorithm}[H] \caption{Multiple testing procedure} \label{alg1}
	\vspace{0.3em}
	\begin{algorithmic}
		\STATE\hspace{-1.3em}
		{\bf Step 1.} Calculate the test statistics $T_n^{(q)}$ in (\ref{T.stat}) and 
its transformation $V_q$ for all $q\in[Q]$.
		\STATE\hspace{-1.3em}
		{\bf Step 2.} Estimate the FDP by $\widehat{\rm FDP}(t)$ in (\ref{est.FDP}).
		\STATE\hspace{-1.3em}
		{\bf Step 3.} For a given $0< \alpha< 1$, choose
             \vspace{-1em}
		\begin{align*}
			\hat{t} = \inf\big\{0\leq t\leq (2\log Q-2\log\log Q)^{1/2}: \, \widehat{\rm FDP}(t)\leq \alpha\big\} \,.
		\end{align*}
		\vspace{-3em}
		\STATE\hspace{2.7em}
		If $\hat{t}$ does not exist, set $\hat{t}=(2\log Q)^{1/2}$.
		\STATE\hspace{-1.3em}
		{\bf Step 4.} For $q\in[Q]$, reject $H_{0}^{(q)}$ with $V_q\geq\hat{t}.$ 
		\vspace{0.08in}
	\end{algorithmic}
\end{algorithm}


Before presenting the theoretical results, we impose some regularity conditions. 
\begin{condition}\label{cond:subG}
For each $j \in [p]$ and $i \in [n],$ $Z_{ij}(\cdot)$ is a mean-zero sub-Gaussian process, i.e., there exists some constant $\check c_j> 0$ such that for all $ x \in \mathbb{H},$
$
    \mathbb{E}
    \big\{\exp(\langle x, Z_{ij} \rangle)\big\} \leq \exp\big\{2^{-1}\check c_j^2 \iint x(u)\Sigma_{jj}(u,v)x(v)\rm{d}u{\rm d}v\big\}.
$
\end{condition}
 

To accommodate the formulation in this section, for any $q  \in \mathcal{H}_0,$ the Karhunen-Lo\`eve expansion allows us to represent 
$z_{iq}(u,v):=Z_{ij_q}(u)Z_{ik_q}(v)=\sum_{r=1}^{\infty}a_{ir}^{(q)}\varphi_{r}^{(q)}(u,v),$ where $\lambda_r^{(q)}:=\lambda_{j_qk_qr}$, $\varphi_r^{(q)}(u,v):=\varphi_{j_qk_qr}(u,v)$, and $a_{ir}^{(q)} =\iint z_{iq}(u,v) \varphi_{r}^{(q)}(u,v)\,{\rm d}u{\rm d}v.$  


\begin{condition}\label{con_eigenvalue}
For any $q  \in \mathcal{H}_0,$  $\lambda_{1}^{(q)} > \lambda_{2}^{(q)}> \cdots >0,$ 
and there exist some constant $\eta_1 > 1$ such that $\lambda_{r}^{(q)} \lesssim r^{-\eta_1}$ for $r \in\mathbb{N}_+.$ 
\end{condition}

Condition~\ref{cond:subG} is an infinite-dimensional analog of the sub-Gaussian condition within Hilbert space, which is satisfied by Gaussian processes. This condition ensures that $z_{iq}(\cdot,\cdot)$ exhibits sub-exponential behavior, as reflected in the associated basis coefficients $\{a_{ir}^{(q)}\}_{r\in {\mathbb N_+}}.$ See Lemma~\ref{lm_norm} of the Supplementary Material. These sub-exponentialities facilitate the Gaussian approximations of the proposed test statistics. 
Condition~\ref{con_eigenvalue} allows the polynomial decay of the upper bounds on eigenvalues with parameter $\eta_1$ determining the decay rate, and can be used to control the truncation error rate 
$\sum_{r=M+1}^{\infty}\lambda_{r}^{(q)} \lesssim M^{-\eta_1+1}$ 
for $M \rightarrow \infty.$ Larger values of $\eta_1$ correspond to faster decay rate and smaller truncation errors. It can be easily verified that Condition~\ref{con_eigenvalue} is fulfilled when 
the eigenvalues for the covariance function of $Z_{ij_q}(\cdot)$ decay polynomially \cite[]{kong2016,qiao2019}.


To investigate the theoretical properties of the proposed multiple testing procedure, we need to give some regularity conditions on the dependency among the test statistics $\{T_n^{(q)}\}_{q=1}^Q$. Notice that $\zeta_q=\Phi^{-1}[F_q\{T_n^{(q)}\}]$ follows the standard normal distribution, where $F_q(\cdot)$ denotes the distribution of $T_n^{(q)}$. We use the correlation between two normal random variables $\zeta_q$ and $\zeta_{q'}$ to characterize the dependence between test statistics $T_n^{(q)}$ and $T_n^{(q')}$. 
For some constant $\gamma>0$, define the set
$$	\mathcal{S}_q(\gamma) = \{q'\in[Q]: q'\neq q, |{\rm Corr}(\zeta_q,\, \zeta_{q'})| \geq(\log Q)^{-2-\gamma}\} \,.$$
For a given $q \in [Q]$, the remaining $Q-1$ test statistics 
$\{ T_n^{(q')}\}_{q'\in[Q]/\{q\}}$ can be categorized into two scenarios: 
(i) for $q' \in \mathcal{S}_q(\gamma)$, 
$T_n^{(q')}$ has a relatively strong dependence with $T_n^{(q)},$ 
and (ii) for $q'\notin \mathcal{S}_q(\gamma)$, 
$T_n^{(q')}$ has a quite weak dependence with $T_n^{(q)}$. 
We establish the theoretical guarantee under these two scenarios separately using different technical tools.

\begin{theorem}\label{thm:nonFDR}
Suppose that Conditions~{\rm \ref{cond:covfun}--\ref{con_eigenvalue}} hold,   $\max_{1\leq q\neq q'\leq Q}|{\rm Corr}(\zeta_q,\zeta_{q'})|\leq c_\zeta$ for some constant $c_{\zeta}\in(0,1),$ and $\max_{q\in[Q]}|\mathcal{S}_q(\gamma)|=o(p^{2\nu})$ for some constants $\gamma>0$ and $0<\nu<(1-c_{\zeta})/(1+c_{\zeta})$. If 
 $p\lesssim n^\kappa$ for some constant $\kappa>0$, then $\limsup_{n,Q\rightarrow\infty}{\rm FDR}(\hat t)\leq \alpha Q_0/Q$ and $\lim_{n,Q\rightarrow\infty}\bbP\{{\rm FDP}(\hat t)\leq \alpha Q_0/Q+\varkappa\}=1$ for any $\varkappa>0$.
\end{theorem}

Theorem~\ref{thm:nonFDR} shows that both FDP and FDR can be controlled below the level $\alpha Q_0/Q$ asymptotically.
The conditions on ${\rm Corr}(\zeta_q,\zeta_{q'})$ and $\mathcal{S}_q(\gamma)$ are imposed to bound the variance of $\sum_{q\in\cH_0}I(V_q\geq t)$ in ${\rm FDP}(t),$ see similar conditions in \cite{Liu2013} and \cite{Chang2023b}. Specifically, the condition $\max_{1\leq q\neq q'\leq Q}|{\rm Corr}(\zeta_q,\zeta_{q'})|\leq c_{\zeta}$ places a constraint on the strength of dependence between different $\zeta_q$ and $\zeta_{q'}$. 
The condition $\max_{q\in[Q]}|\mathcal{S}_q(\gamma)|=o(p^{2\nu})$ controls the number of pairs of $(\zeta_q,\zeta_{q'})$ with relatively strong dependencies. 



\section{Test with error-contaminated functional data} 
\label{sec:error}

In practical applications, $\bX_i(\cdot)$ for $i \in [n]$ are rarely fully observed. Instead, we estimate them using certain methods, resulting in estimated curves $\widehat{\bX}_{i}(\cdot)=\{\widehat{X}_{i1}(\cdot), \ldots, \widehat{X}_{ip}(\cdot)\}^{\T},$ denoted as error-contaminated functional data satisfying (\ref{error.X}). This leads to functional errors  $\be_i(\cdot) = \widehat{\bX}_i(\cdot) - \bX_i(\cdot)$ with $\be_i(\cdot)=\{e_{i1}(u), \ldots, e_{ip}(\cdot)\}^{\T}.$ See below for two examples that are related to applications in Sections~\ref{sec.part} and \ref{sec.fgm}, respectively.

\begin{itemize}
\item When $\{\bX_i(\cdot)\}$ are discretely observed with errors, we can apply nonparametric smoothing methods to the observed data, thereby obtaining functional estimates $\widehat{\bX}_i(\cdot).$
\item When $\{\bX_i(\cdot)\}$ represent the functional regression errors in function-on-function linear regressions, we can estimate functional coefficients and thus obtain functional fitted values, resulting in functional residuals $\widehat\bX_i(\cdot).$
\end{itemize}
In the above examples, it is apparent that each $\be_i(\cdot)$ may not be independent, and the properties of $\be_i(\cdot)$ need to be investigated for specific models. This section is devoted to the corresponding tests for the cross-covariance functions $\{\Sigma_{jk}(u,v)\}_{1 \leq j \neq k \leq p}$ based on the estimated curves $\widehat{\bX}_i(\cdot)$  under some general conditions.

Analogous to Section \ref{sec:test_full}, we first consider the marginal testing problem \eqref{test_cov} for error-contaminated functional data.
To this end, we replace each $\bX_i(\cdot)$ by its error-contaminated version $\widehat\bX_i(\cdot)$ in \eqref{eq:sigmajk}--\eqref{eq:statH0}, and propose a new test statistic for $H_{0,jk}$ for a given $1 \leq j \neq k \leq p:$
\begin{align}\label{eq:statH0d}
\widetilde{T}_{n,jk}=n\|\widetilde{\Sigma}_{jk}\|_{\cS}^2 \,,
\end{align}
where $ \widetilde{\Sigma}_{jk}(u,v)=(n-1)^{-1}\sum_{i=1}^n \{\widehat{X}_{ij}(u)-{\widebar{\widehat{X}}}_j(u)\}\{\widehat{X}_{ik}(v)-\widebar{\widehat{X}}_k(v)\} $ and $\widebar{\widehat{X}}_j(u)=n^{-1}\sum_{i=1}^n \widehat{X}_{ij}(u).$
To control the difference between $\widetilde{T}_{n,jk}$ and $T_{n,jk}$ specified in \eqref{eq:statH0}, we need the cross-covariance functions $\bSigma^{X e}(u,v)={\cov}\{\bX(u),\be(v)\}=\{\Sigma_{jk}^{X e}(u,v)\}_{j,k\in[p]}$ and  $\bSigma^{e e}(u,v)={\cov}\{\be(u),\be(v)\}=\{\Sigma_{jk}^{e e}(u,v)\}_{j,k\in[p]},$ whose sample versions are
\begin{align*}
    \widehat\bSigma^{X e}(u,v)=\{\widehat{\Sigma}_{jk}^{X e}(u,v)\}_{j,k\in[p]} =&~ \frac{1}{n-1}\sum_{i=1}^n\{\bX(u)-\widebar\bX(u)\}\{\be(v)-\bar\be(v)\}^{\T}  \,, \\
    \widehat\bSigma^{e e}(u,v)=\{\widehat{\Sigma}_{jk}^{e e}(u,v)\}_{j,k\in[p]}  =&~ \frac{1}{n-1}\sum_{i=1}^n\{\be(u)-\bar\be(u)\}\{\be(v)-\bar\be(v)\}^{\T} \,,
\end{align*}
respectively.
Before presenting the asymptotic properties of the new test statistic $\widetilde T_{n,jk}$ under the null and alternative hypotheses, we need the following condition, which will be verified by deriving the corresponding elementwise maximum rates of convergence for two applications in Section~\ref{sec:app}.

\begin{condition}
\label{cond:Xe}
There exist some constants $a_1, a_3 >1/2$ and $a_2,a_4>0$ such that
\begin{align}\label{cond.cov}
& \max_{j,k\in[p]}\|\widehat\Sigma^{X e}_{jk}\|_{\cS} = O_{\rm p}\big\{n^{-a_1}(\log p)^{a_2}\big\} 
~~{\rm and}~~
 \max_{j,k\in[p]}\|\widehat\Sigma^{e e}_{jk}\|_{\cS}  =O_{\rm p}\big\{n^{-a_3}(\log p)^{a_4}\big\}.
\end{align}
\end{condition}
\begin{theorem}\label{thm:errH0}
Suppose that Conditions~{\rm \ref{cond:covfun},~\ref{cond:eta}},~{\rm \ref{cond:Xe}} hold, 
then the following two assertions hold.
{\rm (i)} Under the null hypothesis $H_{0,jk}$, we have that
$\widetilde{T}_{n,jk}\overset{d}{\rightarrow}T_{0,jk}$ as $n\rightarrow\infty,$
where $T_{0,jk}\overset{d}{=}\sum_{r=1}^{\infty} \hat\lambda_{jkr}A_{r}$ with $A_{r}\overset{i.i.d.}{\sim}\chi^2_1$.\\
{\rm (ii)} Under the alternative hypothesis 
$H_{1,jk}$ with $\|\Sigma_{jk}\|_{\cS}>0$,  
we have $\bbP\{\widetilde{T}_{n,jk}\geq T_{0,jk}(\alpha)\}\rightarrow 1$ as $n \rightarrow \infty.$
\end{theorem}

%

Analogous to Section \ref{sec:mt}, we next consider the multiple testing problem \eqref{eq:mt} under the functional error-contamination setting. For each $q\in[Q]$, we propose the test statistic for the $q$th hypothesis testing as follows:
\begin{align*}
	\widetilde{T}_n^{(q)}:= \widetilde{T}_{n,j_qk_q} = n\| \widetilde{\Sigma}_{j_qk_q}\|_{\cS}^2\,.
\end{align*}
By Theorem \ref{thm:errH0}(i), we have
$\bbP\{\widetilde{T}_n^{(q)}\leq x\} - \bbP\{T_0^{(q)}\leq x\} = o(1) $
for any $q\in\cH_0$. 
Let the resulting p-value of $H_0^{(q)}$ be $\widetilde{\rm pv}^{(q)}=\bbP\{T_0^{(q)}\geq \widetilde{T}_n^{(q)}\}$, its normal quantile transformation $\widetilde{V}_q=\Phi^{-1}\{1-\widetilde{\rm pv}^{(q)}\},$ and
${t}$ be the threshold level such that $H_{0}^{(q)}$ is rejected if $\widetilde{V}_q\geq {t}$.
Given $\alpha\in(0,1)$, the ideal choice of ${t}$ is the smallest ${t}$ such that 
we would reject as many true positives as possible while controlling the false discovery proportion at the prescribed level $\alpha$. To achieve this, similar to (\ref{est.FDP}) and Step~3 in Algorithm \ref{alg1}, we let
$\widetilde{{\rm FDP}}(t)=Q\{1-\Phi(t)\}/\{1\vee \sum_{q=1}^Q I( \widetilde{V}_q \geq t)\}$ and
\begin{equation*}
\label{t.tilde}
		\tilde{t} = \inf\{0\leq t\leq (2\log Q-2\log\log Q)^{1/2}: \, \widetilde{\rm FDP}(t)\leq \alpha\}.
\end{equation*}
We then propose a new algorithm to implement the multiple testing procedure by replacing $T_n^{(q)},$
$\widehat{{\rm FDP}}(t),$ $\hat t$ and $V_q$ in Algorithm \ref{alg1} with 
$\widetilde{T}_n^{(q)},$ $\widetilde{{\rm FDP}}(t),$ $\tilde t$ and $\widetilde{V}_q,$ respectively.

To impose additional regularity conditions similar to those in Theorem~\ref{thm:nonFDR}, we let $\widetilde{F}_q(\cdot)$ be the distribution of $\widetilde{T}_n^{(q)},$ and write
$\widetilde{\zeta}_q=\Phi^{-1}[\widetilde{F}_q\{\widetilde{T}_n^{(q)}\}]$, which follows the standard normal distribution. Accordingly, we define
\begin{align*}
	\widetilde{\mathcal{S}}_q(\tilde \gamma) = \{{q'}\in[Q]: {q'}\neq q, |{\rm Corr}(\widetilde{\zeta}_q,\, \widetilde{\zeta}_{{q'}})| \geq(\log Q)^{-2-\tilde \gamma}\}~\text{ for } \tilde \gamma>0~\text{and}~q \in [Q].
\end{align*}
Similar to Theorem~\ref{thm:nonFDR} established for fully observed functional data, the next theorem indicates that our proposed multiple testing procedure for error-contaminated functional data is guaranteed to control both FDR and FDP below the nominal level asymptotically.

\begin{theorem}\label{thm:errFDR}
Suppose that Conditions  {\rm \ref{cond:covfun}--\ref{cond:Xe}}  hold, $\max_{q,q'\in[Q]}|{\rm Corr}(\tilde{\zeta}_q,\tilde{\zeta}_{q'})|\leq c_{\tilde \zeta}$ for some constant $c_{\tilde \zeta}\in(0,1)$, and $\max_{q\in[Q]}|\widetilde{\mathcal{S}}_q(\tilde \gamma)|=o(p^{2 \tilde \nu})$ for some constants $\tilde \gamma>0$ and $0 < \tilde \nu<(1-c_{\tilde \zeta})/(1+c_{\tilde \zeta})$.  If 
$p\lesssim n^{\tilde\kappa}$ for some constant $\tilde\kappa>0$, then $\limsup_{n,Q\rightarrow\infty}{\rm FDR}({\tilde{t}})\leq \alpha Q_0/Q$ and $\lim_{n,Q\rightarrow\infty}\bbP\{{\rm FDP}({\tilde{t}})\leq \alpha Q_0/Q+\widetilde \varkappa\}=1$ for any $\widetilde \varkappa>0$.
\end{theorem}

\begin{remark}
\label{remark:general}
Note that our proposed multiple testing procedure considers $Q = p(p-1)/2$ hypotheses of cross-covariance functions induced from $p(p-1)/2$ pairs of $\{\widehat X_{ij}(\cdot), \widehat X_{ik}(\cdot)\}$ based on $p$ estimated curves $\widehat X_{i1}(\cdot), \dots, \widehat X_{ip}(\cdot).$ 
The theoretical guarantee of this framework can be easily extended to a more generalized multiple testing setting, where multiple hypotheses are induced from $Q=p(p-1)/2$ pairs based on $p(p-1)/2$ instead of $p$ estimated curves. For example, in Section~\ref{sec.fgm}, we transform the multiple testing problem of conditional dependence structures to that of the error cross-covariance patterns in (\ref{test.fgm}), which is induced from
$p(p-1)/2$ pairs of functional residuals as specified in (\ref{func.res}), i.e. $\{\hat \varepsilon_{i,jk}(\cdot), \hat\varepsilon_{i,kj}(\cdot)\}$ for $1 \leq j <k \leq p.$ 
Not only does our proposal for the generalized multiple testing problem remain unchanged, but also, following similar arguments as those in the proofs of Theorems~\ref{thm:nonFDR} and \ref{thm:errFDR}, it can still achieve effective control of false discoveries asymptotically. Further details are available in  the application to the functional graphical model in Section~\ref{sec.fgm}.
\end{remark}



\vspace{-0.7cm}
\section{Applications}
\label{sec:app}
\vspace{-0.3cm}
\subsection{Discretely observed functional data}
\label{sec.part}
For any $i \in [n]$ and $j \in [p],$ each $X_{ij}(\cdot)$ is not directly observable in practice. Instead, it is observed with errors at $T_{ij}$ random time points, $U_{ij1}, \dots, U_{ijT_{ij}} \in \cU.$ Let $W_{ijt}$ be the observed value of $X_{ij}(U_{ijt})$ satisfying
\begin{equation} \label{model.partial}
    W_{ijt}=X_{ij}(U_{ijt}) + \varsigma_{ijt} \,,
\end{equation}
where the errors $\varsigma_{ijt}$'s, independent of $X_{ij}(\cdot)$'s, are i.i.d. with $\mathbb{E}(\varsigma_{ijt})=0$ and ${\rm Var}(\varsigma_{ijt})=\sigma_j^2<\infty.$
The sampling frequency $T_{ij}$ plays a vital role in the estimation as it may affect the choice of the estimation procedure. For densely sampled functional data with $T_{ij}$'s larger than some order of $n,$ the conventional approach is to implement nonparametric smoothing to the observations from each subject, thus reconstructing each individual curve before subsequent analysis \cite[]{Zhang2007}. For sparsely sampled functional data with bounded $T_{ij}$'s, pre-smoothing is no longer viable; it is standard practice to consider pooling data from all subjects to borrow strength across all observations \cite[]{zhang2016}. 

Since brain signals in neuroimaging data of our interest are typically recorded at a dense set of points, we make use of the pre-smoothing to facilitate the presentation of methodology and theory. In what follows, denote $K_h(\cdot)=h^{-1}K(\cdot/h)$ for a univariate kernel with bandwidth $h>0.$ We first apply the frequently-adopted local linear smoothers to the observed data $\{(U_{ijt},W_{ijt}):t\in [T_{ij}]\}$ for each $i \in [n]$ and $j \in [p],$ and thus can obtain the estimated curves $\widehat{X}_{ij}(u)=\hat{g}_{0ij},$ where
$$
(\hat g_{0ij}, \hat g_{1ij}) = \arg \min_{g_{0ij}, g_{1ij}}\,\sum_{t=1}^{T_{ij}}\big\{W_{ijt}-g_{0ij}-g_{1ij}(U_{ijt}-u)\big\}^2 K_{h_{ij}}(U_{ijt}-u).
$$
For any $j \in [p],$ individual functions across $i \in [n]$ often admit similar smoothness properties and sometimes similar patterns, it is thus reasonable to use a common bandwidth $h_j$ for all of them.
We then employ the multiple testing procedure in Section \ref{sec:error} to the estimated curves $\{\widehat X_{ij}(\cdot)\},$ which leads to $e_{ij}(\cdot)= \widehat{X}_{ij}(\cdot)-X_{ij}(\cdot).$

To provide the theoretical guarantee for such procedure based on Theorem~\ref{thm:errFDR}, we need to verify Condition~\ref{cond:Xe}, i.e. specifying the convergence rates of $\max_{j,k\in[p]}\|\widehat\Sigma_{jk}^{X e}\|_{\cS}$ and $\max_{j,k\in[p]}\|\widehat{\Sigma}_{jk}^{e e}\|_{\cS}.$ 
To this end, we impose the following regularity conditions. 

\begin{condition}\label{cond.subgauss}
For each $i\in[n],$ $j\in[p], t \in [T_{ij}]$, $\varsigma_{ijt}$ is a sub-Gaussian random variable, i.e., there exists some constant $\tilde c> 0$ such that for all $ \tilde x \in \mathbb{R},$
$
    \mathbb{E}
    \{\exp(\varsigma_{ijt} \tilde x )\} \leq \exp(2^{-1}\tilde c^2 \sigma_{j}^2\tilde x^2).
$
\end{condition}

\begin{condition}\label{cond.Tij}
{\rm(i)} For any $j \in [p],$ the average sampling frequency $\widebar{T}_j=(n^{-1}\sum_{i=1}^n T_{ij}^{-1})^{-1} \rightarrow \infty$, $h_j\rightarrow0$ 
and $\widebar{T}_jh_j  \rightarrow \infty$  as $n\rightarrow\infty;$
{\rm(ii)} $h_j \asymp h$ and $\widebar{T}_j \asymp T$ for each $j$. 
\end{condition}

\begin{condition}\label{cond.fu}
{\rm (i)} Let $\{U_{ijt}; i\in[n], j\in[p], t\in[T_{ij}]\}$ be i.i.d. copies of a random variable $U$ defined on $\cU$ 
with density $f_U(\cdot)$ satisfying $0<m_f\leq \inf_{u\in\cU} f_U(u)\leq \sup_{u\in\cU}f_U(u)\leq M_f<\infty$;
{\rm (ii)} $X_{ij}(\cdot)$, $U_{ijt}$ and $\varsigma_{ijt}$ are independent for each $i,j,t.$
\end{condition}

\begin{condition}\label{cond.kern}
{\rm (i)}
The kernel $K(\cdot)$ is a symmetric probability density function on compact support $[-1,1]$ with $\int u^2K(u)\,{\rm d} u<\infty$ and $\int K^2(u)\,{\rm d}u<\infty$. 
{\rm (ii)} $K(\cdot)$ is Lipschitz continuous, i.e., there exists some positive constant $L$ such that $|K(u)-K(v)|\leq L|u-v|$ for any $u,v\in[-1,1]$.
\end{condition}

\begin{condition}\label{cond.cont}
{\rm(i)} $\partial^2 \mu_{j}(u)/\partial u^2$ is uniformly bounded over $u\in\cU$ and $j\in[p]$. {\rm(ii)} The covariance functions satisfy $\max_{j \in [p]} \sup_{u \in \cU} |\Sigma_{jj}(u,u)| = O(1)$.

\end{condition}



Condition~\ref{cond.subgauss} requires that the random errors are sub-Gaussian.
Conditions \ref{cond.Tij}--\ref{cond.cont} are standard in the literature of nonparametric smoothing for functional data \cite[]{Zhang2007,zhang2016} adaptable to the multivariate setting. The following proposition indicates that Condition~\ref{cond:Xe} holds for discretely observed functional data.

\begin{proposition}
\label{thm.partial}
Under Conditions~{\rm \ref{cond:subG}} and {\rm\ref{cond.subgauss}--\ref{cond.cont}}, we have
$\max_{j,k\in[p]}\|\widehat{\Sigma}_{jk}^{X e}\|_{\cS} = 
O_{\rm p} \big[\log(p \vee n)(Th)^{-1/2} + \{\log(p\vee n)\}^{1/2} h^2\big]$
and
$\max_{j,k\in[p]}\|\widehat{\Sigma}_{jk}^{e e}\|_{\cS} = 
O_{\rm p}\big\{\log(p \vee n)(Th)^{-1} + h^4\big\}$.
\end{proposition}


\vspace{-0.3cm}

\begin{remark}
Recall that $p\lesssim n^{\tilde\kappa}$ for some $\tilde\kappa>0$ is required in Theorem~\ref{thm:errFDR}.
To facilitate further discussion, we consider a simplified common scenario where there exist three constants $\eta_2,\eta_3$ and $\eta_4>0$ such that $ T\asymp n^{\eta_2}$, $h\asymp n^{-\eta_3}$ and  $p\asymp n^{\eta_4}$.
Some calculations show that Condition~\ref{cond:Xe} is satisfied by setting $a_1 = \min\{(\eta_2-\eta_3)/2, 2\eta_3\}$, $a_2=1,$ $a_3 = \min(\eta_2-\eta_3, 4\eta_3)$ and  $a_4=1$ provided that $\eta_3>1/4$ and $\eta_2>\eta_3+1.$ 
With verifiable Condition~\ref{cond:Xe} under the possibly sub-optimal rates established in Proposition \ref{thm.partial} and other conditions required in Theorem~\ref{thm:errFDR}, an application of Theorem~\ref{thm:errFDR} yields that our multiple testing procedure for discretely observed functional data ensures the effective control of both FDR and FDP asymptotically.
\end{remark}

\vspace{-0.3cm}
\subsection{Functional graphical model}
\label{sec.fgm}

Suppose that we observe $p$-vector of functional data $\bY_i(\cdot)=\{Y_{i1}(\cdot),\ldots,Y_{ip}(\cdot)\}^{\T}$ for $i \in [n]$ defined on $\cU$. Assume that  $\bY_i(\cdot)$ follows a mean zero multivariate Gaussian process with covariance function matrix $\bXi(u,v)=\cov\{\bY_i(u),\bY_i(v)\}$ for $(u,v) \in \cU^2.$ Let $\bY_{i,-j,-k}(\cdot)$ be the subvector of $\bY_i(\cdot)$ with $Y_{ij}(\cdot)$ and $Y_{ik}(\cdot)$ being removed. In this section, we formulate the estimation of functional graphical model depicting the conditional dependence structure among $p$ components of $\bY_i(\cdot)$ as simultaneous testing of $Q = p(p-1)/2$ hypotheses,
\begin{align}\label{eq:condH0}
\widetilde H_{0,jk}: Y_{ij}(\cdot) \perp \!\!\! \perp Y_{ik}(\cdot) \,|\, \bY_{i,-j,-k}(\cdot) ~~~\mbox{vs}~~~
\widetilde H_{1,jk}: Y_{ij}(\cdot) \not \perp \!\!\! \perp Y_{ik}(\cdot) \,|\, \bY_{i,-j,-k}(\cdot)
\end{align}
for $1 \leq j<k \leq p$, and we identify a significant edge between nodes $j$ and $k$ if and only if $\widetilde H_{0,jk}$ is rejected. 

For each $1 \leq j<k \leq p$, consider regressing each pair $Y_{ij}(\cdot)$ and $Y_{ik}(\cdot)$ on the remaining $p-2$ functional variables, i.e.,
\begin{align}
Y_{ij}(u)= \sum_{\ell \neq j,k} \int_{\cU} Y_{i\ell}(v) \beta_{jk,\ell}(u,v)\,{\rm d}v+\varepsilon_{i,jk}(u) \,, ~~ u \in \cU\, ,\label{eq_rg} 
\\
Y_{ik}(u)= \sum_{\ell \neq j,k} \int_{\cU} Y_{i\ell}(v) \beta_{kj,\ell}(u,v)\,{\rm d}v+\varepsilon_{i,kj}(u)\,,  ~~ u \in \cU\, ,\nonumber
\end{align}
where $\varepsilon_{i,jk}(\cdot)$ and $\varepsilon_{i,kj}(\cdot)$ represent the functional regression errors, independent of $\bY_{i,-j,-k}(\cdot)$. Under such functional regression setup, we have $\cov\big\{Y_{ij}(u), Y_{ik}(v) \big|\bY_{i,-j-k}(\cdot) \big\}$
\begin{eqnarray*}
 & = & {\mathbb E}\Big(\big[Y_{ij}(u) - {\mathbb E}\{Y_{ij}(u)|\bY_{i,-j-k}(\cdot)\}\big] \big[Y_{ik}(v) - {\mathbb E}\{Y_{ik}(v)|\bY_{i,-j-k}(\cdot)\}\big] \Big| \bY_{i,-j-k}(\cdot)\Big)\\
 & = & {\mathbb E}\big\{\varepsilon_{i,jk}(u) \varepsilon_{i,kj}(v) | \bY_{i,-j,-k}(\cdot)\big\},
\end{eqnarray*}
which possesses the property that
$\cov\big\{Y_{ij}(u), Y_{ik}(v) \big|\bY_{i,-j-k}(\cdot) \big\} = 0$ for any $(u,v) \in \cU^2,$ and hence 
$\Sigma_{jk}^{\varepsilon}(u,v) = {\cov}\{\varepsilon_{i,jk}(u), \varepsilon_{i,kj}(v)\}=0$
for any $(u,v) \in \cU^2$ whenever $Y_{ij}(\cdot) \perp \!\!\! \perp Y_{ik}(\cdot) \,|\, \bY_{i,-j,-k}(\cdot).$
Thus, the hypotheses of conditional independence in \eqref{eq:condH0} reduces to 
\begin{align}\label{test.fgm}
	\widetilde H_{0}^{(q)}: \Sigma_{j_qk_q}^{\varepsilon}(u,v)=0 ~\mbox{for any } (u,v)\in\cU^2\mbox{ and } (j_q,k_q)\in\cH
	~\mbox{vs}~
	\widetilde  H_{1}^{(q)}:\widetilde 
 H_{0}^{(q)} \mbox{ is not true}\,,
\end{align}
for $q \in [Q].$ 

Let $\widehat{\boldsymbol{\beta}}_{jk}(\cdot,\cdot)=\{\hat\beta_{jk,\ell}(\cdot,\cdot)\}_{\ell \neq j,k}^{\T}$ be some consistent estimator of $\boldsymbol{\beta}_{jk}(\cdot,\cdot)= \{\beta_{jk,\ell}(\cdot,\cdot)\}_{\ell \neq j,k}^{\T}$.
Define the functional residuals by
\begin{equation}
\label{func.res}
	\hat{\varepsilon}_{i,jk}(u) := Y_{ij}(u)-\int_{\cU}\bY_{i,-j,-k}^{\T}(v)\widehat{\boldsymbol{\beta}}_{jk}(u,v)\,{\rm d}v = {\varepsilon}_{i,jk}(u) + \text{deviation},
\end{equation}
which takes a similar form as (\ref{error.X}) within the error-contamination framework.
Let the sample covariance functions between residuals be
$\widetilde\Sigma_{jk}^{\varepsilon}(u,v)=(n-1)^{-1}\sum_{i=1}^n\{\hat{\varepsilon}_{i,jk}(u)-\bar{\hat\varepsilon}_{jk}(u)\} \{\hat{\varepsilon}_{i,kj}(v)-\bar{\hat\varepsilon}_{kj}(v)\} $ with $\bar{\hat\varepsilon}_{jk}(u)=n^{-1}\sum_{i=1}^n\hat\varepsilon_{i,jk}(u)$
for any $j,k\in[p]$. Building upon the methodology in Section~\ref{sec:error}, for each $q \in [Q]$, we  propose the test statistic for $\widetilde H_0^{(q)}$ as
\begin{align*}
	\widetilde{T}_n^{\varepsilon,(q)}= n\| \widetilde{\Sigma}_{j_qk_q}^{\varepsilon}\|_{\cS}^2\,.
\end{align*}
The above formulation then allows for the direct adaptation of the notation and the multiple testing procedure from Section~\ref{sec:error} to the functional graphical model setting by substituting the true signal pairs $\{X_{ij}(\cdot), X_{ik}(\cdot)\}$ with $\{\varepsilon_{i,jk}(\cdot), \varepsilon_{i,kj}(\cdot)\}$, and the estimated signal pairs $\{\widehat X_{ij}(\cdot), \widehat X_{ik}(\cdot)\}$ with $\{\hat \varepsilon_{i,jk}(\cdot), \hat \varepsilon_{i,kj}(\cdot)\}$ for any $1 \leq j<k \leq p$.  To save space, we annotate quantities related to the functional graphical model with a superscript $\varepsilon$. For notational consistency, we will continue to use FDP and FDR to denote the false discovery proportion and rate, respectively, where ${\rm FDP}(t)={R_0^{\varepsilon}(t)}/\{R^{\varepsilon}(t)\vee1\}$ in this example.

We then follow the proposed multiple testing procedure in Section \ref{sec:error} based on $\widetilde T_{n}^{\varepsilon,(q)}$. Similar to Section~\ref{sec.part}, we need to verify Condition~\ref{cond:Xe} by establishing the convergence rates in  \eqref{cond.cov} within the functional graphical model framework.
Denote  $\widehat{\bXi}(u,v)=(n-1)^{-1}\sum_{i=1}^n\{\bY_i(u)-\widebar\bY(u)\}\{\bY_i(v)-\widebar\bY(v)\}^{\T}$ and $\widebar\bY(u)=n^{-1}\sum_{i=1}^n\bY_i(u)$.
Let $\bXi_{-j,-k}(u,v)$ and $\widehat{\bXi}_{-j,-k}(u,v)$ denote the submatrices of $\bXi(u,v)$  and $\widehat{\bXi}(u,v)$ by removing the $(j,k)$th rows and $(j,k)$th columns.  Suppose that 
\begin{align} 
\max_{j,k\in[p]} \|\widehat{\boldsymbol{\beta}}_{jk}-\boldsymbol{\beta}_{jk}\|_{\cS,1}  =O_{\rm p}&(\delta_{1n}) \,, \label{a3n}
\\
\max_{j,k\in[p]} \Big\|\int_\cU \widehat{\bXi}_{-j,-k}(\cdot,s) (\widehat{\bbeta}_{jk}-\bbeta_{jk})(\cdot,s)\,&{\rm d}s\Big\|_{\cS,\max}  =O_{\rm p}(\delta_{2n}) \,.\label{a4n}
\end{align}
for some rates $\delta_{1n}$ and $\delta_{2n}$.
Condition~\ref{cond:Xe} then holds by fulfilling
\begin{equation}
\label{a.rate}
\delta_{1n}n^{-1/2} (\log p)^{1/2} \asymp n^{-a_{1}} (\log p)^{a_2} \text{   and   } \delta_{1n}\delta_{2n} \asymp n^{-a_{3}} (\log p)^{a_4}.
\end{equation}
See \eqref{a2.rate} and \eqref{a1.rate} in Section~\ref{supp_th5} of the Supplementary Material for technical details.
Therefore, to verify Condition~\ref{cond:Xe}, we only need  to specify the rates in \eqref{a3n} and \eqref{a4n}. 
To achieve this, we will introduce an estimation procedure to obtain $\widehat \bbeta_{jk}(\cdot,\cdot)$ and
then derive the corresponding rates in \eqref{a3n} and \eqref{a4n} respectively in Propositions~\ref{lm_beta} and \ref{lm_betasigma} below. 

Recall that our framework replaces $\{X_{ij}(\cdot), X_{ik}(\cdot)\}$ and $\{\widehat X_{ij}(\cdot), \widehat X_{ik}(\cdot)\}$ in Section~\ref{sec:error} by $\{\varepsilon_{i,jk}(\cdot), \varepsilon_{i,kj}(\cdot)\}$ and $\{\hat \varepsilon_{i,jk}(\cdot), \hat \varepsilon_{i,kj}(\cdot)\}$, respectively, for any $1 \leq j<k \leq p$.
Applying a  generalized version of Theorem~\ref{thm:errFDR} as discussed in Remark~\ref{remark:general}, we establish the
next theorem, which guarantees that our proposed multiple testing procedure for functional graphical model can effectively control false discoveries asymptotically.

\begin{theorem}\label{thm:noncondH0}
 Suppose that Conditions {\rm \ref{cond:covfun}--\ref{con_eigenvalue}} hold for $\{\varepsilon_{i,jk}(\cdot)\}_{j,k \in [p]},$ $\max_{q,q'\in[Q]}|{\rm Corr}(\tilde{\zeta}_q ^{\varepsilon},\tilde{\zeta}_{q'}^{\varepsilon})|\leq c_{\tilde \zeta ^\varepsilon}$ for some constant $c_{\tilde \zeta ^\varepsilon}\in(0,1)$, and $\max_{q\in[Q]}|\tilde{\mathcal{S}}_q^{\varepsilon}(\tilde \gamma^{\varepsilon})|=o(p^{2\tilde \nu^\varepsilon})$ for some constants $\tilde \gamma^{\varepsilon}>0$ and $0<\tilde \nu^\varepsilon<(1-c_{\tilde \zeta ^\varepsilon})/(1+c_{\tilde \zeta ^\varepsilon})$. 
 If $p\lesssim n^{\tilde \kappa^\varepsilon}$ for some constant ${\tilde \kappa^\varepsilon}>0$ and {\rm (\ref{a.rate})} is satisfied with $a_1, a_3 >1/2$, then it holds that $\limsup_{n,Q\rightarrow\infty}{\rm FDR}({\tilde{t}}^{\varepsilon})\leq \alpha Q_0^\varepsilon/Q$ and $\lim_{n,Q\rightarrow\infty}\bbP\{{\rm FDP}({\tilde{t}}^{\varepsilon})\leq \alpha Q_0^\varepsilon/Q+\widetilde \varkappa^\varepsilon\}=1$ for any $\widetilde\varkappa^\varepsilon>0$.
\end{theorem}

We next develop a three-step procedure to estimate the functional coefficients ${\boldsymbol{\beta}}_{jk}(\cdot,\cdot)$ for $1 \leq j \neq k \leq p$. To ensure a feasible solution under a high-dimensional regime, we impose a sparsity assumption on $\{\bbeta_{jk}(\cdot,\cdot)\}$, i.e., $\bbeta_{jk}(\cdot,\cdot)$ is functional sparse with support set $S_{jk} = \big\{ \ell \in [p] \setminus \{j,k\}:\|\beta_{jk,\ell}\|_\mathcal{S}\neq 0\big\}$ and its cardinality $|S_{jk}| = s_{jk} \ll p.$ 
Let $s = \max_{jk} s_{jk}$.

In the first step, due to the infinite-dimensionality of functional data, we approximate each $Y_{ij}(\cdot)$ for $j \in [p]$ under the Karhunen-Lo\`eve expansion truncated at $d_j$ (to be specified in Section~\ref{sec:sim}), serving as the foundation of functional principal components analysis (FPCA),
$Y_{ij}(\cdot) \approx \sum_{m=1}^{d_j}\xi_{ij  m}\phi_{j  m}(\cdot) = \bxi_{ij}^{\T}\bphi_j(\cdot),$
where $\bxi_{ij} = (\xi_{ij1},\dots,\xi_{ijd_j})^{\T},$ $\bphi_j(\cdot) = \{\phi_{j1}(\cdot),\dots,\phi_{jd_j}(\cdot)\}^{\T}.$ Here $\xi_{ij  m} = \langle Y_{ij }, \phi_{j  m} \rangle,$ namely FPC scores, correspond to a sequence of random variables with $\cov(\xi_{ij  m}, \xi_{ij  m'}) = \omega_{ j  m}I(m=m'),$ where $\omega_{j1} \geq \omega_{j2} \geq \dots >0$ are the eigenvalues of $\Xi_{jj}(\cdot,\cdot),$ the $(j,j)$th entry of $\bXi(\cdot,\cdot)$, and $\phi_{j1}(\cdot), \phi_{j2}(\cdot), \dots$ are the corresponding eigenfunctions. To implement FPCA based on $n$ observations, we perform eigenanalysis of $\widehat\Xi_{jj}(\cdot,\cdot),$ the $(j,j)$th entry of $\widehat\bXi(\cdot,\cdot),$ and obtain estimated eigenvalue/eigenvector pairs $\{\hat\omega_{jm}, \hat\phi_{jm}(\cdot)\}_{m \in [d_j]}.$ The estimated FPC scores are $\hat\xi_{ij  m} = \langle Y_{ij }, \hat\phi_{j  m}\rangle$ for $m \in [d_j].$

In the second step, some calculations lead to the representation of \eqref{eq_rg} as
\begin{align} \nonumber
   \bxi_{ij}^{\T} = \sum_{\ell \neq j,k} \bxi_{i\ell}^{\T} \bPsi_{jk,\ell} + \br_{i,jk}^{\T} + \beps_{i,jk}^{\T},
\end{align}
where 
$\bPsi_{jk,\ell}=\int_{\cU}\int_{\cU}\bphi_{\ell}(v)\beta_{jk,\ell}(u,v)\bphi_{j}(u)^{\T}\, {\rm d}u  {\rm d}v\in\mathbb{R}^{d_{\ell}\times d_j}$ for $\ell \in [p]\setminus\{j,k\},$
$\br_{i,jk}$ and $\beps_{i,jk}$
are formed by truncation and random errors, respectively. See \eqref{supp_matrix} of the Supplementary Material for detailed expressions of $\br_{i,jk}$ and $\beps_{i,jk}$. 
Notice that 
$\bPsi_{jk,\ell} = \boldsymbol{0}$ if $\beta_{jk,\ell}(\cdot,\cdot) = 0,$ we can rely on the block sparsity in $\{\bPsi_{jk,\ell}\}$ to recover the functional sparsity in $\{\beta_{jk,\ell}(\cdot,\cdot)\}.$ Towards this, we let $\widehat \bxi_{ij} \!= \!(\hat \xi_{ij1},\dots,\hat \xi_{ijd_j})^{\T}$ and $\widehat \bphi_j(\cdot) \!=\! \{\hat \phi_{j1}(\cdot),\dots,\hat \phi_{jd_j}(\cdot)\}^{\T},$ and then propose to solve a standardized group lasso \cite[]{simon2012} problem by minimizing the following penalized least squares criterion over $\{\bPsi_{jk,\ell}\}_{\ell \in [p] \setminus \{j,k\}}$:
\begin{equation}
\label{target}
 \frac{1}{2} \sum_{i=1}^n\|\widehat \bxi_{ij}^{\T}- \sum_{\ell \neq j,k}  \widehat \bxi_{i\ell}^{\T} \bPsi_{jk,\ell}\|^2 + \tau_{n}\sum_{\ell \neq j,k}\|\widehat \bV_\ell \bPsi_{jk,\ell}\|_\tF,
\end{equation} 
where $\widehat \bV_{\ell}= (\widehat \bxi_{1\ell},\dots,\widehat \bxi_{n\ell})^{\T}\in \mathbb{R}^{n\times d_\ell}$ 
for $\ell \in [p]\setminus\{j,k\}$
and $\tau_n \geq 0$ is the regularization parameter.
Let $\{\widehat \bPsi_{jk,\ell}\}_{\ell \in [p] \setminus \{j,k\}}$  be the minimizer of (\ref{target}).

In the third step, we estimate functional coefficients by $$\hat{\beta}_{jk,\ell}(u,v) = \widehat{\bphi}_\ell(v)^{\T}\widehat{\bPsi}_{jk,\ell}\widehat{\bphi}_{j}(u)~~\text{for}~(u,v) \in \cU^2~\text{and}~\ell \in [p] \setminus\{j,k\}.$$


Before presenting the theoretical results, we list some regularity conditions.  To
simplify notation, we assume the same $d = d_j$ across $j \in [p]$, but our theoretical results below can be extended naturally to the more general setting where $d_j$’s are different. 

\begin{condition}\label{con_eigen}
For each $j  \in[p]$,  $\omega_{j  1} > \omega_{j  2}> \cdots >0$. There exists some constant  $\varpi > 1$ such that $\omega_{j  m}- \omega_{j (m+1)} \gtrsim m^{-\varpi-1}$ and $\omega_{jm} \lesssim  m^{-\varpi}$ for $m \in\mathbb{N}_+$.
\end{condition}

\begin{condition} \label{con_fof_eigen_min}
Denote the diagonal function matrix by $\widetilde \bD_0 = \text{diag}(\Xi_{11},\dots,\Xi_{pp}).$ The infimum 
$
\underline{\mu} = \inf\limits_{\bPhi \in \Bar{\mathbb{H}}_0^{p}} \frac{\iint\bPhi(u)^{\T}\bXi(u, v)\bPhi(v) \,{\rm d}u{\rm d}v}{\iint\bPhi(u)^{\T}\widetilde{\bD}_0(u, v)\bPhi(v) \,{\rm d}u{\rm d}v}\gtrsim s d^{\varpi +3} \{{\log(pd)}/{n}\}^{1/2},
$
where $\bPhi \in \Bar{\mathbb{H}}_0^{p} = \{\bPhi \in \mathbb{H}^{p}: \iint\bPhi(u)^{\T}\widetilde{\bD}_0(u, v)\bPhi(v) \,{\rm d}u{\rm d}v \in (0,\infty)\}.$
\end{condition}

\begin{condition} \label{con_fof_beta_a}
For each $j,k \in[p]$ and $\ell\in S_{jk},$ $\beta_{jk,\ell}(u,v) = \sum_{m_1,m_2=1}^\infty b_{jk,\ell m_1m_2} \phi_{jm_1}(u) \phi_{\ell m_2}(v) $ and there exists some constant $\upsilon > \varpi  /2 + 1$ such that $|b_{jk,\ell m_1m_2}| \lesssim (m_1+m_2)^{-\upsilon-1/2}$ for $m_1,m_2 \geq 1.$
\end{condition}

Conditions \ref{con_eigen} and \ref{con_fof_beta_a} are standard in functional linear regression literature \cite[]{kong2016,guo2019} with parameter $\varpi$ capturing the decay rate of eigenvalues and parameter $\upsilon$ controlling the level of smoothness in nonzero coefficient functions. Condition~\ref{con_fof_eigen_min} can be understood as requiring a lower bound on the minimum eigenvalue of the correlation function matrix of $\bY_i(\cdot).$
The following propositions verify the rates in \eqref{a3n} and \eqref{a4n}, respectively, for $\widehat \bbeta_{jk}(\cdot,\cdot)$ obtained through the proposed estimation procedure. This implies that Condition~\ref{cond:Xe} holds for the functional graphical model. 

\begin{proposition}[verify the rate in \eqref{a3n}]
\label{lm_beta}
    Suppose that Conditions~{\rm\ref{cond:subG}} and {\rm\ref{con_eigen}}--{\rm\ref{con_fof_beta_a}} hold. If $n \gtrsim \log (pd)d^{4\varpi+4} ,$
    then there exist some positive constants $c_1$ and $c_2$ such that, for any regularization parameter, $\tau_n \gtrsim s  [  d^{\varpi+2}\{{\log(pd)}/{n}\}^{1/2}+ d^{-\upsilon+1}]$, 
the estimate $\widehat\bbeta_{jk}(\cdot,\cdot)$ satisfies
\begin{equation}  \nonumber
\begin{split}
\max_{j,k\in[p]} \|\widehat{\bbeta}_{jk}-\bbeta_{jk}\|_{\cS,1} \lesssim  \underline{\mu}^{-1}{d^{\varpi/2} s \tau_n} \asymp \delta_{1n},
\end{split}
\end{equation}
with probability greater than $1-c_1(pd)^{-c_2}.$
\end{proposition}

\begin{proposition}[verify the rate in \eqref{a4n}]\label{lm_betasigma}
    Suppose that the conditions in Proposition~{\rm\ref{lm_beta}} hold, then with probability greater than $1-c_1(pd)^{-c_2},$
$$\max_{j,k\in[p]} \Big\|\int_\cU \widehat{\bXi}_{-j,-k}(\cdot,s) (\widehat{\bbeta}_{jk}-\bbeta_{jk})(\cdot,s)\,{\rm d}s\Big\|_{\cS,\max} \lesssim \tau_n + \underline{\mu}^{-1} d^{-\varpi/2+1}s\tau_n \asymp \delta_{2n}.$$
\end{proposition}

\begin{remark}
Our procedure can be extended naturally to handle the functional graphical model with discrete observations by applying the methodology developed in Section~\ref{sec.fgm} to reconstructed curves obtained through the pre-smoothing step introduced in Section~\ref{sec.part}. To provide theoretical support for our multiple testing procedure, we still need to verify Condition~\ref{cond:Xe} by establishing the corresponding rates in (\ref{a3n}) and (\ref{a4n}). However, the initial pre-smoothing step results in additional estimation errors, which makes the subsequent verification much more involved. Empirically, we demonstrate the effectiveness of such procedure in controlling false discoveries through extensive simulations in Section~\ref{sim.fgm}.
\end{remark}

\vspace{-0.6cm}
\section{Simulation studies}
\label{sec:sim}
\vspace{-0.2cm}
We conduct a number of simulations to illustrate the proposed multiple testing procedures for the functional covariance model and functional graphical model in Sections~\ref{sim.fcm} and \ref{sim.fgm}, respectively, for both fully and discretely observed functional data.

\vspace{-0.4cm}
\subsection{Functional covariance model} 
\label{sim.fcm}
\vspace{-0.2cm}
In each scenario, 
we generate functional variables by $X_{ij}(u)=\bs(u)^{\T}\btheta_{ij}$ for $i \in [n], j \in [p]$ and $u \in \cU=[0,1],$ where $\bs(u)$ is a $10$-dimensional Fourier basis function and $\btheta_{i}=(\btheta_{i1}^{\T},\dots,\btheta_{ip}^{\T} )^{\T} \in \eR^{10p}$ is generated from a mean zero multivariate Gaussian distribution with block covariance matrix
$\bLambda \in \eR^{10p \times 10p},$ whose $(j,k)$th block is $\bLambda_{jk} \in \eR^{10 \times 10}$ for $j,k \in [p].$
The functional sparsity pattern in
$\bSigma(\cdot,\cdot)=\{\Sigma_{jk}(\cdot,\cdot)\}_{j,k\in[p]}$ with its $(j,k)$th entry $\Sigma_{jk}(u,v) = \bs(u)^{\T}\bLambda_{jk}\bs(v)$ can thus be characterized by the block sparsity structure in $\bLambda.$ Define $\bLambda_{jk} = \Pi_{jk} \bDelta$ for $j, k \in [p]$ with $\bDelta = \text{diag}(1^{-2},\dots, 10^{-2}),$ which results in polynomially decaying eigenvalues of $\Sigma_{jj}(\cdot,\cdot)$ for each $j.$

\begin{figure}[tbp]
	\centering
	\includegraphics[width=14cm, height = 6.7cm]{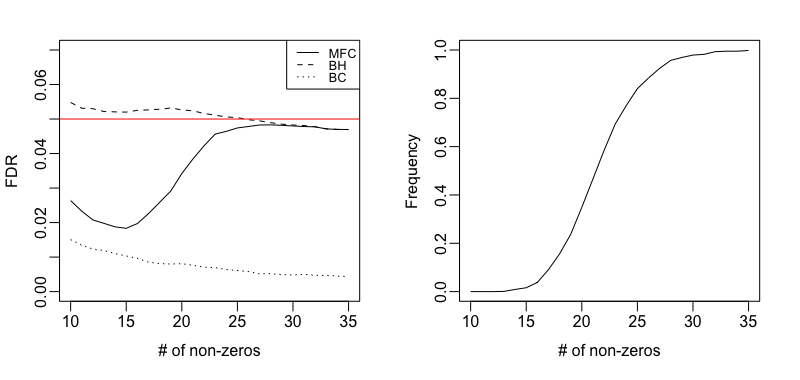}
 \vspace{-0.5cm}
	\caption{\label{plot_BHBC}{Empirical FDRs of MFC, BH and BC (left)  and the frequency that $\hat t$ exists (right) over 1000 simulation runs.}}
\end{figure}

We start by comparing the performance of our proposed {\bf m}ultiple testing procedure for the {\bf f}unctional {\bf c}ovariance model (MFC) with  {\bf B}enjamini--{\bf H}ochberg procedure (BH) \citep{Benjamini1995} and {\bf B}onferroni {\bf c}orrection procedure (BC) under sparse covariance settings. To this end, we generate $\bPi \equiv (\Pi_{jk})_{p \times p} = \bA\odot\bB + c' \bI_{p},$ where $\odot$ denotes the Hadamard product, $\bA = (A_{jk})_{p \times p}$ is a symmetric (0,1) matrix with 1 on the diagonal,    $\bB \in \mathbb{R}^{p\times p}$ is a symmetric matrix with entries being i.i.d. sampled from $\text{Unif}[0.2,0.3]$,  and   $c' = \big(-\lambda_{\min}(\bA\odot\bB)\big)_{+} + 0.01$ to guarantee the positive definiteness of $\bPi,$ with $\lambda_{\min}(\bM)$ denoting the minimum eigenvalue of matrix $\bM.$ Let the set $S_A = \{(j,k) : 1 \leq j <k \leq p, A_{jk} \neq 0 \}$ with its degree $|S_A|=s_A$. We consider $n = 200$, $p = 30$ and examine the empirical FDR by
\begin{align} \label{eq.eFDR}
	\frac{1}{B}\sum_{b=1}^B \frac{\sum_{q\in\cH_0}I\{V_q^{(b)}\geq \hat t^{(b)}\}}{1\vee \sum_{q=1}^Q I\{V_q^{(b)}\geq \hat t^{(b)}\}}\,
\end{align}
over $B = 1000$ replications, where $V_q^{(b)}$ and $\hat t^{(b)}$ are specified in Algorithm~\ref{alg1} for the $b$th replication. Figure~\ref{plot_BHBC} presents the empirical FDRs of the three competing procedures and the frequency that $\hat t$ exists for the proposed method at the $5\%$ nominal level when  $s_A \in\{10,11,\dots,35\}$. It is evident that BH tends to result in inflated FDRs and thus cannot control the FDP if the number of true alternatives is fixed (very sparse case). Similar patterns can also be found in \cite{Liu2014,Xia2018b}. 

We next assess the performance of the proposed MFC, BC and two commonly-used thresholding strategies, i.e., hard and soft functional thresholding methods, under both fully and discretely observed functional scenarios. The $(j,k)$th entries of hard and soft functional thresholding estimators are respectively defined as 
$\widehat\Sigma_{jk} I (\|\widehat\Sigma_{jk}\|_\cS \ge \tilde\tau)$ and $ \widehat\Sigma_{jk}(1-\tilde\tau/\|\widehat\Sigma_{jk}\|_{\cS})_{+},$ where $\widehat\Sigma_{jk}$ is specified in \eqref{eq:sigmajk} and the thresholding parameter $\tilde \tau > 0$ is selected by the cross-validation approach \citep{fang2023}.

For each method, we first generate fully observed curves $X_{ij}(\cdot)$ under different functional sparsity patterns in $\bSigma(\cdot,\cdot)$, i.e., block sparsity structures in $\bPi$ as follows. 
\begin{itemize}
 \item Model~1 (block banded). For $j,k \in [p],$ we generate $\Pi_{jk} = (1-{|j-k|}/{3})_{+}$.

\item Model~2 (block sparse without any special structure).  We generate $\bPi = \bB + c'' \bI_{p},$ where entries of $\bB$ are sampled independently from $\text{Unif}[0.3,0.8]$ with probability $3/p$ or $0$ with probability $1-3/p,$ and $c'' = \big(-\lambda_{\min}(\bB)\big)_{+} + 0.01.$ 
\end{itemize}
We then generate the discretely observed curves satisfying (\ref{model.partial}), where the observational time points $U_{ijt}$ and errors $\varsigma_{ijt}$ are sampled independently from $\text{Unif}[0,1]$ and ${\cal N}(0, 1),$ respectively.
We consider settings of $T_{ij}\in\{25,51\}$ and ultilize the Gaussian kernel with the optimal bandwidths proportional to $T_{ij}^{-1/5},$ as suggested in \cite{Zhang2007}, when implementing the local linear smoothers to obtain estimated curves $\widehat X_{ij}(\cdot)$. Let $\cH_1 = \cH \setminus \cH_0$ be the set of  true alternatives. Define the empirical power as 
\begin{align} \label{eq.epower}
	\frac{1}{100}\sum_{b=1}^{100} \frac{1}{|\cH_1|} \sum_{q\in\cH_1} I\{V_{q}^{(b)}\geq\hat t^{(b)}\} \,.
\end{align}

We finally report the empirical FDRs and the empirical powers over $100$ replications under both fully and discretely observed functional scenarios with $n \in \{100,200\}$ and $p \in \{30,60\}$ at $\alpha = 0.1$ in Table~\ref{table.sim.12} and at $\alpha = 0.05$ in Table~\ref{table.sim.12.0.05} of the Supplementary Material. Note that the empirical FDRs and powers for the discretely observed cases are calculated with $V_q^{(b)}$ and $\hat t^{(b)}$ in  \eqref{eq.eFDR} and \eqref{eq.epower}, respectively, replaced by $\widetilde{V}_q^{(b)}$ and  $\tilde t^{(b)},$ as specified in Section~\ref{sec:error}.
Several conclusions can be drawn from Tables~\ref{table.sim.12} and \ref{table.sim.12.0.05}. 
First, the empirical FDRs of the proposed MFC are well maintained below the target FDR levels for both fully and discretely observed functional data. In contrast, BC tends to be conservative with empirical FDRs substantially lower than the target levels. However, the soft thresholding method deteriorates significantly with highly elevated FDRs. Interestingly, the hard thresholding method seems to work well under these two models, with lower FDRs and powers. 
Second, MFC consistently achieves the highest empirical powers among four competitors across all settings. As expected, we observe enhanced powers as the number of subjects $n$ increases or as the dimension $p$ declines.
Third, all methods exhibit higher powers with increasing sampling frequency $T_{ij}$ and demonstrate comparable performance when curves are densely observed, compared to the fully observed functional case. 

\subsection{Functional graphical model}
\label{sim.fgm}

We now turn to the functional graphical model estimation. Different from the data generating process in Section~\ref{sim.fcm}, we first generate the functional components
$\tilde \varepsilon_{ij}(u) = \bs(u)^{\T}\widetilde \btheta_{ij}$ for $i\in [n]$ and $j=[p],$ 
where each $\widetilde \btheta_{ij}$ is sampled independently from ${\cal N}({\bf 0}, \bDelta),$ and then sequentially generate fully observed curves $Y_{i1}(\cdot), \dots, Y_{ip}(\cdot)$. To be specific, we begin by establishing directed acyclic graphs according to the following functional structural models and then moralize the directed graph to obtain the undirected  graph \cite[]{Cowell2007}. 
\begin{itemize}
    \item Model 3 (banded sparse). We generate $Y_{i1}(u)=\tilde \varepsilon_{i1}(u),$ $Y_{i2}(u)=\tilde \varepsilon_{i2}(u)$ and
$$Y_{ij}(u) = \sum_{k=1}^2\int_{\cU} Y_{i(j-k)}(v) \beta_{j(j-k)}(u,v) dv  + \tilde \varepsilon_{ij}(u) ~~\text{for}~~j\in [p] \setminus [2].$$

\item Model 4 (randomly sparse). 
For $j \in [p/3],$ we generate $Y_{ij}(u)=\tilde \varepsilon_{ij}(u).$  
To determine the directed edge set $E_D,$ we then randomly select one or two directed edges with equal probability from the candidate directed edge set $\{(k,j):1 \leq k<j\}$ for each node $j\in  \{p/3+1, \dots, p\}$ in a sequential way.
Then
$$Y_{ij}(u) =\sum_{(k,j) \in E_{D}} \int_{\cU} Y_{ik}(v) \beta_{jk}(u,v) dv   + \tilde \varepsilon_{ij}(u) ~~\text{for}~~j\in [p]\setminus [p/3],$$ which leads to an overall edge connection probability of approximately $3/p$ in the undirected graph.
\end{itemize}
We generate functional coefficients $\beta_{jk}(u,v) = \bs(u)^{\T} \bB_{jk} \bs(v) $ for $(u,v) \in \cU^2$ with $\bB_{jk} = (B_{jklm})_{l,m \in[10]}$ for $(k,j) \in E_{D}$ and $B_{jklm} = (-1)^{l+m} s_{D,j}^{-1} c_B (l+m)^{-2},$ which results in polynomially decaying basis coefficients. 
Here $c_B$'s are sampled independently from $\text{Unif}[4,6],$ and the normalization term is $s_{D,j} = |E_{D,j}|,$ where $E_{D,j} = \{ (k,j) \in E_{D}, 1 \leq k<j\}$ is defined as the directed edge set for node $j.$
We finally generate discretely observed curves according to (\ref{model.partial}) and the same procedure adopted in Section~\ref{sim.fcm} with sampling frequency of $T_{ij}$ in $\{25,51\}$.

We consider $n=200,400$ subjects of  $p=30,60$ functional variables, and repeat each simulation 100 times. To select the truncated dimension $d_j$ for each $j \in [p]$, we take the standard approach by selecting the largest $d_j$ eigenvalues
of $\widehat \Xi_{jj}(\cdot, \cdot)$ such that the cumulative percentage of selected eigenvalues exceeds $95\%$. Following the proposal of \cite{Xia2018b}, we select the optimal regularization parameter $\hat \tau_n$ with the principle of making $\sum_{q\in\mathcal{H}_0^\varepsilon}I(\widetilde{V}_q^\varepsilon \geq t)$ and $Q\{1-\Phi(t)\}$ as close as possible. Specifically, for a sequence of $\tau_n$ values, we  obtain estimated functional coefficients $\{\widehat \bbeta_{jk, \tau_n}(\cdot,\cdot)\}_{1 \leq j \neq k \leq p}$, construct the corresponding transformed statistics $\widetilde V_{q,\tau_n}^\varepsilon$ for each $\tau_n$ and choose 
$\hat \tau_n$ as the minimizer of
$$
\sum_{l = 1}^{10}\bigg\{\frac{\sum_{q\in\mathcal{H}} I(\widetilde V_{q,\tau_n}^\varepsilon \geq \Phi^{-1} [1-l\{1-\Phi(\sqrt{\log p})\}/10])}{Q \cdot l\{1-\Phi(\sqrt{\log p})\}/10}-1\bigg\}^2.
$$

We examine the performance of our proposed {\bf m}ultiple testing procedure for the {\bf f}unctional {\bf g}raphical model (MFG) based on empirical FDR and power, and compare it with the performance of BC, hard and soft functional thresholding procedures. The numerical results at the  $10\%$ and $5\%$ nominal levels are summarized in Table~\ref{table.sim.34} and Table~\ref{table.sim.34.0.05} of the Supplementary Material, respectively. Similar conclusions can be drawn compared
to the numerical summaries obtained for the functional covariance model in Section~\ref{sim.fcm}. Furthermore, we observe that the hard thresholding approach fails completely under both model settings. One possible explanation is that the hard and soft thresholding methods adopt a uniform threshold level across all pairs $(j,k)$ and thus are not able to handle the heteroscedastic problem of the large-scale cross-covariance function estimation, i.e. $\{\Sigma_{jk}^\varepsilon(u,v)\}_{1 \leq j <k \leq p}$.

\afterpage{%
    \clearpage
    \begin{landscape}
\begin{table}[tbp]
	\caption{Empirical FDRs  ($\%$) and powers ($\%$) of MFC, BC, hard and soft functional thresholding procedures at the $10\%$ nominal level for Models 1 and 2 over 100 simulation runs. \label{table.sim.12}}
	\begin{center}
		\vspace{-0.5cm}
		\resizebox{7.8in}{!}{
\begin{tabular}{cccccccccccccccccccc}
\hline
\multirow{3}{*}{$n$} & \multirow{3}{*}{$p$} & \multirow{3}{*}{Scenario} & \multicolumn{8}{c}{Model 1}                                                                           &  & \multicolumn{8}{c}{Model 2}                                                                           \\ \cline{4-11} \cline{13-20}
                     &                      &                           & \multicolumn{2}{c}{MFC} & \multicolumn{2}{c}{BC} & \multicolumn{2}{c}{Hard} & \multicolumn{2}{c}{Soft} &  & \multicolumn{2}{c}{MFC} & \multicolumn{2}{c}{BC} & \multicolumn{2}{c}{Hard} & \multicolumn{2}{c}{Soft} \\
                     &                      &                           & FDR      & Power     & FDR      & Power     & FDR       & Power      & FDR       & Power      &  & FDR      & Power     & FDR      & Power     & FDR       & Power      & FDR       & Power      \\ \hline
100                  & 30                   & Fully                     & 7.54      & 91.91      & 0.26      & 67.58      & 2.72       & 82.46       & 66.42      & 99.7        &  & 8.07      & 88.18      & 0.37      & 64.66      & 2.65       & 80.07       & 63.02      & 98.55       \\
                     &                      & $T_{ij}$ = 51             & 8.23      & 90.04      & 0.29      & 65.61      & 3.18       & 81.42       & 61.92      & 99.23       &  & 8.25      & 86.14      & 0.36      & 60.55      & 3.67       & 79.41       & 65.07      & 98.09       \\
                     &                      & $T_{ij}$ = 25             & 8.63      & 88.09      & 0.25      & 63.77      & 3.27       & 79.39       & 61.62      & 98.84       &  & 8.19      & 84.84      & 0.29      & 57.73      & 4.34       & 79.02       & 65.04      & 97.89       \\
                     & 60                   & Fully                     & 8.96      & 87.74      & 0.12      & 60.93      & 1.75       & 75.21       & 63.48      & 98.65       &  & 8.44      & 81.94      & 0.15      & 53.48      & 2.28       & 71.61       & 65.51      & 96.01       \\
                     &                      & $T_{ij}$ = 51             & 9.16      & 85.06      & 0.17      & 58.60       & 2.46       & 74.14       & 64.75      & 98.14       &  & 8.85      & 79.92      & 0.14      & 50.42      & 3.35       & 70.84       & 67.38      & 95.80        \\
                     &                      & $T_{ij}$ = 25             & 9.37      & 82.61      & 0.15      & 56.63      & 2.31       & 71.43       & 64.20       & 97.16       &  & 9.18      & 78.21      & 0.15      & 46.81      & 3.57       & 70.01       & 68.14      & 94.93       \\ \hline
200                  & 30                   & Fully                     & 8.63      & 99.86      & 0.14      & 95.49      & 1.14       & 97.95       & 61.98      & 100         &  & 9.08      & 98.89      & 0.20       & 92.05      & 0.96       & 94.36       & 64.02      & 99.95       \\
                     &                      & $T_{ij}$ = 51             & 8.49      & 99.68      & 0.17      & 93.28      & 1.72       & 97.44       & 57.68      & 99.98       &  & 9.47      & 98.39      & 0.27      & 90.43      & 1.54       & 93.68       & 61.93      & 99.89       \\
                     &                      & $T_{ij}$ = 25             & 8.50       & 99.40       & 0.25      & 90.47      & 1.70        & 96.53       & 58.20       & 100         &  & 8.65      & 98.07      & 0.18      & 88.66      & 1.64       & 93.07       & 63.01      & 99.86       \\
                     & 60                   & Fully                     & 9.27      & 99.66      & 0.08      & 90.37      & 1.41       & 96.54       & 75.94      & 100         &  & 9.38      & 96.67      & 0.18      & 83.4       & 1.16       & 90.09       & 63.68      & 99.57       \\
                     &                      & $T_{ij}$ = 51             & 9.12      & 99.33      & 0.11      & 86.91      & 1.86       & 95.45       & 74.12      & 99.99       &  & 9.58      & 95.92      & 0.16      & 81.11      & 1.51       & 89.21       & 67.19      & 99.40        \\
                     &                      & $T_{ij}$ = 25             & 9.10       & 98.84      & 0.09      & 83.68      & 2.41       & 94.31       & 73.51      & 99.99       &  & 9.77      & 95.03      & 0.14      & 78.68      & 1.66       & 87.99       & 67.04      & 99.23       \\ \hline

\end{tabular}
		}	
	\end{center}
	 \vspace{-0.5cm}
\end{table}
    \end{landscape}
    \clearpage
}

\afterpage{%
    \clearpage
    \begin{landscape}
\begin{table}[tbp]
 \caption{Empirical sizes  ($\%$) and empirical FDRs ($\%$)  of MFG, BC, hard and soft functional thresholding procedures at the $10\%$ nominal level for Models 3 and 4 over 100 simulation runs. \label{table.sim.34}}
	\begin{center}
	\vspace{-0.2cm}
		\resizebox{8in}{!}{
\begin{tabular}{ccccccccccclcccccccc}
\hline
\multirow{3}{*}{$n$} & \multirow{3}{*}{$p$} & \multirow{3}{*}{Scenario} & \multicolumn{8}{c}{Model 3}                                                                           &  & \multicolumn{8}{c}{Model   4}                                                                         \\ \cline{4-11} \cline{13-20}
                     &                      &                           & \multicolumn{2}{c}{MFG} & \multicolumn{2}{c}{BC} & \multicolumn{2}{c}{Hard} & \multicolumn{2}{c}{Soft} &  & \multicolumn{2}{c}{MFG} & \multicolumn{2}{c}{BC} & \multicolumn{2}{c}{Hard} & \multicolumn{2}{c}{Soft} \\
                     &                      &                           & FDR      & Power     & FDR      & Power     & FDR       & Power      & FDR       & Power      &  & FDR      & Power     & FDR      & Power     & FDR       & Power      & FDR       & Power      \\ \hline
200& 30                   & Fully                     & 8.46& 70.44& 0.45& 30.84& 12.53& 7.72& 83.11& 92.61
&  & 8.87      & 76.89      & 0.40       & 50.42      & 0.70        & 15.42       & 88.52      & 88.39       \\
                     &                      & $T_{ij}$ = 51             & 8.52& 69.54& 0.47& 29.77& 22.49& 20.47& 83.17& 93.12
&  & 9.03      & 76.33      & 0.60       & 49.03      & 1.99       & 16.17       & 88.55      & 90.39       \\
                     &                      & $T_{ij}$ = 25             & 9.02& 68.65& 0.35& 29.12& 34.65& 35.88& 83.14& 93.44
&  & 9.60       & 75.11      & 0.32      & 48.64      & 1.32       & 15.89       & 88.53      & 90.61       \\
                     & 60                   & Fully                     & 9.41& 60.12& 0.47& 18.62& 14.96& 14.79& 92.79& 94.27
&  & 9.72      & 72.84      & 0.49      & 41.43      & 0.45       & 13.56       & 95.15      & 94.83       \\
                     &                      & $T_{ij}$ = 51             & 9.44& 58.73& 0.59& 17.96& 48.06& 49.20& 92.80& 94.76
&  & 9.96      & 71.67      & 0.44      & 40.20       & 1.40        & 14.36       & 95.18      & 95.96       \\
                     &                      & $T_{ij}$ = 25             & 9.10& 57.26& 0.58& 16.72& 69.92& 72.75& 92.82& 94.53
&  & 10.05     & 70.30       & 0.32      & 39.23      & 16.72      & 28.18       & 95.16      & 95.78       \\ \hline
400& 30                   & Fully                     & 8.09& 94.58& 0.08& 67.12& 18.27& 28.37& 82.01& 99.16
&  & 9.55      & 96.17      & 0.52      & 83.03      & 5.55       & 31.92       & 89.75      & 98.94       \\
                     &                      & $T_{ij}$ = 51             & 8.19& 94.16& 0.11& 65.26& 35.28& 47.72& 82.14& 99.16
&  & 9.46      & 95.56      & 0.47      & 81.44      & 11.26      & 39.64       & 89.20       & 99.53       \\
                     &                      & $T_{ij}$ = 25             & 8.65& 93.42& 0.08& 63.16& 47.81& 61.40& 82.00& 99.16
&  & 9.88      & 95.25      & 0.40       & 81.22      & 18.42      & 45.47       & 89.24      & 99.17       \\
                     & 60                   & Fully                     & 9.45& 92.18& 0.12& 57.36& 67.38& 72.37& 92.59& 99.78
&  & 9.88      & 94.45      & 0.12      & 75.04      & 1.87       & 20.59       & 95.09      & 99.56       \\
                     &                      & $T_{ij}$ = 51             & 9.78& 91.50& 0.11& 54.86& 91.24& 97.93& 92.56& 99.78
&  & 9.77      & 93.71      & 0.07      & 73.09      & 64.13      & 72.73       & 95.10       & 99.47       \\
                     &                      & $T_{ij}$ = 25             & 9.57& 89.93& 0.10& 53.13& 92.98& 99.91& 92.55& 99.75&  & 9.64      & 92.84      & 0.16      & 71.60       & 88.04      & 93.32       & 95.02      & 99.37       \\ \hline
\end{tabular}
		}	
	\end{center}
\end{table}
    \end{landscape}
    \clearpage
}

\vspace{-0.5cm}

\section{Real data analysis}
\label{sec:real}
\vspace{-0.2cm}
\subsection{HCP data}
\label{sec:real_hcp}
In this section, we apply our proposed MFC to functional magnetic resonance imaging (fMRI) data from the Human Connectome Project (HCP).
This dataset includes  resting-state fMRI scans of subjects, recorded  every 0.72 seconds at $L=1200$ measurement locations (14.4 minutes),  along with their corresponding fluid intelligence \textit{gF} scores. These scores are measured by Raven's Progressive Matrices and serve as indicators of the subjects' ability to process new information, learn, and solve problems \cite[]{cattell1987}. 
Specifically, we consider $n_\low = 73$ subjects with $\textit{gF} \leq 8 $ and $n_\high = 85$ subjects with $\textit{gF} \geq 23 $, and employ a $30$-second observation window ($T_{ij} = 43$) starting from the 12th minute to investigate the relationship between the short time scale brain connectivity \citep{HUTCHISON2013} and  \textit{gF}.
The preprocessing of the raw fMRI data  involved the ICA-FIX preprocessed pipeline \cite[]{glasser2013} and a standard band-pass filter at $[0.01,0.08]$ Hz to filter out frequency bands irrelevant to resting state brain connectivity \cite[]{biswal1995}.   
Following the proposals of \cite{zapata2019} and \cite{miao2022}, we treat signals from different regions of interest (ROIs) as multivariate functional data.
We then implement MFC to construct the resting state brain networks depicting the marginal correlation structures among $p = 83$ ROIs of three well-acknowledged modules in neuroscience study \cite[]{finn2015}: the medial frontal (29 ROIs), frontoparietal (34 ROIs) and default mode modules (20 ROIs), for two groups of subjects, respectively. 

\begin{figure}[!htbp]

\centering
\begin{subfigure}{0.49\linewidth}
  \includegraphics[width=7cm]{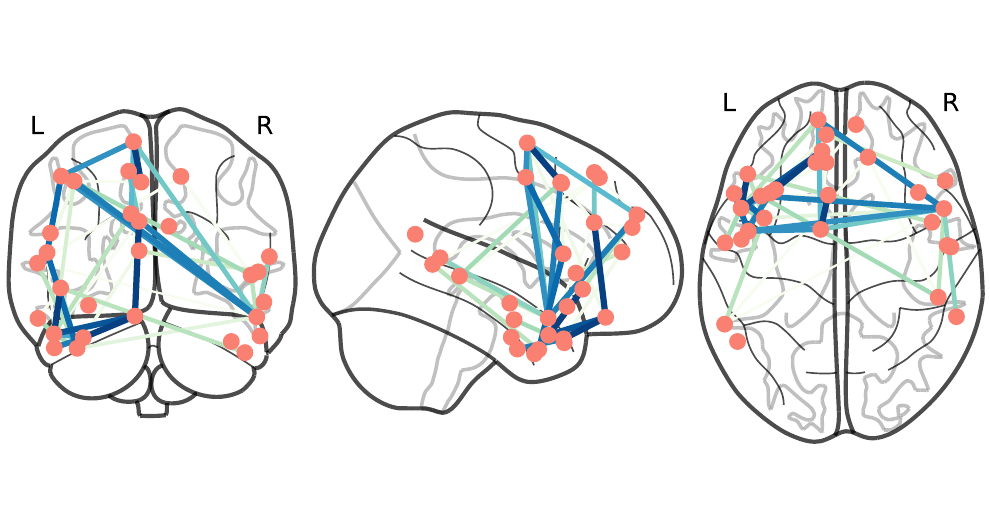}
  \vspace{-0.2cm}
  \caption{$\textit{gF}  \leq 8$: the medial frontal module} 
  \par\smallskip 
    \vspace{-0.2cm}
  \includegraphics[width=7cm]{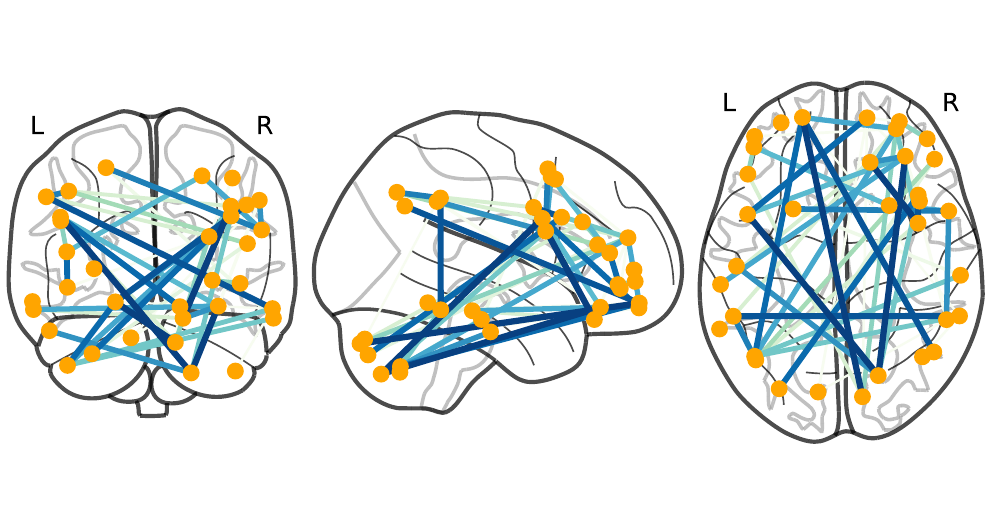}
    \vspace{-0.2cm}
  \caption{$\textit{gF}  \leq 8$: the frontoparietal module} 
    \vspace{-0.2cm}
  \par\smallskip 
  \includegraphics[width=7cm]{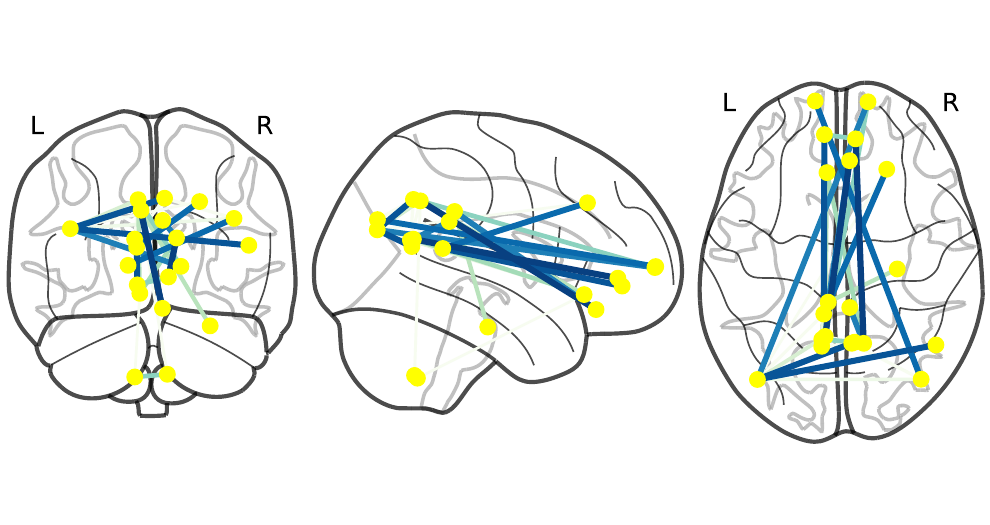}
    \vspace{-0.2cm}
  \caption{$\textit{gF}  \leq 8$: the default mode module} 
\end{subfigure}
\centering
\begin{subfigure}{0.49\linewidth}
  \includegraphics[width=7cm]{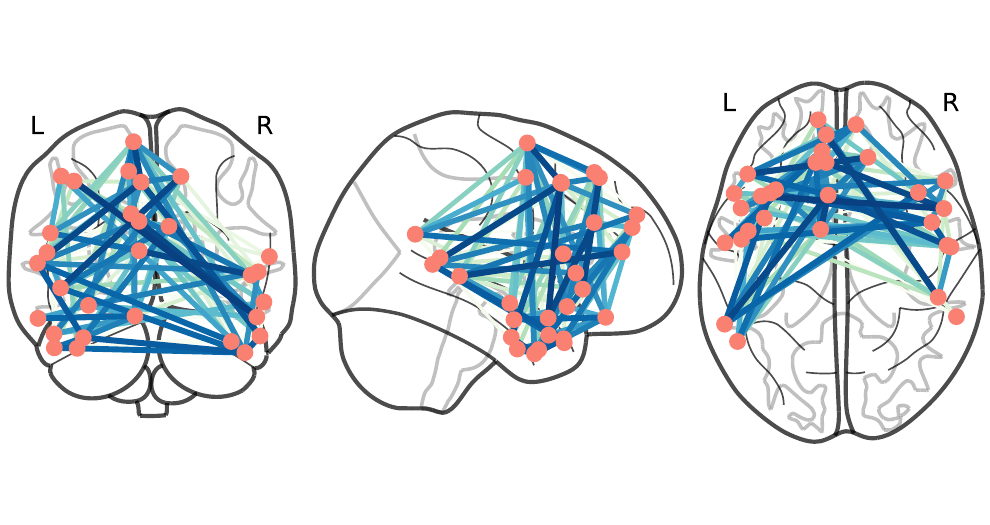}
    \vspace{-0.2cm}
  \caption{$\textit{gF}   \geq 23$: the medial frontal module} 
  \par\smallskip 
    \vspace{-0.2cm}
  \includegraphics[width=7cm]{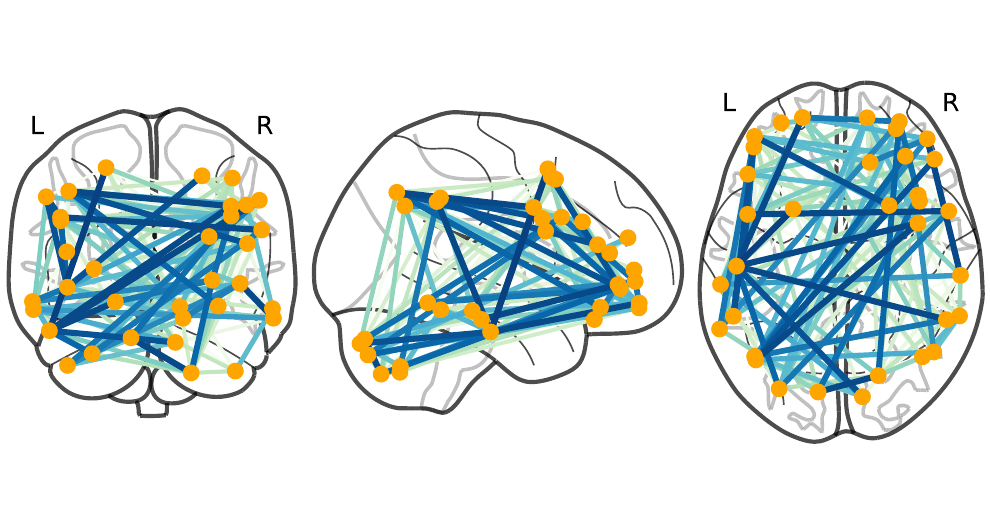}
    \vspace{-0.2cm}
  \caption{$\textit{gF}   \geq 23$: the frontoparietal module} 
    \vspace{-0.2cm}
  \par\smallskip 
  \includegraphics[width=7cm]{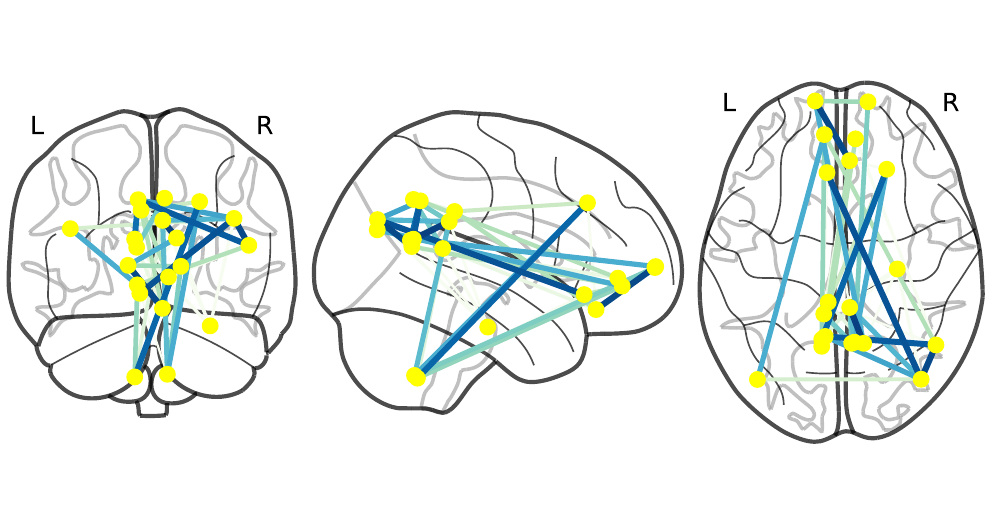}
    \vspace{-0.2cm}
  \caption{$\textit{gF}   \geq 23$: the default mode module} 
\end{subfigure}

\centering
\caption{\label{hcp_network}{ The connectivity strengths at fluid intelligence $\textit{gF} \leq 8$ and $\textit{gF} \geq 23$.
Salmon, orange and yellow nodes represent the ROIs in the medial frontal, frontoparietal and  default mode modules, respectively. 
The edge color from cyan to blue corresponds to the $p$-value from large to small.}}
\end{figure}

Figure~\ref{hcp_network} displays the brain networks based on MFC at the target FDR level $\alpha = 10\%$. Out of a total of $3403$ pairs, we identify 826 connections in subjects with $\textit{gF} \geq 23$ and $365$ in the group with $\textit{gF} \leq 8$. We use the p-value, i.e. $\widetilde{\rm pv}^{(l)}=\bbP\{T_0^{(l)}\geq \widetilde T_n^{(l)}\}$ defined in Section~\ref{sec:error} as the connectivity measure for the $l$th identified pairs, with a larger p-value indicating reduced connectivity strength. 
It is observable that individuals with higher \textit{gF} scores have significantly stronger connectivity within the medial frontal and frontoparietal modules. Notably, despite having a larger number of connections, the brain connectivity in the default mode network appears to be weaker for the higher \textit{gF} group compared to those with lower \textit{gF} scores. These patterns nicely align with the findings in \cite{finn2015}.

\vspace{-0.5cm}
\subsection{EEG data}
\label{sec:real_eeg}
In this section, we illustrate our proposed MFG using the electroencephalogram (EEG) data from an alcoholism study \cite[]{zhang1995}. The study consists
of 122 subjects, out of which $n_a = 77$ are in the alcoholic group and $n_c = 45$ are in the control group. Each
subject was exposed to a picture stimulus while brain activities were measured at 256 time points over a one-second time interval at $64$ electrodes/nodes.  We follow the preprocessing procedure in \cite{zhu2016} and \cite{qiao2019}, which involved averaging and filtering the EEG signals at $\alpha$ frequency bands across all trials. In our analysis, we removed three anatomical landmarked electrodes, X, Y and n, leading to $p=61$. Since the duration for the stimulus in each trial was 300 ms and the subjects' average response time for each stimulus was 653.45 ms, we consider the time from 300--655 ms, which covers the time period of interest, with $T_{ij} = 90$ sampled time points in total. 

\begin{figure}[!htbp]
\centering
\begin{subfigure}{0.49\linewidth}
  \includegraphics[width=7.7cm]{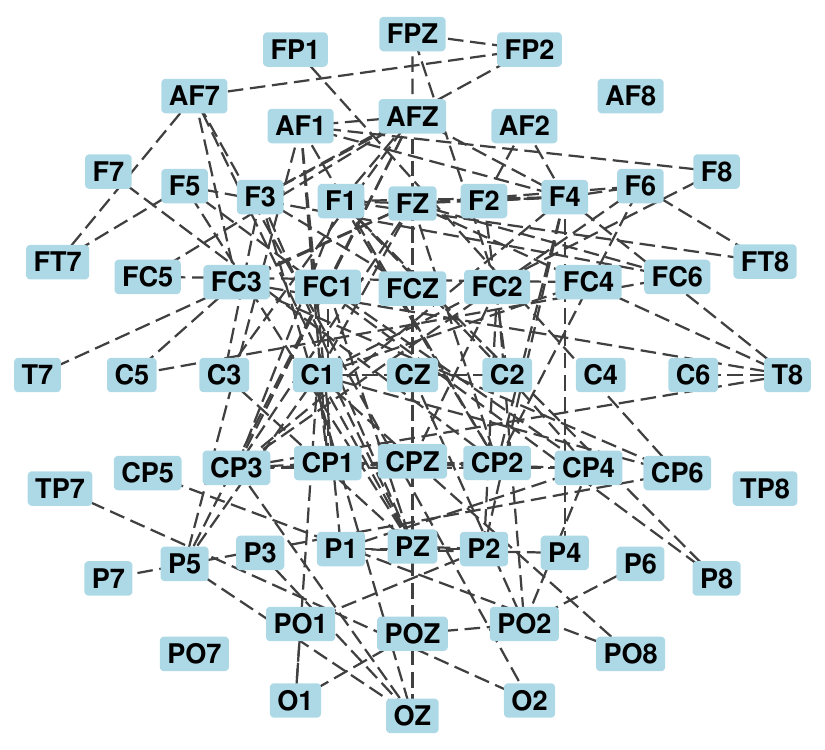}
  \caption{alcohol group} 

\end{subfigure}
\centering
\begin{subfigure}{0.49\linewidth}
  \includegraphics[width=7.7cm]{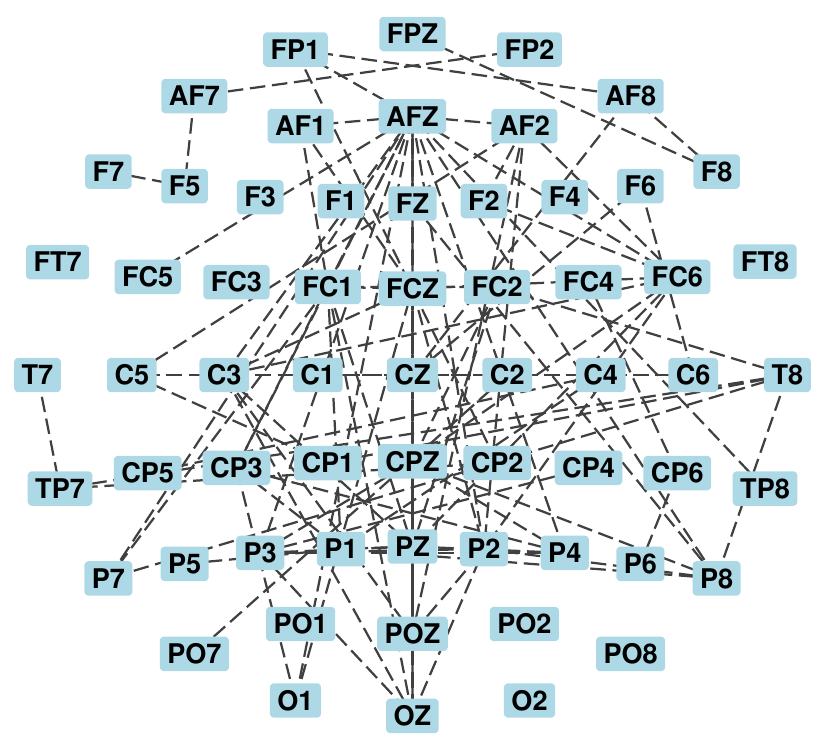}
  \caption{control group} 
\end{subfigure}

\centering
\caption{\label{eeg_network}{Brain networks based on MFG for the alcoholic and control groups.
}}
\end{figure}

We apply MFG to build brain networks. To approximate the $\chi^2$-type mixture in Theorem~\ref{thm:errH0}, we adopted a fast dimension-reduction-based approach \cite[]{rice2019}. Figure~\ref{eeg_network} plots the brain networks constructed at the target FDR level of $\alpha = 2\%$, revealing that $6\%$ and $5.3\%$ of node pairs are recognized as significant connections for the alcohol and control groups, respectively. Unlike the symmetric connectivity observed in the control group, the alcohol group exhibits enhanced connectivity in the left frontal lobe area. This is consistent with existing neuroscience
literature documenting a compensatory elevation in brain activation and functional connectivity, especially in the left frontal region, among individuals with alcohol dependence to achieve comparable levels of task performance to non-alcoholics \citep{desmond2003,
chanraud2011}. 
Additionally, the corpus callosum area (midline nodes in Figure~\ref{eeg_network}) demonstrates a significant number of connections for both groups, highlighting its crucial role in facilitating information transfer between the two hemispheres of the brain \citep{bloom2005,van2011}. These observations again validate the scientific reliability of our proposed method.

 



\section{Discussion}
\label{sec:discuss}
Our test statistic $T_{n,jk}$ (or $\widetilde{T}_{n,jk}$) is constructed based on the $L_2$-norm of the difference between the estimated cross-covariance function $\widehat\Sigma_{jk}(u,v)$ (or $\widetilde{\Sigma}_{jk}(u,v)$) and the true one $\Sigma_{jk}(u,v)$.  For functional data with local spikes, one may use the supremum-norm based test statistics,
${T}_{n,jk}’=n^{1/2}\sup_{(u,v)\in\cU^2}|\widehat\Sigma_{jk}(u,v)|$ and $\widetilde{T}_{n,jk}'=n^{1/2}\sup_{(u,v)\in\cU^2}|\widetilde{\Sigma}_{jk}(u,v)|,
$
upon which our multiple testing procedures in Sections \ref{sec:mt} and \ref{sec:error} remain valid under the associated technical analysis. We next provide a sketch for establishing the theoretical properties of the supremum-norm based test statistics. For fully observed functional data, under Conditions {\rm \ref{cond:covfun}--\ref{cond:eta}} and regular continuity condition for $Z_{ij}(u)$, we can use the central limit theorem for empirical process and continuous mapping theorem to obtain the asymptotic properties of ${T}_{n,jk}'$ under the null and alternative hypotheses. It is not difficult to know the limiting null distribution of ${T}_{n,jk}'$ is $\sup_{(u,v)\in\cU^2} |G_{jk}(u,v)|$, where $G_{jk}(u,v)$ follows a mean-zero Gaussian process with covariance function $\Gamma_{jk}(u_1,v_1,u_2,v_2)$. As $\sup_{(u,v)\in\cU^2} |G_{jk}(u,v)|$ is usually intractable, a parametric bootstrap method could be employed to determine critical values, which would also be used to construct the multiple testing. To investigate the theoretical properties of the multiple testing, we could apply the partition technique to reducing the problem from supremum over $\cU^2$ to the maximum over a grid of pairs,   and then adopt the idea of Gaussian approximation in \cite{CCK2017} and the theory of multiple testing in \cite{Chang2023b} to derive the theoretical results.
Furthermore, with the aid of initial pre-smoothing, this procedure can be extended naturally to the densely observed functional scenario. Although similar technical tools could be used, we need to address the additional complexity that arises from the estimation errors in the pre-smoothing step.



\spacingset{1.08}
\bibliography{paperbib}
\bibliographystyle{dcu}

\appendix
\newpage
\begin{center}
	{\noindent \bf \large Supplementary material to ``Large-scale multiple testing of cross-covariance functions with applications to functional network models"}\\
\end{center}
\begin{center}
	{\noindent Qin Fang, Qing Jiang and Xinghao Qiao
	}
\end{center}
\bigskip

\setcounter{page}{1}
\setcounter{section}{0}
\renewcommand\thesection{\Alph{section}}
\setcounter{lemma}{0}
\renewcommand{\thelemma}{\Alph{section}\arabic{lemma}}
\setcounter{equation}{0}
\renewcommand{\theequation}{S.\arabic{equation}}

\setcounter{table}{0}
\setcounter{figure}{0}
\renewcommand{\thefigure}{S\arabic{figure}}
\renewcommand{\thetable}{S\arabic{table}}

\spacingset{1.68}
This supplementary material contains proofs of main theoretical results in Section~\ref{supp.sec_th}, additional technical proofs in Section~\ref{supp.sec_ad}, and additional empirical results in Section~\ref{supp.sec_emp}.

\section{Proofs of main theoretical results}
\label{supp.sec_th}

\subsection{Proof of Theorem \ref{thm:nonH0}} 
\textbf{Proof of Theorem \ref{thm:nonH0}(i).}
Notice that for any $j,k\in[p]$,
\begin{align*}
	\widehat\Sigma_{jk}(u,v)
	=&~ \frac{1}{n-1}\sum_{i=1}^n \{X_{ij}(u)-\widebar{X}_j(u)\}\{X_{ik}(v)-\widebar{X}_k(v)\} \\
	=&~ \frac{1}{n-1}\sum_{i=1}^n\{Z_{ij}(u)-\bar Z_j(u)\}\{Z_{ik}(v)-\bar{Z}_k(v)\} \\
	=&~ \frac{1}{n-1}\sum_{i=1}^nZ_{ij}(u)Z_{ik}(v)-\frac{n}{n-1}\bar Z_j(u)\bar Z_k(v) \\
	\overset{\triangle}{=}&~ \frac{1}{n-1}\sum_{i=1}^n z_{ijk}(u,v) - \frac{n}{n-1}\bar Z_j(u)\bar Z_k(v)
\end{align*}
where $\widebar{X}_j(u)=n^{-1}\sum_{i=1}^n X_{ij}(u)$, $\bar{Z}_j(u)=n^{-1}\sum_{i=1}^nZ_{ij}(u)$ and $z_{ijk}(u,v),i\in[n]$ are i.i.d. with $\bbE\{z_{ijk}(u,v)\}=\Sigma_{jk}(u,v)$ and
\begin{align*}
	\Gamma_{jk}(u_1,v_1,u_2,v_2)
	=&~ {\cov}\{z_{1jk}(u_1,v_1),z_{1jk}(u_2,v_2)\} \\
	=&~ \bbE\{Z_{1j}(u_1)Z_{1k}(v_1)Z_{1j}(u_2)Z_{1k}(v_2)\} - \Sigma_{jk}(u_1,v_1)\Sigma_{jk}(u_2,v_2)\,.
\end{align*}
By Cauchy--Schwarz inequality and Condition \ref{cond:eta}, it holds that
\begin{align*}
	\bbE(\|z_{1jk}\|_{\cS}^2)
	=&~ \bbE\bigg\{\iint z_{1jk}^2(u,v)\,\md u\md v\bigg\}
	=\bbE\bigg\{\iint Z_{1j}^2(u)Z_{1k}^2(v)\,\md u\md v\bigg\} \\
	\leq&~ \bigg(\bbE\bigg[\bigg\{ \int Z_{1j}^2(u)\,\md u \bigg\}^{2}\bigg]\bigg)^{1/2}
	\bigg(\bbE\bigg[\bigg\{ \int Z_{1k}^2(u)\,\md u \bigg\}^{2}\bigg]\bigg)^{1/2}  \\
	=&~ \{\bbE(\|Z_{1j}\|^4)\}^{1/2} \{\bbE(\|Z_{1k}\|^4)\}^{1/2} <\infty \,.
\end{align*}
Let  ${\rm GP}(\mu,\Sigma)$  denotes  a Gaussian process with mean function $\mu(\cdot)$ and covariance function $\Sigma(\cdot,\cdot)$.
By the central limit theorem of i.i.d. stochastic processes [Theorem 4.12 of \textcolor{blue}{Zhang (2013)}], we have 
\begin{align}\label{eq:zjk}
	\frac{1}{\sqrt{n}}\sum_{i=1}^n \{z_{ijk}(u,v)-\Sigma_{jk}(u,v)\}\overset{d}{\rightarrow} {\rm GP}(0,\Gamma_{jk})\,,
\end{align}
as $n\rightarrow\infty$.
Notice that $\bbE\{\bar Z_j(u)\}=0$ for any $j\in[p]$. Cauchy--Schwarz inequality and Condition \ref{cond:covfun}(ii) imply that for any $j\in[p]$,
\begin{align*}
\sup_{u,v\in\cU}{\cov}\{\bar Z_j(u),\bar Z_j(v)\} = \frac{1}{n} \sup_{u,v\in\cU}\Sigma_{jj}(u,v)\leq \frac{\sup_{u\in\cU}\Sigma_{jj}(u,u)}{n} \rightarrow 0, ~~~\mbox{as}~~ n\rightarrow\infty\,.
\end{align*}
Together with Condition \ref{cond:covfun}(i), it holds that $\bar Z_j(u)$ converges to 0 in probability uniformly over $\cU$ for any $[p]$ and $\sqrt{n}\,\bar{Z}_j(u)\overset{d}{\rightarrow}{\rm GP}(0,\Sigma_{jj})$. Combining with \eqref{eq:zjk}, we have 
\begin{align*}
	\sqrt{n}\{\widehat\Sigma_{jk}(u,v)-\Sigma_{jk}(u,v)\} \overset{d}{\rightarrow} {\rm GP}(0,\Gamma_{jk})
\end{align*}
as $n\rightarrow\infty$. Following from Condition \ref{cond:eta}, it holds that
\begin{align*}
	{\rm tr}(\Gamma_{jk})
	=&~ \iint\Gamma_{jk}\{(u,v),(u,v)\}\, \md u\md v \\
	\leq &~ \iint [\bbE\{Z_{1j}^2(u)Z_{1k}^2(v)\} - \Sigma_{jk}^2(u,v)] \, \md u\md v \\
	\leq &~ 2\iint \bbE\{Z_{1j}^2(u)Z_{1k}^2(v)\}  \, \md u\md v <\infty \,.
\end{align*}
Then applying the continuous mapping theorem for random elements taking values in a Hilbert space \textcolor{blue}{(Billingsley, 1999)} and Theorem 4.2 of \textcolor{blue}{Zhang (2013)}, it holds that
\begin{align}\label{eq:asymTn}
	n\iint\{\widehat\Sigma_{jk}(u,v)-\Sigma_{jk}(u,v)\}^2\,\md u\md v \overset{d}{\rightarrow} \sum_{r=1}^{\infty} \lambda_{jkr}A_{r}~~~{\rm as}~~~n\rightarrow\infty\,,
\end{align}
where $ A_{r}\overset{i.i.d.}{\sim}\chi^2_1 $. The results of \textcolor{blue}{Dauxois et al. (1982)} imply that $|\hat{\lambda}_{jkr}-\lambda_{jkr}|\lesssim \lambda_{jkr}O_{\rm p}(n^{-1/2})$ for any $r\in\mathbb{N}_+$ and $j,k\in[p]$. Together with Condition \ref{cond:eta}, it holds that $|\sum_{r=1}^{\infty} \hat\lambda_{jkr}A_{r} - \sum_{r=1}^{\infty} \lambda_{jkr}A_{r}|\leq \sum_{r=1}^\infty \lambda_{jkr}O_{\rm p}(n^{-1/2})=O_{\rm p}(n^{-1/2})=o_{\rm p}(1)$.
Notice that $\Sigma_{jk}(u,v)=0$ for any $u,v\in\cU$ under $H_{0,jk}$. Thus we  complete the proof of Theorem \ref{thm:nonH0}(i).
$\hfill\Box$\\
\textbf{Proof of Theorem \ref{thm:nonH0}(ii).}
Notice that for any $j,k\in[p]$,
\begin{align*}
	T_{n,jk}
	=&~ \iint[\sqrt{n}\{\widehat\Sigma_{jk}(u,v)-\Sigma_{jk}(u,v)\}]^2 \,\md u\md v
	+n\iint\Sigma_{jk}^2(u,v)\, \md u\md v \\
	&~ + 2 \sqrt{n}\iint\Sigma_{jk}(u,v) [\sqrt{n}\{\widehat\Sigma_{jk}(u,v)-\Sigma_{jk}(u,v)\}] \,\md u\md v \,,
\end{align*}
Since $\iint\Sigma_{jk}^2(u,v)\,\md u\md v >0$, then
following from \eqref{eq:asymTn}, Cauchy--Schwarz inequality and $H_{1,jk}$, it holds that $\iint[\sqrt{n}\{\widehat\Sigma_{jk}(u,v)-\Sigma_{jk}(u,v)\}]^2 \,\md u\md v=O_{\rm p}(1)$, $n\iint\Sigma_{jk}^2(u,v)\, \md u\md v=O(n)$ and $ 2 \sqrt{n}\iint\Sigma_{jk}(u,v) [\sqrt{n}\{\widehat\Sigma_{jk}(u,v)-\Sigma_{jk}(u,v)\}] \,\md u\md v=O_{\rm p}(n^{1/2})$, which implies that
\begin{align*}
	\bbP\{T_{n,jk}\geq T_{0,jk}(\alpha)\} \rightarrow 1
\end{align*}
as $n\rightarrow\infty$ under $H_{1,jk}$.    $\hfill\Box$

\subsection{ Proof of Theorem \ref{thm:nonFDR}}
To prove Theorem \ref{thm:nonFDR}, we need the following lemmas whose proof are given in Section \ref{sec.pf.lm_norm} and \ref{sec.pf.lm_diff}, respectively.

\begin{lemma}
    \label{lm_norm}
Suppose Condition~{\rm\ref{cond:subG}} hold.
For any $q  \in \mathcal{H}_0,$ the Karhunen-Lo\`eve expansion allows us to represent 
 $z_{iq}(u,v):=Z_{ij_q}(u)Z_{ik_q}(v)=\sum_{r=1}^{\infty}a_{ir}^{(q)}\varphi_{r}^{(q)}(u,v),$ where $\lambda_r^{(q)}:=\lambda_{j_qk_qr}$, $\varphi_r^{(q)}(u,v):=\varphi_{j_qk_qr}(u,v)$, and the coefficients $a_{ir}^{(q)} =\iint z_{iq}(u,v) \varphi_{r}^{(q)}(u,v)\,{\rm d}u{\rm d}v$  
are uncorrelated random variables with mean zero and ${\rm Cov}\{a_{ir}^{(q)}, a_{ir'}^{(q)}\} = \lambda_{r}^{(q)}I(r=r')$. 
Then 
there exist some constant $M>0$ such that $\max_{i\in[n]}\max_{q\in \mathcal{H}_0}\sup_{r\in\bbN_+}\bbE\big(\exp[M|a_{ir}^{(q)}|/\{\lambda_{r}^{(q)}\}^{1/2}]\big)=O(1)$.
\end{lemma}

\begin{lemma}\label{lm_diff}
Under Conditions {\rm\ref{cond:covfun}} and {\rm\ref{cond:subG}}, if $p\lesssim n^{\kappa}$ for some constant $\kappa>0$, we have
$\max_{q\in\cH_0}|T_{n}^{(q)}-\mathring{T}_{n}^{(q)}|=O_{\rm p}\{n^{-1/2}(\log p)^{3/2}\}$, where $\mathring{T}_n^{(q)}=\iint\{n^{-1/2}\sum_{i=1}^n z_{iq}(u,v)\}^2\,{\rm d}u{\rm d}v$ for any $i\in[n]$ and $q\in[Q]$. 
\end{lemma}

Recall $Q=p(p-1)/2$ and  $p\leq Cn^{\kappa}$ for some constants $C>0$ and $\kappa>0$. 
By  Lemma \ref{lm_diff}, we have $\max_{q\in\cH_0}|T_n^{(q)}-\mathring{T}_n^{(q)}| =o_{\rm p}\{(\log p)^{-2}\}$. 
Notice that $z_{iq}(u,v)=\sum_{r=1}^\infty a_{ir}^{(q)}\varphi_{r}^{(q)}(u,v)$, which implies
\begin{align*}
  \mathring{T}_n^{(q)}= \iint\bigg\{\frac{1}{\sqrt{n}}\sum_{i=1}^n\sum_{r=1}^\infty a_{ir}^{(q)}\varphi_{r}^{(q)}(u,v)\bigg\}^2\,{\rm d}u{\rm d}v
  =\frac{1}{n}\sum_{i_1=1}^n\sum_{i_2=1}^n \sum_{r=1}^\infty a_{i_1r}^{(q)}a_{i_2r}^{(q)}
  =\sum_{r=1}^\infty \bigg\{\frac{1}{\sqrt{n}}\sum_{i=1}^n a_{ir}^{(q)}\bigg\}^2 \,.
\end{align*}
Denote $\mathring{T}_n^{(q,K)}=\sum_{r=1}^K\{n^{-1/2}\sum_{i=1}^n a_{ir}^{(q)}\}^2$. 
By Markov's inequality and Condition \ref{con_eigenvalue}, we have
\begin{align*}
 &~\bbP\bigg\{ \max_{q\in\cH_0}|\mathring{T}_n^{(q)}-\mathring{T}_n^{(q,K)}|\geq (\log p)^{-1}\bigg\}
 = \bbP\bigg[\max_{q\in\cH_0} \sum_{r=K+1}^\infty \bigg\{\frac{1}{\sqrt{n}}\sum_{i=1}^n a_{ir}^{(q)}\bigg\}^2\geq (\log p)^{-1} \bigg] \\
 \lesssim&~ p^2(\log p)\max_{q\in\cH_0}\sum_{r=K+1}^\infty \bbE\bigg[\bigg\{\frac{1}{\sqrt{n}}\sum_{i=1}^n a_{ir}^{(q)}\bigg\}^2\bigg]
 \lesssim \frac{p^2(\log p)}{n}\max_{q\in\cH_0}\sum_{r=K+1}^\infty \sum_{i=1}^n \lambda_{r}^{(q)} \\
 \lesssim&~ \frac{p^2(\log p)}{K^{\eta_1 - 1}}
 = O(p^{-1})
\end{align*}
provided that $K=p^{(3+\epsilon)/(\eta_1-1)}$ for some sufficiently small constant $\epsilon>0$. 
Let $\lambda_0=\max_{q\in\cH_0}\sup_{r\in\mathbb{N}_+}\lambda_{r}^{(q)}\lesssim 1$.
Define $\hat{a}_{ir}^{(q)}=a_{ir}^{(q)}I\{|a_{ir}^{(q)}|\leq \pi_n\} - \bbE[a_{ir}^{(q)} I\{|a_{ir}^{(q)}| \leq\pi_n\}]$ with $\pi_n=\lambda_0^{1/2}\log(p+n)(4+4\eta_1+2\epsilon')/\{(\eta_1-1)M\}$ for some sufficiently small constant $\epsilon'>\epsilon$. 
Let $\hat{T}_n^{(q,K)}=\sum_{r=1}^K\{ n^{-1/2}\sum_{i=1}^n\hat{a}_{ir}^{(q)}\}^2$. By Cauchy--Schwarz inequality and Lemma \ref{lm_norm},
\begin{align*}
 &\max_{q\in\cH_0}\sup_{r\in\bbN_+} \frac{1}{\sqrt{n}}\sum_{i=1}^n \bbE[|a_{ir}^{(q)}|I\{|a_{ir}^{(q)}|\geq \pi_n\}] \\
 &~~~~~~~~~~~~
 =\max_{q\in\cH_0}\sup_{r\in\bbN_+} \frac{1}{\sqrt{n}}\sum_{i=1}^n \{\lambda_{r}^{(q)}\}^{1/2} \bbE\bigg(\frac{|a_{ir}^{(q)}|}{ \{ \lambda_{r}^{(q)}\}^{1/2}} I\bigg[\frac{|a_{ir}^{(q)}|}{ \{\lambda_{r}^{(q)}\}^{1/2}}\geq \frac{\lambda_0^{1/2}(4+4\eta_1+2\epsilon')}{\{\lambda_{r}^{(q)}\}^{1/2}(\eta_1-1)M}\log(p+n)\bigg]\bigg)  \\
 &~~~~~~~~~~~~
 \lesssim n^{1/2}\max_{q\in\cH_0}\sup_{r\in\bbN_+}\max_{i\in[n]}\bbE\bigg(\frac{|a_{ir}^{(q)}|}{\{\lambda_{r}^{(q)}\}^{1/2}} I\bigg[\frac{|a_{ir}^{(q)}|}{\{\lambda_{r}^{(q)}\}^{1/2}} \geq \frac{4+4\eta_1+2\epsilon'}{(\eta_1-1)M}\log(p+n)\bigg]\bigg) \\
 &~~~~~~~~~~~~
 \lesssim \frac{n^{1/2}}{(p+n)^{(2+2\eta_1+\epsilon')/(\eta_1-1)}} \max_{q\in\cH_0}\sup_{r\in\bbN_+}\max_{i\in[n]}\bbE\bigg(\frac{|a_{ir}^{(q)}|}{\{\lambda_{r}^{(q)}\}^{1/2}}\exp\bigg[\frac{M|a_{ir}^{(q)}|}{2\{\lambda_{r}^{(q)}\}^{1/2}}\bigg]\bigg) \\
 &~~~~~~~~~~~~
 \lesssim \frac{n^{1/2}}{(p+n)^{(2+2\eta_1+\epsilon')/(\eta_1-1)}} \,.
\end{align*}
Hence,
\begin{align*}
 & \bbP\bigg\{\max_{q\in\cH_0} \max_{r\in[K]}\bigg|\frac{1}{\sqrt{n}}\sum_{i=1}^n a_{ir}^{(q)}- \frac{1}{\sqrt{n}}\sum_{i=1}^n \hat{a}_{ir}^{(q)}\bigg|\geq p^{-(3+\eta_1+\epsilon')/(\eta_1-1)}\bigg\} \\
  &~~~~~~~
  \leq \bbP\bigg\{\max_{i\in[n]}\max_{q\in\cH_0} \max_{r\in[K]}|a_{ir}^{(q)}|\geq \pi_n \bigg\} \\
  &~~~~~~~
  \lesssim np^2K\max_{i\in[n]}\max_{q\in\cH_0}\max_{r\in[K]}\bbP\bigg[\frac{|a_{ir}^{(q)}|}{ \{ \lambda_{r}^{(q)}\}^{1/2}}\geq \frac{4+4\eta_1+2\epsilon'}{(\eta_1-1)M}\log(p+n) \bigg]  \\
  &~~~~~~~
  \lesssim \frac{np^2K}{{(p+n)^{(4+4\eta_1+2\epsilon')/(\eta_1-1)}} } \max_{i\in[n]}\max_{q\in\cH_0}\max_{r\in[K]}\bbE\bigg(\exp\bigg[\frac{M|a_{ir}^{(q)}|}{ \{ \lambda_{r}^{(q)}\}^{1/2}}\bigg]\bigg)
  =O(p^{-1})\,.
\end{align*}
Triangle inequality yields
\begin{align*}
  \max_{q\in\cH_0}|\mathring{T}_n^{(q,K)}-\hat{T}_{n}^{(q,K)}|
  \lesssim&~ K\max_{q\in\cH_0}\max_{r\in[K]}\bigg|\frac{1}{\sqrt{n}}\sum_{i=1}^n \hat{a}_{ir}^{(q)}\bigg| \cdot \max_{q\in\cH_0} \max_{r\in[K]}\bigg|\frac{1}{\sqrt{n}}\sum_{i=1}^n a_{ir}^{(q)}- \frac{1}{\sqrt{n}}\sum_{i=1}^n \hat{a}_{ir}^{(q)}\bigg| \\
  &~ + K\max_{q\in\cH_0}\max_{r\in[K]}\bigg|\frac{1}{\sqrt{n}}\sum_{i=1}^n a_{ir}^{(q)}- \frac{1}{\sqrt{n}}\sum_{i=1}^n \hat{a}_{ir}^{(q)}\bigg|^2 \\
  =&~ o_{\rm p}\{(\log p)^{-2}\} \,.
\end{align*}
Therefore, we have $\bbP[\max_{q\in\cH_0}|{T}_n^{(q)}-\hat{T}_n^{(q,K)}|=o\{(\log p)^{-1}\}]=1$. Recall $T_0^{(q)}\overset{d}{=}\sum_{r=1}^{\infty}\hat\lambda_{r}^{(q)}A_r$ with $A_r\overset{i.i.d}{\sim}\chi^2_1$. Denote $T_0^{(q,K)}\overset{d}{=}\sum_{r=1}^{K}\hat\lambda_{r}^{(q)}A_r$ and $\mathring{T}_0^{(q,K)}\overset{d}{=}\sum_{r=1}^{K}\lambda_{r}^{(q)}A_r$. Since $|\hat{\lambda}_{r}^{(q)}-\lambda_{r}^{(q)}|=\lambda_{r}^{(q)}O_{\rm p}(n^{-1/2})$ for any $r\in\mathbb{N}_+$ and $q\in[Q]$ based on \textcolor{blue}{Dauxois et al. (1982)}, by Markov's inequality, triangle inequality  and Lemma \ref{lm_norm}, it holds that 
\begin{align*}
 \bbP\bigg\{\max_{q\in\cH_0}|T_{0}^{(q,K)}-\mathring{T}_0^{(q,K)}|\geq (\log p)^{-1}\bigg\}
   \lesssim (\log p)\mathbb{E}\bigg[\max_{q\in\cH_0}\bigg|\sum_{r=1}^K\{\hat{\lambda}_{r}^{(q)}-\lambda_{r}^{(q)}\}A_r\bigg|\bigg]
   \lesssim \frac{\log p}{n^{1/2}}=o(1) 
\end{align*}
provided that $p\lesssim n^{\kappa}$ for some positive constant $\kappa$.
Following the similar arguments, we can also obtain $\bbP[\max_{q\in\cH_0}|{T}_0^{(q)}-{T}_0^{(q,K)}|=o\{(\log p)^{-1}\}]=1$. By Theorem 1.1 in 
\textcolor{blue}{Zaitsev (1987)}, for $q\in\cH_0$, we have
\begin{align*}
\bbP\{\hat{T}_n^{(q,K)}\geq t\} \leq&~ \bbP\{\mathring{T}_0^{(q,K)}\geq t-\epsilon_n(\log p)^{-1/2}\} + \tilde  c_1 K^{5/2}\exp\bigg\{-\frac{n^{1/2}\epsilon_n}{\tilde c_2 K^3\pi_n(\log p)^{1/2}}\bigg\} \,, \\
\bbP\{\hat{T}_n^{(q,K)}\geq t\} \geq&~ \bbP\{\mathring{T}_0^{(q,K)}\geq t+\epsilon_n(\log p)^{-1/2}\} - \tilde c_1 K^{5/2}\exp\bigg\{-\frac{n^{1/2}\epsilon_n}{\tilde c_2 K^3\pi_n(\log p)^{1/2}}\bigg\} \,,
\end{align*}
where $\tilde c_1,\tilde c_2>0$ are constants and $\epsilon_n=o(1)$.
 Recall $V_{q}=\Phi^{-1}\{1-{\rm pv}^{(q)}\}$ with ${\rm pv}^{(q)}=\bbP\{T_0^{(q)}\geq T_n^{(q)}\}$, and $F_q(z)=\mathbb{P}\{T_n^{(q)}\leq z\}$ for any $z\in\mathbb{R}$.
Then on the event $\{\max_{q\in\cH_0}|{T}_n^{(q)}-\hat{T}_n^{(q,K)}|=o\{(\log p)^{-1}\}\}\bigcap \{\max_{q\in\cH_0}|{T}_0^{(q)}-{T}_0^{(q,K)}|=o\{(\log p)^{-1}\}\}\bigcap \{\max_{q\in\cH_0}|{T}_0^{(q,K)}-\mathring{T}_0^{(q,K)}|=o\{(\log p)^{-1}\}\}$, by the fact that $G(t+o\{(\log p)^{-1/2}\})/G(t)=1+o(1)$ uniformly in $0\leq t\leq 2(\log p)^{1/2}$ with $G(t)=1-\Phi(t)$, it holds that for $q\in\cH_0$,
 \begin{align}\label{eq:FVq}
 		\bbP\{V_{q}\geq t\}
 		=\bbP\big[\bbP\{T_0^{(q)}   \leq T_n^{(q)}\} \geq \Phi(t)\big]
 = G(t)\{1+o(1)\} \,.
 \end{align}
 Let $t_{\max} = (2\log Q - 2\log\log Q)^{1/2}$. We consider two cases: (i) there exists $t\in [0,t_{\max}]$ such that $\widehat{{\rm FDP}}(t)\leq\alpha$, and (ii) $\widehat{{\rm FDP}}(t)>\alpha$ for any $t\in [0,t_{\max}]$.

 In Case (i), by the definition of $\hat{t}$, it holds that $\widehat{\mathrm{FDP}}(t) > \alpha$ for any $t < \hat{t}$. Then we have
 \begin{align*}
 	\frac{Q G(t)}{1 \vee \sum_{q\in [Q]}I(V_q \ge \hat{t})} \geq \frac{Q G(t)}{1 \vee \sum_{q\in [Q]}I(V_q \ge t)}=\widehat{\mathrm{FDP}}(t)> \alpha\,,
 \end{align*}
 based on the fact that  $I(V_q > t ) \ge I(V_q > \hat{t})$ for $t < \hat{t}$.
 By letting $t \uparrow \hat{t}$ in the numerator of the first term in above inequality, we have $\widehat{\mathrm{FDP}}(\hat{t}) \ge \alpha$. On the other hand, based on the definition of $\hat{t}$, there exists a sequence $\{ t_i\}$ with $t_i \ge \hat{t}$ and $t_i \downarrow \hat{t}$ such that $\widehat{\mathrm{FDP}}(t_i) \le \alpha$. Thus we have $I(V_q \geq \hat{t}) \ge I(V_{q} \geq t_i)$, which implies that
 \begin{align*}
 	\frac{Q G(t_i)}{1 \vee \sum_{q\in [Q]}I(V_q \geq \hat{t})}\leq \frac{Q G(t_i)}{1 \vee \sum_{q\in [Q]}I(V_q \geq t_i)} \le \alpha\,.
 \end{align*}
 Letting $t_i \downarrow \hat{t}$ in the numerator of the first term in above inequality, we have $\widehat{\mathrm{FDP}}(\hat{t}) \le \alpha$. Therefore, we have $\widehat{\mathrm{FDP}}(\hat{t}) = \alpha$ in Case (i).

 In Case (ii), we first show that  the threshold of $\hat{t}$ at $(2\log Q)^{1/2}$ leads to no false rejection with probability tending to 1. Notice that for any $t>0$, $G(t)<\phi(t)/t$,  where $\phi(\cdot)$ is the density function of the standard normal distribution $N(0,1)$. Then if $\hat{t} = (2\log Q)^{1/2}$, following from  {\color{red}\eqref{eq:FVq}}, we have
 	\begin{align*}
 		\bbP\Bigg\{ \sum_{q \in \mathcal{H}_0}I( V_q \ge \hat{t} )\geq1 \Bigg\}
 		\leq Q_0 \max_{q \in \mathcal{H}_0}\bbP( V_q \ge \hat{t} )
 		\leq Q_0 G\{(2\log Q)^{1/2} \}\{1+o(1)\} =o(1)
 \end{align*}
 as $Q\rightarrow\infty$, which implies that $\bbP\{{\rm FDP}(\hat{t})=0\}\rightarrow1$ in Case (ii), that is, the probability of false rejection is tending to 0 asymptotically.
 Then
 \begin{align}\label{eq:tm3p1}
 	\mathbb{P}\bigg\{{\rm FDP}(\hat{t})\leq\frac{\alpha Q_0}{Q}+\varepsilon\bigg\}=&~\mathbb{P}\bigg\{{\rm FDP}(\hat{t})\leq\frac{\alpha Q_0}{Q}+\varepsilon,~{\rm Case~(i)~holds}\bigg\}\notag\\
 	&+\mathbb{P}\bigg\{{\rm FDP}(\hat{t})\leq\frac{\alpha Q_0}{Q}+\varepsilon,~{\rm Case~(ii)~holds}\bigg\}\notag\\
 	=&~\mathbb{P}\big\{{\rm Case~(i)~holds}\big\}+\mathbb{P}\big\{{\rm Case~(ii)~holds}\big\} \notag\\
 	&-\mathbb{P}\bigg\{{\rm FDP}(\hat{t})>\frac{\alpha Q_0}{Q}+\varepsilon,~{\rm Case~(i)~holds}\bigg\}\notag\\
 	=&~1-\mathbb{P}\bigg\{{\rm FDP}(\hat{t})>\frac{\alpha Q_0}{Q}+\varepsilon,~{\rm Case~(i)~holds}\bigg\}\,.
 \end{align}
 Notice that $\widehat{{\rm FDP}}(\hat{t})=\alpha$ and $\hat{t}\in [0,t_{\max}]$ under Case (i).
 If we can show
 \begin{align}\label{eq:tm3p2}
 	\sup_{t \in [0,t_{\max}]}\bigg| \frac{ \mathrm{FDP}(t)}{\widehat{\mathrm{FDP}}(t)} - \frac{Q_0}{Q} \bigg| \to 0~~\mbox{in probability}
 \end{align}
 as $n,Q\rightarrow\infty$, then we have $\lim_{n,Q\rightarrow\infty}\mathbb{P}\{{\rm FDP}(\hat{t})\leq\alpha Q_0/Q+\varepsilon\}=1$ by \eqref{eq:tm3p1}.

 To prove \eqref{eq:tm3p2}, it suffices to show that
 \begin{align*}
 	\sup_{t \in [0,t_{\max}]}\Bigg| \frac{ \sum_{q\in \mathcal{H}_0} \{I(V_q \geq t) - G(t)\} }{Q G(t)} \Bigg|  \to 0~~\mbox{in probability}
 \end{align*}
 as $n,Q\rightarrow\infty$. Let $0 = t_0 < t_1 < \cdots < t_s = t_{\max}$ such that $t_i - t_{i-1} = \tilde{v}$ for $i \in [s-1]$ and $t_s - t_{s-1} \le \tilde{v}$, where $\tilde{v} = \{(\log Q)(\log\log Q)^{1/2}\}^{-1/2}$. Then we have $s \asymp t_{\max}/\tilde{v}$. For any $t\in[t_{i-1},t_i]$, it holds that
 \begin{align*}
 	\frac{ \sum_{q\in \mathcal{H}_0} I(V_q \geq t_{i}) }{Q G(t_{i})}\frac{G(t_{i})}{G(t_{i-1})}  \leq \frac{ \sum_{q\in \mathcal{H}_0} I(V_q \geq t) }{Q G(t)}
 	\leq \frac{ \sum_{q\in \mathcal{H}_0} I(V_q \geq t_{i-1}) }{Q G(t_{i-1})}\frac{G(t_{i-1})}{G(t_i)}\,.
 \end{align*}
 Notice that there exists a universal constant $C_2>0$ such that $e^{-t^2/2}\leq \max(C_2,2t)\int_t^\infty e^{-x^2/2}\,{\rm d}x$ for any $t>0$. Then
 \begin{align*}
 	0<1-\frac{G(t_i)}{G(t_{i-1})}=\frac{\int_{t_{i-1}}^{t_i}e^{-x^2/2}\,{\rm d}x}{\int_{t_{i-1}}^\infty e^{-x^2/2}\,{\rm d}x}\leq\frac{\tilde{v}e^{-t_{i-1}^2/2}}{\int_{t_{i-1}}^\infty e^{-x^2/2}\,{\rm d}x}\leq \tilde{v}\max(C_2,2t_{i-1})\,,
 \end{align*}
 which implies that
 $\max_{ i\in [s]}|1-G(t_i)/G(t_{i-1})|\leq 2\tilde{v}t_{\max}\rightarrow0$
 as $Q\rightarrow\infty$. Thus, to prove \eqref{eq:tm3p2}, it suffices to show that
 \begin{align*}
 	\max_{0\leq i \leq s} \Bigg| \frac{ \sum_{q\in \mathcal{H}_0} \{I(V_q \geq t_i) - G(t_i)\} }{Q G(t_i)} \Bigg|  \to 0 ~~\mbox{in probability}
 \end{align*}
 as $n,Q\rightarrow\infty$. Following from  Bonferroni inequality and  Markov's inequality, we have
 \begin{align}\label{eq:tm3p3}
 	& \bbP\bigg[\max_{0 \leq i \leq s} \bigg| \frac{ \sum_{q\in \mathcal{H}_0} \{I(V_q \geq t_i) - G(t_i)\} }{Q G(t_i)} \bigg| \ge \varkappa \bigg]\notag \\
 	&~~~~~~~~~~\leq \sum_{i=0}^{s}\bbP \bigg[ \bigg| \frac{ \sum_{q\in \mathcal{H}_0} \{I(V_q \geq t_i) - G(t_i)\} }{Q G(t_i)} \bigg| \ge \varkappa \bigg] \notag\\
 	&~~~~~~~~~~\leq \frac{C}{\tilde{v}}\int_{0}^{t_{\max}} \bbP \bigg[\bigg| \frac{ \sum_{q\in \mathcal{H}_0} \{I(V_q \geq t) - G(t)\} }{Q G(t)} \bigg| \ge \varkappa \bigg]\,\mathrm{d}t \notag\\
 	&~~~~~~~~~~\leq \frac{C}{\tilde{v}}\int_{0}^{t_{\max}} \bbP \bigg[\bigg| \frac{ \sum_{q\in \mathcal{H}_0} \{I(V_q \geq t) - \bbP(V_q\geq t)\} }{Q G(t)} \bigg| \ge \frac{\varkappa}{2} \bigg]\,\mathrm{d}t  \notag\\
 	&~~~~~~~~~~~~~~ + \frac{C}{\tilde{v}}\int_{0}^{t_{\max}} \bbP \bigg[\bigg| \frac{ \sum_{q\in \mathcal{H}_0} \{\bbP(V_q\geq t)-G(t)\} }{Q G(t)} \bigg| \ge \frac{\varkappa}{2} \bigg]\,\mathrm{d}t \notag\\
 	&~~~~~~~~~~\leq \frac{C}{\tilde{v}\varkappa^2}\int_{0}^{t_{\max}}\mathbb{E}\bigg[\bigg| \frac{ \sum_{q\in \mathcal{H}_0} \{I(V_q \geq t) - \bbP(V_q\geq t)\} }{Q G(t)} \bigg|^2\bigg]\,{\rm d}t
 \end{align}
 for any $\varkappa>0$ and sufficiently large $n$, where $C$ is some positive constant. The second step is based on the relationship between integration and its associated Riemann sum, and the forth step is due to \eqref{eq:FVq}.
 Thus it suffices to prove the following statement is true, 
 \begin{align}\label{eq:toprove5}
 	&~ \int_{0}^{t_{\max}}\mathbb{E}\bigg[\bigg| \frac{ \sum_{q\in \mathcal{H}_0} \{I(V_q \geq t) - \bbP(V_q\geq t)\} }{Q G(t)} \bigg|^2\bigg]\,{\rm d}t \notag\\
 	&~~~~~~~~~~~~~~ =
 	\int_0^{t_{\max}}\sum_{q,q' \in \mathcal{H}_0} \frac{\bbP( V_q \geq t, V_{q'} \geq t) - \bbP( V_q \geq t)\bbP( V_{q'} \geq t)}{Q^2G^2(t)} \,{\rm d}t
 	= o(\tilde{v})\,.
 \end{align}
 We define three subsets
 \begin{align*}
 	\mathcal{H}_{01} =&~ \{ (q, q'): q,q'\in \mathcal{H}_0,\, q = q' \} \,, \\
 	\mathcal{H}_{02} =&~ \{ (q,q'):q, q' \in \mathcal{H}_0,\, q \neq q', q \in \mathcal{S}_{q'}(\gamma) \mbox{~or~} q' \in \mathcal{S}_{q}(\gamma) \} \,,\\
 	\mathcal{H}_{03} =&~ \{(q,q'):q,q'\in\mathcal{H}_{0}\}\backslash (\mathcal{H}_{01} \cup \mathcal{H}_{02}) \,.
 \end{align*}
 Then $\sum_{q,q'\in\mathcal{H}_0}=\sum_{(q,q')\in\mathcal{H}_{01}}+\sum_{(q,q')\in\mathcal{H}_{02}}+\sum_{(q,q')\in\mathcal{H}_{03}}$.

 When $(q,q')\in\mathcal{H}_{01}$, $\bbP( V_q \geq t, V_{q'} \geq t) - \bbP( V_{q} \geq t) \bbP( V_{q'} \geq t) =\bbP( V_{q} \geq t) \{1- \bbP( V_q \geq t)\}\leq \bbP( V_q \geq t)$.  By \eqref{eq:FVq}, it holds that
 	\begin{align*}
 		\sum_{(q,q')\in \mathcal{H}_{01}}\frac{\bbP( V_q \geq t, V_{q'} \geq t) - \bbP( V_{q} \geq t) \bbP( V_{q'} \geq t) }{Q^2G^2(t)} \lesssim \frac{1}{Q G(t)}\,.
 \end{align*}
 Due to $\int_0^{t_{\max}}\{G(t)\}^{-1}\,{\rm d}t\lesssim \exp(2^{-1}t_{\max}^2)=Q(\log Q)^{-1}$, we have
 	\begin{align}
 		\int_0^{t_{\max}}\sum_{(q,q')\in \mathcal{H}_{01}}\frac{\bbP( V_{q} \geq t, V_{q'} \geq t) - \bbP( V_{q} \geq t)\bbP( V_{q'} \geq t) }{Q^2G^2(t)}\,{\rm d}t
 \lesssim&~ \int_0^{t_{\max}}\frac{{\rm d}t}{QG(t)} \notag\\
 \lesssim&~ \frac{1}{\log Q} =o(\tilde{v})\,.\label{eq:toprove51}
 \end{align}

 Since $F_q\{T_{n}^{(q)}\}\sim U[0,1]$ and $F_{q'}\{T_{n}^{(q')}\}\sim U[0,1]$,
 we have
 $		\bbP ( V_q \geq t,  V_{q'} \geq t )
 		\leq \bbP\big[  F_q\{T_{n}^{(q)}\} \geq \Phi(t) , F_{q'}\{T_{n}^{(q')}\} \geq \Phi(t) \big]\{1+o(1)\}
 		\leq \bbP( \zeta_q \geq t , \zeta_{q'} \geq t )\{1+o(1)\} $,
 where $\zeta_q = \Phi^{-1}[F_{q}\{T_{n}^{(q)}\}]\sim N(0,1)$ and $\zeta_{q'} = \Phi^{-1}[F_{q'}\{T_{n}^{(q')}\}]\sim N(0,1)$.
 When $(q,q')\in\mathcal{H}_{02}$, since $\max_{1\leq q\neq q' \leq Q}|\mathrm{Corr}(\zeta_{q}, \zeta_{q'})| \leq c_{\zeta}<1$,   Lemma 2 in \textcolor{blue}{Berman (1962)} implies that $\bbP( \zeta_q\ge t , \zeta_{q'} \ge t )\lesssim t^{-2}\exp\{-t^2/(1+c_{\zeta})\}$ for any $t>C_3$, where $C_3>0$ is a universal constant. Notice that $e^{-t^2/2}\leq \max(C_2,2t)\int_t^\infty e^{-x^2/2}\,{\rm d}x$ for any $t>0$.
 Then
 	$\bbP( V_q \geq t, V_{q'} \geq t) \lesssim t^{-2} \exp\{-t^2/(1+c_{\zeta})\} \lesssim t^{-2} \{\max(C_2,2t)\}^{2/(1+c_{\zeta})} \{G(t)\}^{2/(1+c_{\zeta})}$ for any $t>\max\{C_3,C_2/2\}$.
 Since $\max_{q \in [Q]}|\mathcal{S}_q(\gamma)| = o(Q^{\nu})$, we have $|\mathcal{H}_{02}| = O(Q^{1+\nu})$. Note that $\nu<(1-c_{\zeta})/(1+c_{\zeta})<1$ and $\bbP( V_q \geq t, V_{q'} \geq t) \leq1$ for any $0<t<\max\{C_3,C_2/2\}$. We have
 	\begin{align}
 		&\int_0^{t_{\max}}\sum_{(q,q')\in\mathcal{H}_{02}}\frac{\bbP(V_q\geq t,V_{q'}\geq t) -\bbP(V_q\geq t)\bbP(V_{q'}\geq t)}{Q^2G^2(t)}\,{\rm d}t\notag
 		\leq \int_0^{t_{\max}}\sum_{(q,q')\in\mathcal{H}_{02}}\frac{\bbP(V_{q}\geq t,V_{q'}\geq t)} {Q^2G^2(t)}\,{\rm d}t\notag \\
 		&~~~~~~~~ \lesssim  \int_0^{\max\{1,C_5,C_3/2\}} \frac{{\rm d}t}{Q^{1-\nu}G^2(t)}
 		+ \int_{\max\{1,C_3,C_2/2\}}^{t_{\max}} \frac{{\rm d}t}{Q^{1-\nu}\{G(t)\}^{2c_{\zeta}/(1+c_{\zeta})}}
 		\notag\\
 		&~~~~~~~~ \lesssim \frac{1}{Q^{1-\nu}} + \frac{Q^{(c_{\zeta}-1)/(1+c_{\zeta})+\nu}}{(\log Q)^{(1+3c_{\zeta})/(2+2c_{\zeta})}} =o(\tilde{v})\,.\label{eq:toprove52}
 	\end{align}

 When $(q,q')\in\mathcal{H}_{03}$,  denote $\rho_{q,q'}={\rm Corr}(\zeta_{q},\zeta_{q'})$. By Theorem 2.1.e of \textcolor{blue}{Lin and Bai (2010)},
 \begin{equation}\label{eq:tm3p11}
 	\bbP( \zeta_q\geq t, \zeta_{q'}\geq t )
 	\leq
 	\left\{\begin{aligned}
 		G(t)G\bigg\{\frac{(1-\rho_{q,q'})t}{(1-\rho_{q,q'}^2)^{1/2}}\bigg\}\,, ~~~~~~~~~~&\mbox{if} ~ -1 < \rho_{q,q'} \le 0\,; \\
 		(1+\rho_{q,q'})G(t)G\bigg\{\frac{(1-\rho_{q,q'})t}{(1-\rho_{q,q'}^2)^{1/2}}\bigg\}\,, ~~~~&\mbox{if} ~ 0 \le \rho_{q,q'} < 1\,.
 	\end{aligned}\right.
 \end{equation}
 Note that $|\rho_{q,q'}|\leq(\log Q)^{-2-\gamma}$ for any $(q,q')\in\mathcal{H}_{03}$. When $-(\log Q)^{-2-\gamma} \le \rho_{q,q'} \le 0$, due to $(1-\rho_{q,q'})/(1-\rho_{q,q'}^2)^{1/2} \ge 1$, we have $G\{(1-\rho_{q,q'})t/(1-\rho_{q,q'}^2)^{1/2}\} \le G(t)$, which implies that $\bbP( \zeta_q \geq t, \zeta_{q'} \geq t ) \le G^2(t) \le \{1+(\log Q)^{-1-\gamma}\}G^2(t)$.
 When $0 <\rho_{q,q'} \le (\log Q)^{-2-\gamma}$, by the mean-value theorem,
 there exists $t'$ satisfying $(1-\rho_{q,q'})t/(1-\rho_{q,q'}^2)^{1/2}<t'<t$ such that
 \begin{align*}
 	G\bigg\{\frac{(1-\rho_{q,q'})t}{(1-\rho_{q,q'}^2)^{1/2}}\bigg\} = G(t)+\phi(t')\bigg\{t-\frac{(1-\rho_{q,q'})t}{(1-\rho_{q,q'}^2)^{1/2}} \bigg\}\,.
 \end{align*}
 Then we have
 \begin{align*}
 	G\bigg\{\frac{(1-\rho_{q,q'})t}{(1-\rho_{q,q'}^2)^{1/2}}\bigg\}\{G(t)\}^{-1} = &~ 1+\frac{\phi(t')}{G(t)}\bigg\{t-\frac{(1-\rho_{q,q'})t}{(1-\rho_{q,q'}^2)^{1/2}} \bigg\} \\
 	\leq &~ 1+\frac{t\phi(t')}{t\phi(t)/(1+t^2)}\bigg\{1-\frac{1-\rho_{q,q'}}{(1-\rho_{q,q'}^2)^{1/2}}\bigg\} \\
 	\leq &~ 1+\frac{\phi(t')}{\phi(t)}\cdot2\rho_{q,q'}(1+t^2)
 \end{align*}
 for any $t>0$. For $0<t<t_{\max}$, it holds that
 $
 \phi(t')/\phi(t)=\exp\{(t-{t}')(t+{t}')/2\}
 \leq\exp[t_{\max}^2\{1-(1-\rho_{q,q'})/(1+\rho_{q,q'}^2)^{1/2}\}]
 \leq\exp(2\rho_{q,q'}t_{\max}^2)\leq \exp\{4(\log Q)^{-1-\gamma}\}\lesssim 1$,
 which implies that $G\{(1-\rho_{q,q'})t/(1-\rho_{q,q'}^2)^{1/2}\} \le G(t)[1+O\{(\log Q)^{-1-\gamma}\}]$ for any $t\in[0,t_{\max}]$ if $0<\rho_{q,q'}\leq(\log Q)^{-2-\gamma}$, where the term $O\{(\log Q)^{-1-\gamma}\}$ holds uniformly over $t\in[0,t_{\max}]$. Then $(1+\rho_{q,q'})G(t)G\{(1-\rho_{q,q'})t/(1-\rho_{q,q'}^2)^{1/2}\}\leq G^2(t)[1+O\{(\log Q)^{-1-\gamma}\}]$ for any $t\in[0,t_{\max}]$ if $0<\rho_{q,q'}\leq(\log Q)^{-2-\gamma}$. By \eqref{eq:tm3p11},  $\bbP( \zeta_{q} \geq t, \zeta_{q'}\geq t ) \le [1+O\{(\log Q)^{-1-\gamma}\}]G^2(t)$ for $t \in [0,t_{\max}]$. Thus we have
$
 		\max_{(q,q') \in \mathcal{H}_{03}}\bbP( V_{q}\geq t, V_{q'}\geq t) \le [1+O\{(\log Q)^{-1-\gamma}\}]G^2(t)$
 for $t \in [0,t_{\max}]$. It then holds that
 	\begin{align*}
 		&\int_0^{t_{\max}}\sum_{(q,q')\in \mathcal{H}_{03}}\frac{\bbP(V_{q}\geq t,V_{q'}\geq t) -\bbP(V_{q}\geq t)\bbP(V_{q'}\geq t) }{Q^2G^2(t)}\,{\rm d}t \\
 		&~~~~~~~~~~~~~ \leq O\{(\log Q)^{-1-\gamma}\}\cdot\int_0^{t_{\max}}1\,{\rm d}t  =o(\tilde{v})\,.
 \end{align*}
 Together with \eqref{eq:toprove51} and \eqref{eq:toprove52}, we know \eqref{eq:toprove5} holds. Then $\lim_{n,Q\rightarrow\infty}\mathbb{P}\{{\rm FDP}(\hat{t})\leq\alpha Q_0/Q+\varkappa\}=1$ for any $\varkappa>0$. Since ${\rm FDR}(t)={\bbE}\{{\rm FDP}(t)\}$, it holds that $\limsup_{n,Q\to\infty}{\rm FDR}(\hat{t})\leq\alpha Q_0/Q$. We complete the proof of Theorem \ref{thm:nonFDR}. $\hfill\Box$

\subsection{Proof of Theorem \ref{thm:errH0}}
\textbf{Proof of Theorem \ref{thm:errH0}(i).}
Notice that $\widetilde{T}_{n,jk}=n\|\widetilde{\Sigma}_{jk}\|_{\cS}^2 $ 
and $\widetilde{\Sigma}_{jk}(u,v)=(n-1)^{-1}\sum_{i=1}^n\{\widehat{X}_{ij}-\widebar{\widehat{X}}_j(u)\}\{\widehat{X}_{ik}(v) -\widebar{\widehat{X}}_k(v)\}$ with $\widehat{X}_{ij}(u)=X_{ij}(u)+e_{ij}(u)$. Then it holds that
\begin{align} \label{eq.tildeSigma}
  \widetilde{\Sigma}_{jk}(u,v)=&~ \frac{1}{n-1}\sum_{i=1}^n\{X_{ij}(u)+e_{ij}(u) -\widebar{X}_j(u)-\bar{e}_j(u)\}\{X_{ik}(v)+e_{ik}(v) -\widebar{X}_k(v)-\bar{e}_k(v)\} \notag\\
  =&~ \frac{1}{n-1}\sum_{i=1}^n\{X_{ij}(u)-\widebar{X}_j(u)\}\{X_{ik}(v)-\widebar{X}_k(v)\} \notag\\
  &~ + \frac{1}{n-1}\sum_{i=1}^n\{X_{ij}(u)-\widebar{X}_j(u)\}\{e_{ik}(v)-\bar{e}_k(v)\} \notag\\
  &~ + \frac{1}{n-1}\sum_{i=1}^n\{e_{ij}(u)-\bar{e}_j(u)\}\{X_{ik}(v)-\widebar{X}_k(v)\} \notag\\
  &~ +\frac{1}{n-1}\sum_{i=1}^n\{e_{ij}(u)-\bar{e}_j(u)\}\{e_{ik}(v)-\bar{e}_k(v)\}  \notag\\
  \overset{\triangle}{=}&~ \widehat\Sigma_{jk}(u,v) +  \widehat{\Sigma}^{X e}_{jk}(u,v) +  \widehat{\Sigma}^{X e}_{kj}(u,v)  +  \widehat{\Sigma}^{e e}_{jk}(u,v)  
\end{align}
for any $u,v\in\cU$. 
By Condition \ref{cond:Xe} with $p=1$, we have
\begin{align*}
   n\|\widehat{\Sigma}^{X e}_{jk}\|_{\cS}^2   =  O_{\rm p}(n^{1-2a_1}) =  n\|\widehat{\Sigma}^{X e}_{kj}\|_{\cS}^2   \,,
\end{align*}
and $n\|\widehat{\Sigma}^{e e}_{jk}\|_{\cS}^2=O_{\rm p}(n^{1-2a_3})$, where $a_1,a_3>1/2$. Following from Cauchy--Schwarz inequality, it holds that $\widetilde{T}_{n,jk}=T_{n,jk}+o_{\rm p}(1)$. Under  Conditions {\rm \ref{cond:covfun}--\ref{con_eigenvalue}} and $H_{0,jk}$, following the same arguments as the proof of Theorem \ref{thm:nonH0}(i), we can derive the results of Theorem \ref{thm:errH0}(i).  $\hfill\Box$\\
\textbf{Proof of Theorem \ref{thm:errH0}(ii).}
Following the similar arguments as  the proof of Theorem \ref{thm:errH0}(i), based on \eqref{cond.cov} and Cauchy--Schwarz inequality, we have $\widetilde{T}_{n,jk}=T_{n,jk}+o_{\rm p}(1)$.
Then it holds that
\begin{align*}
	\widetilde{T}_{n,jk}  
	=&~ \iint[\sqrt{n}\{\widehat\Sigma_{jk}(u,v)-\Sigma_{jk}(u,v)\}]^2 \,\md u\md v
	+n\iint\Sigma_{jk}^2(u,v)\, \md u\md v \\
	&~ + 2 \sqrt{n}\iint\Sigma_{jk}(u,v) [\sqrt{n}\{\widehat\Sigma_{jk}(u,v)-\Sigma_{jk}(u,v)\}] \,\md u\md v +o_{\rm p}(1) \,,
\end{align*}
Since $\iint\Sigma_{jk}^2(u,v)\,\md u\md v >0$, then
following from \eqref{eq:asymTn}, Cauchy--Schwarz inequality and $H_{1,jk}$,  it holds that $\iint[\sqrt{n}\{\widehat\Sigma_{jk}(u,v)-\Sigma_{jk}(u,v)\}]^2 \,\md u\md v=O_{\rm p}(1)$, $n\iint\Sigma_{jk}^2(u,v)\, \md u\md v=O(n)$ and $ 2 \sqrt{n}\iint\Sigma_{jk}(u,v) [\sqrt{n}\{\widehat\Sigma_{jk}(u,v)-\Sigma_{jk}(u,v)\}] \,\md u\md v=O_{\rm p}(n^{1/2})$, which implies that
$
	\bbP\{\widetilde{T}_{n,jk}\geq T_{0,jk}(\alpha)\} \rightarrow 1$
as $n\rightarrow\infty$ under $H_{1,jk}$.    $\hfill\Box$

\subsection{Proof of Theorem \ref{thm:errFDR}}
Recall $\widetilde{T}_{n}^{(q)}=n\|\widetilde{\Sigma}_{j_qk_q}\|_{\cS}^2 $ and $\mathring{T}_n^{(q)}=\iint\{n^{-1/2}\sum_{i=1}^n Z_{ij_q}(u)Z_{ik_q}(v)\}^2\,{\rm d}u{\rm d}v$ for each $q\in[Q]$.
Triangle inequality and   Lemma \ref{lm_diff} yield
\begin{align*}
   \max_{q\in\cH_0}|\widetilde{T}_n^{(q)}-\mathring{T}_n^{(q)}|
   \leq&~  \max_{q\in\cH_0}|\widetilde{T}_n^{(q)}-{T}_n^{(q)}| + \max_{q\in\cH_0}|{T}_n^{(q)}-\mathring{T}_n^{(q)}| = \max_{q\in\cH_0}|\widetilde{T}_n^{(q)}-{T}_n^{(q)}| + O_{\rm p}\{n^{-1/2}(\log p)^{3/2}\} \,.
\end{align*}
Following the similar arguments as the proof of Lemma \ref{lm_diff}, based on \eqref{eq.tildeSigma}, triangle inequality, Lemma 1 in the supplementary material of {\color{blue}  Zapata et al. (2022)} and Condition \ref{cond:Xe}, we also have
\begin{align*}
 &~\max_{q\in\cH_0}|\widetilde{T}_n^{(q)}-{T}_n^{(q)}| \\
 &~~~~~~~~~~\lesssim  n\max_{q\in\cH_0}\|\widehat{\Sigma}_{j_qk_q}^{Xe}+\widehat{\Sigma}_{k_qj_q}^{Xe}+\widehat{\Sigma}_{j_qk_q}^{ee}\|_{\cS}^2   + n\max_{q\in\cH_0} \|\widehat{\Sigma}_{j_qk_q}\|_{\cS} \|\widehat{\Sigma}_{j_qk_q}^{Xe}+\widehat{\Sigma}_{k_qj_q}^{Xe}+\widehat{\Sigma}_{j_qk_q}^{ee}\|_{\cS} \\
 &~~~~~~~~~~ = O_{\rm p}\{n^{1-2a_1}(\log p)^{2a_2}+n^{1-2a_3}(\log p)^{2a_4}+n^{1/2-a_1}(\log p)^{1/2+a_2}+n^{1/2-a_3}(\log p)^{1/2+a_4}\} \\
  &~~~~~~~~~~ = O_{\rm p}\{n^{1/2-a_1}(\log p)^{1/2+a_2}+n^{1/2-a_3}(\log p)^{1/2+a_4}\} 
\end{align*}
provided that $p\lesssim n^\kappa$ for some constant $\kappa>0$.
Hence, we obtain that $  \max_{q\in\cH_0}|\widetilde{T}_n^{(q)}-\mathring{T}_n^{(q)}|= O_{\rm p}\{n^{1/2-a_1}(\log p)^{1/2+a_2}+n^{1/2-a_3}(\log p)^{1/2+a_4}+n^{-1}\log p\} =o_{\rm p}\{(\log p)^{-2}\}$. 
Then following the same arguments as the proof of Theorem \ref{thm:nonFDR},
it is easy to derive the results of Theorem \ref{thm:errFDR}.$\hfill\Box$

\subsection{Proof of Proposition \ref{thm.partial}}
Denote $\be_0=(1,0)^{\T}$, $\widetilde{\bU}_{ijt}=\{1,(U_{ijt}-u)/h_j\}^{\T}$. Define
\begin{align*}
 \widehat{\bS}_{ij}(u)=\frac{1}{T_{ij}}\sum_{t=1}^{T_{ij}} \widetilde{\bU}_{ijt}\widetilde{\bU}_{ijt}^{\T}K_{h_j}(U_{ijt}-u)\,, ~~~~
 \widehat{\bR}_{ij}(u)=\frac{1}{T_{ij}}\sum_{t=1}^{T_{ij}} \widetilde{\bU}_{ij}\{W_{ijt} - Z_{ij}(u)\}K_{h_j}(U_{ijt}-u) \,.
\end{align*}
Let $\widetilde{\bX}_i=(\widetilde{X}_{i1},\ldots,\widetilde{X}_{ip})^{\T}$ with $\widetilde{X}_{ij}(u)=\be_0^{\T}[\bbE\{\widehat{\bS}_{ij}(u)\}]^{-1}\bbE\{\widehat{\bR}_{ij}(u)\}$.
To construct Proposition \ref{thm.partial}, we need the following lemma whose proof is given in Sections \ref{pf.lm_vb} .
\begin{lemma}\label{partial.vb}
   Under  Conditions~{\rm \ref{cond:subG}, \ref{cond.subgauss}--\ref{cond.cont}}, 
we have that {\rm  (i)} $\max_{i\in [n], j \in[p]}
\|\widehat X_{ij} - Z_{ij} - \widetilde{X}_{ij}\| = O_p\{\log(p \vee n)^{1/2} (Th)^{-1/2}\}
$; {\rm (ii)} $\|\widetilde{X}_{ij}-\mu_{j}\| = O(h^2)$.
\end{lemma}

Notice that
\begin{align*}
\widehat{\Sigma}_{jk}^{Xe}(u,v)
=&~ \frac{1}{n-1}\sum_{i=1}^n\{X_{ij}(u)-\widebar{X}_j(u)\} [\{\widehat{X}_{ik}(v)-X_{ik}(v)\}-\{\widebar{\widehat{X}}_k(v)-\widebar{X}_k(v)\}] \\
=&~ \frac{1}{n-1}\sum_{i=1}^n Z_{ij}(u)\{\widehat{X}_{ik}(v)-X_{ik}(v)\} 
- \bigg\{\frac{1}{n}\sum_{i=1}^n Z_{ij}(u)\bigg\}\bigg[\frac{1}{n-1}\sum_{i=1}^n\{\widehat{X}_{ik}(v)-X_{ik}(v)\}\bigg]. 
\end{align*}
By Cauchy--Schwarz inequality, we obtain that
\begin{align*}
    \max_{j\neq k } \|\widehat \Sigma_{jk}^{Xe}\| & \lesssim \max_{j \in [p]}\sqrt{\frac{1}{n} \sum_{i = 1}^n \|Z_{ij}\|^2}\max_{k \in [p]}\sqrt{\frac{1}{n} \sum_{i = 1}^n \|\widehat X_{ik}-X_{ik}\|^2}
    \\& \lesssim \max_{i \in [n], j \in [p]} \|Z_{ij}\| \max_{i \in [n], k \in [p]} \|\widehat X_{ik}-X_{ik}\|
\end{align*}
By Condition~\ref{cond:subG}  and Lemma 1 in the supplementary material of {\color{blue}  Zapata et al. (2022)}, we have $\max_{i \in [n], j \in [p]} \|Z_{ij}\| = O_{\rm p}\{\log(p \vee n)^{1/2}\}$.
Lemma~\ref{partial.vb} implies that
\begin{align} 
    \max_{i \in [n], j \in [p]}\|\widehat X_{ij} - X_{ij}\| &\leq \max_{i \in [n], j \in [p]}\|\widehat X_{ij} - Z_{ij} -\widetilde{X}_{ij}\| + \max_{i\in[n],j\in[p]}\|\widetilde{X}_{ij} - \mu_j\| \nonumber
    \\&=O_{\rm p}\{\log(p \vee n)^{1/2} (Th)^{-1/2}+h^2\}\label{eq.hatX}.
\end{align}
Combining the above results, we obtain that $\max_{j,k\in[p]}\|\widehat{\Sigma}_{jk}^{X e}\|_{\cS} = 
O_{\rm p} \big[\log(p \vee n)(Th)^{-1/2} + \{\log(p\vee n)\}^{1/2} h^2\big]$. Observe that
\begin{align*}
\widehat{\Sigma}_{jk}^{ee}(u,v)
=&~ \frac{1}{n-1}\sum_{i=1}^n\{\widehat{X}_{ij}(u)-X_{ij}(u)\}\{\widehat{X}_{ik}(v)-X_{ik}(v)\} \\
&~ -\frac{n}{n-1}\bigg[\frac{1}{n}\sum_{i=1}^n\{\widehat{X}_{ij}(u)-X_{ij}(u)\}\bigg] \bigg[\frac{1}{n}\sum_{i=1}^n\{\widehat{X}_{ik}(v)-X_{ik}(v)\}\bigg].
\end{align*} 
This, together with \eqref{eq.hatX}, implies that
\begin{align*}
\max_{j\neq k}\|\widehat{\Sigma}_{jk}^{ee}\|_\cS
\lesssim  \Big(\max_{i \in [n], j \in [p]} \|\widehat X_{ij}-X_{ij}\|\Big)^2 = 
O_{\rm p}\big\{\log(p \vee n)(Th)^{-1} + h^4\big\}.
\end{align*}
Hence, we complete the proof of  Proposition \ref{thm.partial}.
$\hfill\Box$

\subsection{ Proof of Theorem \ref{thm:noncondH0}}
\label{supp_th5}
Notice that $
\hat{\varepsilon}_{i,jk}(u)=Y_{ij}(u)-\int_{\cU}\bY_{i,-j,-k}^{\T}(v)\widehat{\boldsymbol{\beta}}_{jk}(u,v)\,{\rm d}v$
and
 $\bar{\hat\varepsilon}_{jk}(u)=n^{-1}\sum_{i=1}^n\hat\varepsilon_{i,jk}(u)$. Then
\begin{align*}
	\hat{\varepsilon}_{i,jk}(u)-\bar{\hat\varepsilon}_{jk}(u)
	=&~ Y_{ij}(u) - \widebar{Y}_{j}(u)- \int_{\cU}\{\bY_{i,-j,-k}(v)-\widebar{\bY}_{-j,-k}(v)\}^{\T} \widehat{\boldsymbol{\beta}}_j(u,v)\,{\rm d}v \\
	=&~ \int_{\cU}\{\bY_{i,-j,-k}(v)-\widebar{\bY}_{-j,-k}(v)\}^{\T} \{\boldsymbol{\beta}_{jk}(u,v)-\widehat{\boldsymbol{\beta}}_{jk}(u,v)\} \,{\rm d}v + \varepsilon_{i,jk}(u)-\bar{\varepsilon}_{jk}(u)  \\
	=&~ \int_{\cU}\{\bY_{i,-j,-k}(v)-\widebar{\bY}_{-j,-k}(v)\}^{\T} \{\boldsymbol{\beta}_{jk}(u,v)-\widehat{\boldsymbol{\beta}}_{jk}(u,v)\} \,{\rm d}v + \tilde{\varepsilon}_{i,jk}(u)  \,,
\end{align*}
where $\widebar{Y}_j(u)=n^{-1}\sum_{i=1}^n Y_{ij}(u)$, $\widebar\bY_{-j,-k}=n^{-1}\sum_{i=1}^n\bY_{i,-j,-k}$, $\bar\varepsilon_{jk}(u)=n^{-1}\sum_{i=1}^n\varepsilon_{i,jk}(u)$ and $\tilde{\varepsilon}_{i,jk}(u)=\varepsilon_{i,jk}(u)-\bar{\varepsilon}_{jk}(u)$ for any $u\in\cU$.
Recall $
\widetilde\Sigma_{jk}^{\varepsilon}(u,v)=(n-1)^{-1}\sum_{i=1}^n\{\hat{\varepsilon}_{i,jk}(u)-\bar{\hat\varepsilon}_{jk}(u)\} \{\hat{\varepsilon}_{i,kj}(v)-\bar{\hat\varepsilon}_{kj}(v)\} $. Hence it holds that
\begin{align}\label{eq:sjksplit}
&\{\hat{\varepsilon}_{i,jk}(u)-\bar{\hat\varepsilon}_{jk}(u)\} \{\hat{\varepsilon}_{i,kj}(v)-\bar{\hat\varepsilon}_{kj}(v)\}  \notag\\
&~~~~~~~~~~~  = \bigg[\tilde{\varepsilon}_{i,jk}(u)-\int_{\cU}\{\bY_{i,-j,-k}(t)-\widebar\bY_{-j,-k}(t)\}^{\T} \{\widehat{\boldsymbol{\beta}}_{jk}(u,t)-\boldsymbol{\beta}_{jk}(u,t)\}\,{\rm d}t \bigg]   \notag\\
&~~~~~~~~~~~~~~ \times  \bigg[\tilde{\varepsilon}_{i,kj}(v)-\int_{\cU}\{\bY_{i,-j,-k}(s)-\widebar\bY_{-j,-k}(s)\}^{\T} \{\widehat{\boldsymbol{\beta}}_{kj}(v,s)-\boldsymbol{\beta}_{kj}(v,s)\}\,{\rm d}s \bigg]   \notag\\
&~~~~~~~~~~~  =\tilde{\varepsilon}_{i,jk}(u)\tilde{\varepsilon}_{i,kj}(v)
-\tilde{\varepsilon}_{i,jk}(u)\int_{\cU}\{\bY_{i,-j,-k}(s)-\widebar\bY_{-j,-k}(s)\}^{\T} \{\widehat{\boldsymbol{\beta}}_{kj}(v,s)-\boldsymbol{\beta}_{kj}(v,s)\}\,{\rm d}s  \notag\\
&~~~~~~~~~~~~~~ - \tilde{\varepsilon}_{i,kj}(v)\int_{\cU}\{\bY_{i,-j,-k}(t)-\widebar\bY_{-j,-k}(t)\}^{\T} \{\widehat{\boldsymbol{\beta}}_{jk}(u,t)-\boldsymbol{\beta}_{jk}(u,t)\}\,{\rm d}t  \notag\\
&~~~~~~~~~~~~~~ +\int_{\cU}\int_{\cU} \{\widehat{\boldsymbol{\beta}}_{jk}(u,t)-\boldsymbol{\beta}_{jk}(u,t)\}^{\T} \{\bY_{i,-j,-k}(t)-\widebar\bY_{-j,-k}(t)\}\{\bY_{i,-j,-k}(s)-\widebar\bY_{-j,-k}(s)\}^{\T} \notag\\ &~~~~~~~~~~~~~~~~~~~~~~~~~~~~ \times \{\widehat{\boldsymbol{\beta}}_{kj}(v,s)-\boldsymbol{\beta}_{kj}(v,s)\}\, {\rm d} t{\rm d}s   \notag\\
&~~~~~~~~~~~  \overset{\triangle}{=} {\rm I}_i(j,k,u,v) -{\rm II}_i(j,k,u,v) - {\rm III}_i(j,k,u,v) + {\rm IV}_i(j,k,u,v) \,.
\end{align}
Notice that $\widetilde T_{n,jk}^{\varepsilon}= n\|\widetilde\Sigma_{jk}^{\varepsilon}\|_{\cS}^2$ and  $\widetilde\Sigma_{jk}^{\varepsilon}(u,v)$ exhibits a similar decomposition as $\widetilde\Sigma_{jk}(u,v)$ in \eqref{eq:statH0d}, where $(n-1)^{-1}\sum_{i=1}^n\{{\rm II}_i(j,k,u,v)+ {\rm III}_i(j,k,u,v)\}$ and $(n-1)^{-1}\sum_{i=1}^n{\rm IV}_i(j,k,u,v)$ corresponds to $\widehat \Sigma_{jk}^{Xe}(u,v)$ and $\widehat \Sigma_{jk}^{ee}(u,v)$ in Section~\ref{sec:error}, respectively. 
In the sequel, we will consider the convergence rates of $\max_{j,k\in[p]}\int_{\cU}\int_{\cU}\{(n-1)^{-1}\sum_{i=1}^n{\rm IV}_i(j,k,u,v)\}^2\,{\rm d}u{\rm d}v$, $\max_{j,k\in[p]}\int_{\cU}\int_{\cU}\{(n-1)^{-1}\sum_{i=1}^n{\rm II}_i(j,k,u,v)\}^2\,{\rm d}u{\rm d}v$ and $\max_{j,k\in[p]}\int_{\cU}\int_{\cU}\{(n-1)^{-1}\sum_{i=1}^n{\rm III}_i(j,k,u,v)\}^2\,{\rm d}u{\rm d}v$, respectively, which will in turn specify the rates defined in \eqref{cond.cov}.

\underline{\it Convergence rate of $\max_{j,k\in[p]}\int_{\cU}\int_{\cU}\{(n-1)^{-1}\sum_{i=1}^n{\rm IV}_i(j,k,u,v)\}^2\,{\rm d}u{\rm d}v$.} By   Inequality \ref{lemma_l1_inequality}, Propositions~\ref{lm_beta} and \ref{lm_betasigma},
\begin{align*}
& \max_{j,k\in[p]}\int_{\cU}\int_{\cU} \bigg\{\frac{1}{n-1}\sum_{i=1}^n {\rm IV}_i(j,k,u,v)\bigg\}^2\,{\rm d}u{\rm d}v  \\
&~~~~ = \max_{j,k\in[p]}\int_{\cU}\int_{\cU}\bigg[\int_{\cU}\int_{\cU}\{\widehat{\boldsymbol{\beta}}_{jk}(u,t)-\boldsymbol{\beta}_{jk}(u,t) \}^{\T}\widehat\bXi_{-j,-k}(t,s)\{\widehat{\boldsymbol{\beta}}_{kj}(v,s)-\boldsymbol{\beta}_{kj}(v,s)\} \,{\rm d}t {\rm d}s  \bigg]^2 \,{\rm d}u{\rm d}v  \\
&~~~~ \lesssim \max_{j,k\in[p]}\bigg\|\int_{\cU} \widehat\bXi_{-j,-k}(\cdot,s)\{\widehat{\boldsymbol{\beta}}_{kj}(\cdot,s)-\boldsymbol{\beta}_{kj}(\cdot,s)\}\, {\rm d}s\bigg\|_{\cS,\max}^2\max_{j,k\in[p]}\|\widehat{\boldsymbol{\beta}}_{jk}-\boldsymbol{\beta}_{jk}\|_{\cS,1}^2 \\
&~~~~  =O_{\rm p}(\delta_{1n}^2\delta_{2n}^2) \,.
\end{align*}
Hence, \begin{equation}
\label{a2.rate}
\delta_{1n}\delta_{2n} \asymp n^{-a_{3}}(\log p)^{a_4}.
\end{equation}

\underline{\it Convergence rate of $\max_{j,k\in[p]}\int_{\cU}\int_{\cU}\{(n-1)^{-1}\sum_{i=1}^n{\rm II}_i(j,k,u,v)\}^2\,{\rm d}u{\rm d}v$
.}
Recall that  $ \varepsilon_{i,jk}(\cdot)$ and $\bY_{i,-j,-k}(\cdot)$ are independent. 
By Theorem 2 in \textcolor{blue}{Fang et al. (2022)}, we  obtain that 
 $\max_{j,k\in[p]}\|(n-1)^{-1}\sum_{i=1}^n\tilde \varepsilon_{i,jk}(\bY_{i,-j,-k}-\widebar\bY_{-j,-k})\|_{\cS,\max} = O_{\rm p}(\sqrt{\log p/n})$. This together with Inequality~\ref{lemma_l1_inequality} and Proposition~\ref{lm_beta} yields that
\begin{align*}
&\max_{j,k\in[p]}\int_{\cU}\int_{\cU}\bigg\{\frac{1}{n-1}\sum_{i=1}^n{\rm II}_i(j,k,u,v)\bigg\}^2\,{\rm d}u{\rm d}v \\
&~~~~ \leq \max_{j,k\in[p]}\bigg\|\frac{1}{n-1}\sum_{i=1}^n\tilde{\varepsilon}_{i,jk}(\bY_{i,-j,-k}-\widebar\bY_{-j,-k})\bigg\|_{\cS,\max}^2 \max_{j,k\in[p]}\|{\boldsymbol{\beta}}_{jk}-\boldsymbol{\beta}_{jk}\|_{\cS,1}^2
=O_{\rm p}(\delta_{1n}^2  n^{-1}\log p) \,.
\end{align*}
Hence, \begin{equation}
\label{a1.rate}
\delta_{1n}\sqrt{\log p/n} \asymp n^{-a_{1}}(\log p)^{a_2}.
\end{equation}
Recall $\widetilde T_{n,jk}^{\varepsilon}= n\|\widetilde\Sigma_{jk}^{\varepsilon}\|_{\cS}^2$ and  \eqref{eq:sjksplit}. 
Following the similar arguments as the proof of Theorem \ref{thm:nonFDR} with $X_{ij}(\cdot)$ and $X_{ik}(\cdot)$ replaced by $\varepsilon_{i,jk}(\cdot)$ and $\varepsilon_{i,kj}(\cdot)$, respectively, 
we complete the proof of Theorem \ref{thm:noncondH0}. $\hfill\Box$

\subsection{Proof of Proposition \ref{lm_beta}}
The proof follows from Theorem 4 in \textcolor{blue}{Fang et al. (2022)} and Theorem 5 in \textcolor{blue}{Guo and Qiao (2023)}. Note that since $\bY_1(\cdot), \dots, \bY_n(\cdot)$ are i.i.d. sampled from a multivariate Gaussian process, the functional stability measure takes the value of $1$ and hence does not show up in the final rate.
$\hfill\Box$

\subsection{Proof of Proposition \ref{lm_betasigma}}
By Lemma~\ref{lm_guo_th}, we have that
\begin{align*} 
    &\max_{j,k\in[p]} \bigg\|\int_\cU \widehat{\bXi}_{-j,-k}(\cdot,s) (\widehat{\bbeta}_{kj}-\bbeta_{kj})(\cdot,s)\,{\rm d}s\bigg\|_{\cS,\max} \notag
    \\&~~~\leq \max_{j,k\in[p]} \bigg\|\int_\cU \{\widehat{\bXi}_{-j,-k}(\cdot,s)\widehat{\bbeta}_{kj}(\cdot,s)- {\bXi}_{-j,-k}(\cdot,s)\bbeta_{kj}(\cdot,s)\}\,{\rm d}s\bigg\|_{\cS,\max} 
    \\&~~~~~~~~~+ \max_{j,k\in[p]} \|\widehat{\bXi}_{-j,-k} -{\bXi}_{-j,-k}\|_{\cS,\max} \max_{j,k\in[p]} \|\bbeta_{kj}\|_{\cS,1} \notag
    \\&~~~\lesssim \max_{j,k\in[p]} \bigg\|\int_\cU \{\widehat{\bXi}_{-j,-k}(\cdot,s)\widehat{\bbeta}_{kj}(\cdot,s)- {\bXi}_{-j,-k}(\cdot,s)\bbeta_{kj}(\cdot,s)\}\,{\rm d}s\bigg\|_{\cS,\max} + O_p(s \sqrt{{\log p}/{n}}).
\end{align*}
Recall that  $\hat \beta_{kj,\ell}(v,s) = \sum_{m_1,m_2\in [d]} \hat b_{kj,\ell m_1m_2} \hat \phi_{km_1}(v) \hat \phi_{\ell m_2}(s)$ and for $\ell\in S_{kj},$  $\beta_{kj,\ell}(v,s) = \sum_{m_1,m_2=1}^\infty b_{kj,\ell m_1m_2} \phi_{km_1}(v) \phi_{\ell m_2}(s).$ Rewriting this problem under a FPCA framework leads to
\begin{align*}
    ~~~&\int_\cU \{\widehat \Xi^{(-j,-k)}_{\ell' \ell}(t,s)\hat \beta_{kj,\ell}(v,s) -   \Xi^{(-j,-k)}_{\ell' \ell} (t,s)\beta_{kj,\ell}(v,s)  \} \,{\rm d}s
    \\&~~~= \sum_{m_1,m_2 \in [d]} \bigg[ \hat b_{kj,\ell m_1m_2} \hat \phi_{km_1}(v) \bigg\{\sum_{m_2'=1}^{\infty} \bigg(\frac{\sum_{i=1}^n \xi_{i \ell' m_2'} \hat \xi_{i \ell m_2}}{n}\bigg) \phi_{\ell' m_2'}(t)\bigg\}
    \\&~~~~~~~~~~~~- b_{kj,\ell m_1m_2}  \phi_{km_1}(v) \bigg\{\sum_{m_2'=1}^{\infty} \mathbb{E}(\xi_{i \ell' m_2'}  \xi_{i \ell m_2}) \phi_{\ell' m_2'}(t)\bigg\}\bigg]
    \\&~~~~~~-\underbrace{ \sum_{m_1 = d+1}^{\infty}\sum_{m_2 \in [d]}  b_{kj,\ell m_1m_2}  \phi_{km_1}(v) \bigg\{\sum_{m_2'=1}^{\infty} \mathbb{E}(\xi_{i \ell' m_2'}  \xi_{i \ell m_2}) \phi_{\ell' m_2'}(t)\bigg\}}_{R_{1, kj\ell'\ell}(t,v)}
    \\&~~~~~~- \underbrace{\sum_{m_1 \in [d]} \sum_{m_2 = d+1}^{\infty}  b_{kj,\ell m_1m_2}  \phi_{km_1}(v) \bigg\{\sum_{m_2'=1}^{\infty} \mathbb{E}(\xi_{i \ell' m_2'}  \xi_{i \ell m_2}) \phi_{\ell' m_2'}(t)\bigg\}}_{R_{2, kj\ell'\ell}(t,v)}
    \\&~~~~~~-\underbrace{\sum_{m_1 = d+1}^{\infty} \sum_{m_2 = d+1}^{\infty}  b_{kj,\ell m_1m_2}  \phi_{km_1}(v) \bigg\{\sum_{m_2'=1}^{\infty} \mathbb{E}(\xi_{i \ell' m_2'}  \xi_{i \ell m_2}) \phi_{\ell' m_2'}(t)\bigg\}}_{R_{3, kj\ell'\ell}(t,v)}
     \\&= \underbrace{\sum_{m_1,m_2 \in [d]} \bigg[ (\hat b_{kj,\ell m_1m_2}-b_{kj,\ell m_1m_2}) \hat \phi_{km_1}(v) \bigg\{\sum_{m_2'=1}^{\infty} \bigg(\frac{\sum_{i=1}^n \xi_{i \ell' m_2'} \hat \xi_{i \ell m_2}}{n}\bigg) \phi_{\ell' m_2'}(t)\bigg\}\bigg]}_{I_{1,kj\ell'\ell}(t,v)}
     \\&~~~+\underbrace{\sum_{m_1,m_2 \in [d]} \bigg(  b_{kj,\ell m_1m_2} \hat \phi_{km_1}(v) \bigg[\sum_{m_2'=1}^{\infty} \phi_{\ell' m_2'}(t) \bigg\{\frac{\sum_{i=1}^n \xi_{i \ell' m_2'} \hat \xi_{i \ell m_2}}{n} -  \mathbb{E}(\xi_{i \ell' m_2'}  \xi_{i \ell m_2})\bigg\} \bigg]\bigg)}_{I_{2,kj\ell'\ell}(t,v)}
     \\&~~~+\underbrace{\sum_{m_1,m_2 \in [d]} \bigg[b_{kj,\ell m_1m_2} \{\hat  \phi_{km_1}(v) -\phi_{km_1}(v) \}\bigg\{\sum_{m_2'=1}^{\infty} \mathbb{E}(\xi_{i \ell' m_2'}  \xi_{i \ell m_2}) \phi_{\ell' m_2'}(t)\bigg\}\bigg]}_{I_{3,kj\ell'\ell}(t,v)} - {R_{kj\ell'\ell}(t,v)},
\end{align*}
where $R_{kj\ell'\ell}(t,v) = \sum_{r\in [3]}R_{r,kj\ell'\ell}(t,v)$. Hence,  we obtain that
\begin{align}
    &\max_{j,k\in[p]} \bigg\|\int_\cU \{\widehat{\bXi}_{-j,-k}(\cdot,s)\widehat{\bbeta}_{kj}(\cdot,s)- {\bXi}_{-j,-k}(\cdot,s)\bbeta_{kj}(\cdot,s)\}\,{\rm d}s\bigg\|_{\cS,\max} \notag
    \\&~~~= \max_{j,k\in[p]} \max_{\ell' \in [p-2]} \bigg\|\sum_{\ell \in [p-2]} \int_\cU \{\widehat \Xi^{(-j,-k)}_{\ell' \ell}(\cdot,s)\hat \beta_{kj,\ell}(\cdot,s) -   \Xi^{(-j,-k)}_{\ell' \ell} (\cdot,s)\beta_{kj,\ell}(\cdot,s)  \} \,{\rm d}s\bigg\|_\cS \notag
    \\&~~~\leq \max_{j,k\in[p]} \max_{\ell' \in [p-2]} \big\|\sum_{\ell \in [p-2]}I_{1,kj\ell'\ell} \big\|_\cS + s \cdot \max_{j,k\in[p]} \max_{(\ell',\ell) \in [p-2]^2} \|I_{2,kj\ell'\ell}\|_\cS  \notag
    \\&~~~~~~~~~+ s \cdot \max_{j,k\in[p]} \max_{(\ell',\ell) \in [p-2]^2} \|I_{3,kj\ell'\ell}\|_\cS +s \cdot \max_{j,k\in[p]} \max_{(\ell',\ell) \in [p-2]^2} \|R_{kj\ell'\ell}\|_\cS. \notag
\end{align}
Let $\omega_0 = \max_{ \ell \in [p] }\int_\cU \Xi_{\ell\ell}(u,u){\rm d}u.$ By Cauchy--Schwarz inequality, the orthonormality of $\phi_{jm}(\cdot)$ for $j \in [p]$ and $m \in [d]$ and Condition~\ref{con_fof_beta_a}, we obtain that
\begin{align*}
    &\max_{j,k\in[p]} \max_{(\ell',\ell) \in [p-2]^2} \|R_{1,kj\ell'\ell}\|_\cS^2 \\&~~~~~~ = \max_{j,k\in[p]} \max_{(\ell',\ell) \in [p-2]^2}  \bigg\| \sum_{m_1 = d+1}^{\infty}\sum_{m_2 \in [d]}  b_{kj,\ell m_1m_2}  \phi_{km_1} \bigg\{\sum_{m_2'=1}^{\infty} \mathbb{E}(\xi_{i \ell' m_2'}  \xi_{i \ell m_2}) \phi_{\ell' m_2'}\bigg\}\bigg\|_\cS^2
    \\&~~~~~~ = \max_{j,k\in[p]} \max_{(\ell',\ell) \in [p-2]^2}\sum_{m_2 \in [d]}\sum_{m_2'' \in [d]}\sum_{m_1 =d+1}^{\infty} b_{kj,\ell m_1m_2}b_{kj,\ell m_1m_2''} \sum_{m_2'=1}^{\infty}  \mathbb{E}(\xi_{i \ell' m_2'}  \xi_{i \ell m_2}) \mathbb{E}(\xi_{i \ell' m_2'}  \xi_{i \ell m_2''})
    \\&~~~~~~ \leq \max_{j,k\in[p]} \max_{(\ell',\ell) \in [p-2]^2}\sum_{m_2 \in [d]}\sum_{m_2'' \in [d]}\sum_{m_1 =d+1}^{\infty} b_{kj,\ell m_1m_2}b_{kj,\ell m_1m_2''} \sum_{m_2'=1}^{\infty}   \mathbb{E}(\xi_{i \ell' m_2'} ^2  )\sqrt{\mathbb{E}(\xi_{i \ell' m_2} ^2)  \mathbb{E}(\xi_{i \ell' m_2''} ^2)}
    \\&~~~~~~ \leq \omega_0 \max_{j,k\in[p]} \max_{(\ell',\ell) \in [p-2]^2} \sqrt{\sum_{m_2 \in [d]}\sum_{m_2'' \in [d]} \bigg(\sum_{m_1 =d+1}^{\infty} b_{kj,\ell m_1m_2}b_{kj,\ell m_1m_2''}\bigg)^2 }  \sqrt{\sum_{m_2 \in [d]}\sum_{m_2'' \in [d]} \mathbb{E}(\xi_{i \ell' m_2}^2)  \mathbb{E}(\xi_{i \ell' m_2''}^2)}
    \\&~~~~~~ \lesssim \max_{\ell \in [p-2]} \max_{j,k\in[p]} \sum_{m_2 \in [d]}\sum_{m_1 =d+1}^{\infty}b_{kj,\ell m_1m_2}^2 \lesssim \sum_{m_2 \in [d]}\sum_{m_1 =d+1}^{\infty} (m_1+m_2)^{-2\upsilon-1} = O(d^{-2\upsilon+1}).
\end{align*}
Following the same argument, we can also have  $\max_{j,k\in[p]} \max_{(\ell',\ell) \in [p-2]^2} \|R_{kj\ell'\ell}\|_\cS = O(d^{-\upsilon+1/2}).$
We next focus on the terms $I_{r,kj\ell' \ell}$ for $r \in [3]$. By the orthonormality of $\phi_{jm}(\cdot)$ for $j \in [p]$ and $m \in [d]$ and Condition~\ref{con_fof_beta_a},
\begin{align*}
    &\max_{j,k\in[p]} \max_{(\ell',\ell) \in [p-2]^2} \|I_{3,kj\ell' \ell}\|_\cS^2 \\&~~~~~~= \max_{j,k\in[p]} \max_{(\ell',\ell) \in [p-2]^2}  \bigg\|\sum_{m_1,m_2 \in [d]} \bigg[b_{kj,\ell m_1m_2} \{\hat  \phi_{km_1} -\phi_{km_1} \}\bigg\{\sum_{m_2'=1}^{\infty} \mathbb{E}(\xi_{i \ell' m_2'}  \xi_{i \ell m_2}) \phi_{\ell' m_2'}\bigg\}\bigg]\bigg\|_\cS^2
    \\&~~~~~~ = \max_{j,k\in[p]} \max_{(\ell',\ell) \in [p-2]^2}\sum_{m_2'=1}^{\infty} \int_\cU\bigg[\sum_{m_1 \in [d]}\{\hat  \phi_{km_1}(v) -\phi_{km_1}(v) \} \cdot \bigg\{\sum_{m_2 \in [d]}\mathbb{E}(\xi_{i \ell' m_2'}  \xi_{i \ell m_2}) b_{kj,\ell m_1m_2}
      \bigg\}\bigg]^2{\rm d} v
    \\&~~~~~~ \leq  d \max_{k \in [p], m_1 \in [d]}\|\hat \phi_{km_1}-\phi_{km_1}\|^2 \cdot \max_{j,k\in[p]} \max_{(\ell',\ell) \in [p-2]^2} \sum_{m_2'=1}^{\infty} \sum_{m_1 \in [d]}\bigg\{\sum_{m_2 \in [d]}\mathbb{E}(\xi_{i \ell' m_2'}  \xi_{i \ell m_2}) b_{kj,\ell m_1m_2}
      \bigg\}^2
      \\&~~~~~~ \lesssim d \max_{k \in [p],m_1 \in [d]}\|\hat \phi_{km_1}-\phi_{km_1}\|^2.
\end{align*}
This, together with Lemma~\ref{lm_guo_th}, implies that $\max_{(\ell',\ell) \in [p-2]^2} \|I_{3,kj\ell' \ell}\|_\cS = O_p(d^{\varpi+3/2}\sqrt{{\log (pd)}/{n}}).$ Recall that
$X_{i\ell}(\cdot)=\sum_{m=1}^{\infty}\xi_{i\ell m}\phi_{\ell m}(\cdot),$ where the coefficients $\xi_{i\ell m} = \langle X_{i\ell}, \phi_{\ell m} \rangle,$ and the corresponding estimated FPC scores $\hat \xi_{i\ell m} = \langle X_{i\ell}, \hat \phi_{\ell m} \rangle.$ By the orthonormality of $\hat \phi_{jm}(\cdot)$ for $j \in [p]$ and $m \in [d],$  Condition~\ref{con_fof_beta_a} and the fact that $\|\Xi_{\ell'\ell}\|_{\cS} \leq \omega_0,$ we obtain that
\begin{align*}
    &\max_{j,k\in[p]} \max_{(\ell',\ell) \in [p-2]^2} \|I_{2,kj\ell' \ell}\|_\cS^2
    \\&~~~~~~= \max_{j,k\in[p]} \max_{(\ell',\ell) \in [p-2]^2}  \bigg\|\sum_{m_1,m_2 \in [d]} \bigg(  b_{kj,\ell m_1m_2} \hat \phi_{km_1} \bigg[\sum_{m_2'=1}^{\infty} \phi_{\ell' m_2'} \bigg\{\frac{\sum_{i=1}^n \xi_{i \ell' m_2'} \hat \xi_{i \ell m_2}}{n} -  \mathbb{E}(\xi_{i \ell' m_2'}  \xi_{i \ell m_2})\bigg\} \bigg]\bigg)\bigg\|_\cS^2
     \\&~~~~~~ = \max_{j,k\in[p]} \max_{(\ell',\ell) \in [p-2]^2}\sum_{m_1 \in [d]}\bigg\|\sum_{m_2 \in [d]} \bigg(   b_{kj,\ell m_1m_2} \cdot \bigg[\int_\cU \big\{\widehat \Xi_{\ell' \ell}(\cdot,s) \hat \phi_{\ell m_2}(s)  -  \Xi_{\ell' \ell}(\cdot,s)  \phi_{\ell m_2}(s)\big\} {\rm d} s \bigg]\bigg)\bigg\|^2
     \\&~~~~~~ \lesssim \max_{(\ell',\ell) \in [p-2]^2} \sum_{m_2 \in [d]} \int_\cU \bigg[ \int_\cU \big\{\widehat \Xi_{\ell' \ell}(t,s) \hat \phi_{\ell m_2}(s)  -  \Xi_{\ell' \ell}(t,s)  \phi_{\ell m_2}(s)\big\} {\rm d} s\bigg]^2
     {\rm d} t
    \\&~~~~~~ \lesssim  d \max_{(\ell',\ell) \in [p-2]^2} \|\widehat \Xi_{\ell' \ell}-  \Xi_{\ell' \ell}\|_\cS^2 + d \max_{\ell \in [p-2],\,m_2 \in [d]}\|\hat \phi_{\ell m_2}-\phi_{\ell m_2}\|^2.
\end{align*}
Then by Lemma~\ref{lm_guo_th}, we have that $\max_{j,k\in[p]} \max_{(\ell',\ell) \in [p-2]^2} \|I_{2,kj\ell' \ell}\|_\cS = O_p( d^{\varpi+3/2}\sqrt{{\log(pd)}/{n}}).$
We now focus on the term $I_{1,kj\ell' \ell}(t,v)$.
We further decompose $I_{1,kj\ell' \ell}(t,v)$ as
\begin{align*}
     I_{1,kj\ell' \ell}(t,v) &= \underbrace{\sum_{m_1,m_2 \in [d]}  \bigg[ (\hat b_{kj,\ell m_1m_2}-b_{kj,\ell m_1m_2}) \hat \phi_{km_1}(v) \bigg\{\sum_{m_2'\in[d]} \bigg(\frac{\sum_{i=1}^n \xi_{i\ell' m_2'} \hat \xi_{i\ell m_2}}{n}\bigg) \phi_{\ell' m_2'}(t)\bigg\}\bigg]}_{I_{1,kj\ell' \ell}^{(1)}(t,v)}
     \\&~~
     +\underbrace{\sum_{m_1,m_2 \in [d]}\bigg[ (\hat b_{kj,\ell m_1m_2}-b_{kj,\ell m_1m_2}) \hat \phi_{km_1}(v) \bigg\{\sum_{m_2' = d+1}^{\infty} \bigg(\frac{\sum_{i=1}^n \xi_{i\ell' m_2'} \hat \xi_{i\ell m_2}}{n}\bigg) \phi_{\ell' m_2'}(t)\bigg\}\bigg]}_{I_{1,kj\ell' \ell}^{(2)}(t,v)}.
\end{align*}

We next introduce some matrix forms for the convenience of theoretical development regarding $I_{1,jk\ell'\ell}^{(1)}$. For a block matrix $\bB = (\bB_{jk})_{j,k \in [p]} \in \eR^{pd \times pd}$ with the $(j,k)$th block $\bB_{jk} \in {\eR}^{d \times d},$  we define its $d$-block versions of elementwise $\ell_{\infty}$ and matrix $\ell_1$ norms by $\|\bB\|_{\max}^{(d)} = \max_{j,k} \|\bB_{jk}\|_{\tF}$ and
$
\|\bB\|_{1}^{(d)} = {\max}_{k}\sum_{j} \|\bB_{jk}\|_{\tF},
$
respectively.
Define $\bV_{j}= (\bxi_{1j},\dots,\bxi_{nj})^{\T}\in \mathbb{R}^{n\times d},$ $\bV = (\bV_{1}, \dots, \bV_{p})\in \mathbb{R}^{n \times pd}$, 
and $\bPsi_{jk} = (\bPsi_{jk,\ell}^{\T})_{\ell \in [p] \setminus \{j,k\}}^{\T}\in\mathbb{R}^{(p-2)d\times d}$.
We can rewrite \eqref{eq_rg} as
\begin{align} \label{supp_matrix}
    \bV_j = \bV_{-j,-k}\bPsi_{jk} + \bR_{jk} +\bE_{jk},
\end{align}
where $\bV_{-j,-k}$ denotes the submatrix of $\bV$ by removing the $j$th and $k$th blocks, namely $\bV_{j}$ and $\bV_{k}$.  Note $\bR_{jk} $ and $\bE_{jk} $ are $n\times d$ matrices whose row vectors are formed by truncation errors $\{\br_{i,jk}\in\mathbb{R}^{d}: i \in [n]\}$ and random errors $\{\beps_{i,jk}\in\mathbb{R}^{d}: i \in [n]\}$, where
$\br_{i,jk} = (r_{i1,jk},\dots,r_{id,jk})^{\T}$ with $r_{im_1,jk} = \sum_{\ell \in [p] \setminus\{j,k\} }\sum_{m_2=d+1}^\infty\langle\langle\phi_{jm_1},\beta_{jk,\ell}\rangle,\phi_{\ell m_2}\rangle\xi_{i\ell m_2}$ and $\beps_{i,jk} = (\epsilon_{i1,jk},\dots,\epsilon_{id,jk})^{\T}$
with $\epsilon_{im_1,jk} = \langle\phi_{jm_1},\varepsilon_{i,jk}\rangle$ for $m_1 \in [d]$, respectively.

Let $\widehat{\bQ}_{\ell} = \{(n)^{-1}\widehat{\bV}_{\ell}^{\T}\widehat{\bV}_{\ell}\}^{1/2}\in\mathbb{R}^{d\times d}$ with  $\widehat \bV_{\ell}= (\widehat \bxi_{1\ell},\dots,\widehat \bxi_{n\ell})^{\T}\in \mathbb{R}^{n\times d}$ 
for $\ell \in [p].$  
Define $\widehat \bV = (\widehat \bV_{1}, \dots, \widehat \bV_{p})\in \mathbb{R}^{n\times pd}$ and  $\widehat{\bD} =  \text{diag}(\widehat{\bQ}_{1},\dots,\widehat{\bQ}_{p})\in \mathbb{R}^{pd\times pd}$.  Let $\widehat \bV_{-j,-k}$ denotes the submatrix of $\widehat \bV$ by removing the $j$th and $k$th blocks, namely $\widehat\bV_{j}$ and $\widehat\bV_{k}$ and $\widehat{\bD}_{-j,-k}$ denotes the submatrix of $\widehat{\bD}$ by removing the $j$th and $k$th diagonal blocks, namly $\widehat{\bQ}_{j}$ and $\widehat{\bQ}_{k}$. 
$\{\widehat \bPsi_{jk,\ell}\}_{\ell \in [p] \setminus \{j,k\}}$ can be equivalently obtained by computing
$\widehat{\bPsi}_{jk} = \widehat{\bD}_{-j,-k}^{-1}\widehat{\bB}_{jk} = (\widehat{\bPsi}_{jk,\ell})_{\ell \in [p] \setminus \{j,k\}},$  where
\begin{equation}
\nonumber
 \widehat \bB_{jk} = \underset{\bB _{jk}\in \mathbb{R}^{(p-2)d\times d}}{ \text{arg min}} \left\{ \frac{1}{2n}\|\widehat{\bV}_j-\widehat{\bV}_{-j,-k}\widehat{\bD}_{-j,-k}^{-1}\bB_{jk}\|_\tF^2 + \tau_n\|\bB_{jk}\|^{(d)}_1 \right\}.
\end{equation} 

Let ${\bD} =  \text{diag}({\bQ}_{1},\dots,{\bQ}_{p})\in \mathbb{R}^{pd\times pd}$ with ${\bQ}_{j} = \{(n)^{-1}{\bV}_{j}^{\T}{\bV}_{j}\}^{1/2}\in\mathbb{R}^{d\times d}$ for $j \in [p]$ and  ${\bB}_{jk}  = {\bD}_{-j,-k}{\bPsi}_{jk},$ where ${\bD}_{-j,-k}$ denotes the submatrix of ${\bD}$ by removing the $j$th and $k$th diagonal blocks, namely ${\bQ}_{j}$ and ${\bQ}_{k}$. 
Notice that
\begin{align*}
     &\max_{\ell' \in [p-2]} \big\|\sum_{\ell \in [p-2]}I_{1,kj\ell' \ell}^{(1)} \big\|_\cS = n^{-1}\|\bD_{-k,-j}\bD_{-k,-j}^{-1} \bV_{-k,-j}^{\T}\widehat\bV_{-k,-j}(\widehat \bPsi_{kj} - \bPsi_{kj})\|_{\max}^{(d)}
     \\&~~~~~~\leq n^{-1}\|\bD_{-k,-j}\|_{\max}\|\bD_{-k,-j}^{-1} \bV_{-k,-j}^{\T}\widehat\bV_{-k,-j}(\widehat \bPsi_{kj} - \bPsi_{kj})\|_{\max}^{(d)}
     \\&~~~~~~\lesssim n^{-1}\|\widehat\bD_{-k,-j}^{-1} \widehat\bV_{-k,-j}^{\T}\widehat\bV_{-k,-j} \widehat \bD_{-k,-j}^{-1}(\widehat \bB_{kj} - \bB_{kj})\|_{\max}^{(d)}
     \\&~~~~~~~~~+  \|\bD_{-k,-j}^{-1} \bV_{-k,-j}^{\T}\widehat\bV_{-k,-j} \widehat \bD_{-k,-j}^{-1}\|_{\max}^{(d)}\|\bD_{-k,-j}-\widehat \bD_{-k,-j}\|_{\max}
 \|\bPsi_{kj}\|_{1}^{(d)}
     \\&~~~~~~~~~+ n^{-1}\|\bD_{-k,-j}^{-1} \bV_{-k,-j}^{\T}\widehat\bV_{-k,-j} \widehat \bD_{-k,-j}^{-1}-\widehat\bD_{-k,-j}^{-1} \widehat\bV_{-k,-j}^{\T}\widehat\bV_{-k,-j} \widehat \bD_{-k,-j}^{-1}\|_{\max}^{(d)}\|\widehat \bB_{kj} - \bB_{kj}\|_{1}^{(d)}.
\end{align*}
By the KKT condition, we have $\max_{j,k\in[p]} n^{-1}\|\widehat{\bD}_{-k,-j}^{-1}\widehat{\bV}_{-k,-j}^{\T}(\widehat{\bV}_k-\widehat{\bV}_{-k,-j}\widehat{\bD}_{-k,-j}^{-1} \widehat \bB_{kj})\|_{\max}^{(d)} \leq \tau_n$. This, together with Lemma~\ref{lm_fang_pr7} implies that $\max_{j,k\in[p]} n^{-1}\|\widehat\bD_{-k,-j}^{-1} \widehat\bV_{-k,-j}^{\T}\widehat\bV_{-k,-j} \widehat \bD_{-k,-j}^{-1}(\widehat \bB_{kj} - \bB_{kj})\|_{\max}^{(d)} =O_p(\tau_n)$. Then, by Lemmas~\ref{lm_guo_th}--\ref{lm_fang_lm17} and the fact that $\max_{j,k\in[p]} \|\bPsi_{kj}\|_{1}^{(d)} = O(s)$, we obtain that $\max_{j,k\in[p]} \max_{\ell' \in [p-2]} \|\sum_{\ell \in [p-2]}I_{1,kj\ell' \ell}^{(1)} \|_\cS = O_p(\tau_n)$. Furthermore, by Cauchy--Schwarz inequality,
\begin{align*}
    &\sum_{m_2' = d+1}^{\infty} \bigg(\frac{\sum_{i=1}^n \xi_{i\ell' m_2'} \hat \xi_{i\ell m_2}}{n}\bigg) \phi_{\ell' m_2'}(t) = \frac{\sum_{i=1}^n \{\hat \xi_{i\ell m_2} (\sum_{m_2' = d+1}^{\infty}\xi_{i\ell' m_2'}\phi_{\ell' m_2'}(t)\} }{n}
    \\&~~~~~~ \leq \sqrt{\hat \omega_{\ell m_2}}\sqrt{\frac{\sum_{i=1}^n\{\sum_{m_2' = d+1}^{\infty}\xi_{i\ell' m_2'}\phi_{\ell' m_2'}(t)\}^2}{n}}.
\end{align*}
Combining this with  the orthonormality of $\hat \phi_{jm}(\cdot)$ for $j \in [p]$ and $m \in [d]$ thus leads to
\begin{align*}
    &\max_{j,k\in[p]} \max_{\ell' \in [p-2]} \big\| I_{1,kj\ell' \ell}^{(2)} \big\|_\cS^2
    \\&~~~~~~ =\max_{j,k\in[p]}  \max_{\ell' \in [p-2]}  \bigg\|  \sum_{m_1, m_2 \in [d]} \bigg[ (\hat b_{kj,\ell m_1m_2}-b_{kj,\ell m_1m_2}) \hat \phi_{km_1} \bigg\{\sum_{m_2'=d+1}^{\infty} \bigg(\frac{\sum_{i=1}^n \xi_{i\ell' m_2'} \hat \xi_{i\ell m_2}}{n}\bigg) \phi_{\ell' m_2'}\bigg\}\bigg] \bigg\|_\cS^2
     \\&~~~~~~ \leq \max_{j,k\in[p]} \max_{\ell' \in [p-2]}\sum_{m_1 \in [d]}\bigg\| \sum_{m_2 \in [d]} \bigg\{  (\hat b_{kj,\ell m_1m_2}-b_{kj,\ell m_1m_2})  \cdot  \sqrt{\hat \omega_{\ell m_2}}\sqrt{\frac{\sum_{i=1}^n\{\sum_{m_2' = d+1}^{\infty}\xi_{i\ell' m_2'}\phi_{\ell' m_2'}(t)\}^2}{n}}\bigg\}\bigg\|^2
     \\&~~~~~~ \leq d \max_{j,k\in[p]} \sum_{m_1,m_2 \in [d]}\bigg\{   (\hat b_{kj,\ell m_1m_2}-b_{kj,\ell m_1m_2})\sqrt{\hat\omega_{\ell m_2}} \bigg\}^2\cdot \max_{\ell' \in [p-2]}\frac{\sum_{m_2' = d+1}^{\infty}\sum_{i=1}^n\xi_{i\ell' m_2'}^2}{n}
\end{align*}
By the fact that ${n^{-1}\sum_{i=1}^n\xi_{i\ell' m_2'}^2} - \mathbb{E}(\xi_{i\ell' m_2'}^2) =  \langle \phi_{\ell' m_2'}, \langle \widehat \Xi_{\ell' \ell'} -\Xi_{\ell' \ell'}, \phi_{\ell' m_2'}\rangle\rangle \leq  \|\widehat \Xi_{\ell' \ell'} -\Xi_{\ell' \ell'}\|_\cS = O_p\{ n^{-1/2} \}$,  $\sum_{m_2' = d+1}^{\infty}\mathbb{E}(\xi_{i\ell' m_2'}^2) = O(d^{-\varpi+1})$ and Lemma~\ref{lm_fang_c1},
we thus have $$\max_{j,k\in[p]} \max_{\ell' \in [p-2]} \bigg\|\sum_{\ell \in [p-2]}I_{1,kj\ell' \ell}^{(2)} \bigg\|_\cS = O_p( d^{-\varpi/2+1}s\tau_n/\underline{\mu}).$$
Combining the above results, we obtain that $$\max_{j,k\in[p]} \bigg\|\int_\cU \widehat{\bXi}_{-j,-k}(\cdot,s) (\widehat{\bbeta}_{kj}-\bbeta_{kj})(\cdot,s)\,{\rm d}s\bigg\|_{\cS,\max} =O_p(\tau_n) + O_p(d^{-\varpi/2+1}s\tau_n/\underline{\mu}).$$

\section{Additional technical proofs}
\label{supp.sec_ad}

\subsection{Proof of Lemma \ref{lm_norm}} \label{sec.pf.lm_norm}
For each $j \in [p]$, let 
$Z_{ij }(\cdot)=\sum_{m=1}^{\infty}\vartheta_{ij  m}\psi_{j  m}(\cdot),$ where the coefficients $\vartheta_{ij  m} = \langle Z_{ij }, \psi_{j  m} \rangle,$ namely FPC scores, are uncorrelated random variables with mean zero and $\cov(\vartheta_{ij  m}, \vartheta_{ij  m'}) = \varrho_{ j  m}I(m=m').$ Recall that $q \in \mathcal{H}_0$, we have $\|\Sigma_{j_qk_q}\|_{\cS} = 0$, implying that $\cov(\vartheta_{ij_qm}, \vartheta_{ik_qm'}) = 0$ for any $m,m'\geq 1$. We thus obtain that
\begin{align}
    &\Gamma_{j_qk_q}(u_1,v_1,u_2,v_2) \nonumber
    \\&~~~= \cov\{Z_{ij_q}(u_1)Z_{ik_q}(v_1),Z_{ij_q}(u_2)Z_{ik_q}(v_2)\} \nonumber
    \\&~~~=\cov\bigg\{\sum_{m_1=1}^{\infty}\sum_{m_2 = 1}^\infty\vartheta_{ij_qm_1} \vartheta_{ik_qm_2}\psi_{j_qm_1}(u_1) \psi_{k_qm_2}(v_1)  ,\sum_{m_3=1}^{\infty}\sum_{m_4 = 1}^\infty\vartheta_{ij_qm_3} \vartheta_{ik_qm_4}\psi_{j_qm_3}(u_2) \psi_{k_qm_4}(v_2) \bigg\}\nonumber
    \\& ~~~= \sum_{m_1=1}^{\infty}\sum_{m_2 = 1}^\infty\cov\{\vartheta_{ij_qm_1} \vartheta_{ik_qm_2}, \vartheta_{ij_qm_1} \vartheta_{ik_qm_2}\} \psi_{j_qm_1}(u_1) \psi_{k_qm_2}(v_1)  \psi_{j_qm_1}(u_2) \psi_{k_qm_2}(v_2)\nonumber
    \\& ~~~= \sum_{m_1=1}^{\infty}\sum_{m_2 = 1}^\infty \varrho_{j_qm_1}\varrho_{k_qm_2} \psi_{j_qm_1}(u_1) \psi_{k_qm_2}(v_1)  \psi_{j_qm_1}(u_2) \psi_{k_qm_2}(v_2), \label{eq_decom}
\end{align}
for $q \in \mathcal{H}_0$. Since $\iint\psi_{j_qm_1}(u) \psi_{k_qm_2}(v)\psi_{j_qm_3}(u) \psi_{k_qm_4}(v) {\rm d}u {\rm d}v = I( m_1 = m_3)I(m_2 =m_4)$, $\{ \psi_{j_qm_1} \otimes \psi_{k_qm_2}\}_{m_1,m_2 = 1}^{\infty}$ follows an orthonormal basis of $\mathbb{S}$. Hence, \eqref{eq_decom} gives the spectral decomposition of $\Gamma_{j_qk_q}(\cdot,\cdot,\cdot,\cdot)$.
 Then for each $q \in \mathcal{H}_0$, there always exist a bijection $\tilde \pi_q: \mathbb{N}_+ \times \mathbb{N}_+  \to \mathbb{N}_+$ such that $a_{i\tilde \pi_{q}(m_1,m_2)}^{(q)} = \vartheta_{ij_qm_1}\vartheta_{ik_qm_2},$ $\lambda_{\tilde \pi_{q}(m_1,m_2)}^{(q)} = \varrho_{j_qm_1}\varrho_{k_qm_2},$ and $\phi_{j_q k_q\tilde \pi_{q}(m_1,m_2)}^{(q)}(u,v) = \psi_{j_qm_1}(u)\psi_{k_qm_2}(v).$

This,  together with  Condition \ref{cond:subG} and Lemma 4 in \textcolor{blue}{Fang et al. (2022)}, implies that 
 there exist a constant $\tilde M$ such that
 $\mathbb{E}\{\exp(\tilde M \vartheta_{ijm}^2/{\varrho_{jm}})\} \leq 2$ for $j \in [p]$ and $m \geq 1$. By Young's inequality, there exists a constant $M$, for any $m_1, m_2 \geq 1$, $\mathbb{E}\exp[Ma_{i\tilde \pi_{q}(m_1,m_2)}^{(q)}/\{\lambda_{\tilde \pi_{q}(m_1,m_2)}^{(q)}\}^{1/2}] \leq  \mathbb{E}\{\exp(\tilde M \vartheta_{ij_qm_1}^2/{\varrho_{j_qm_1}}) + \exp(\tilde M \vartheta_{ik_qm_2}^2/{\varrho_{k_qm_2}})\} /2 \leq 2$. 
 Hence, the proof is completed.
$\hfill\Box$

\subsection{Proof of Lemma  \ref{lm_diff}}\label{sec.pf.lm_diff}
Recall $T_n^{(q)}=n\iint \{(n-1)^{-1}\sum_{i=1}^n z_{iq}(u,v)-n(n-1)^{-1}\bar Z_{j_q}(u)\bar Z_{k_q}(v)\}^2\,{\rm d}u{\rm d}v$ and $\mathring{T}_n^{(q)}=n\iint \{n^{-1}\sum_{i=1}^n z_{iq}(u,v)\}^2\,{\rm d}u{\rm d}v$.  Triangle inequality implies
\begin{align*}
\max_{q\in\cH_0}|T_n^{(q)}-\mathring{T}_n^{(q)}| 
\lesssim&~ n\max_{q\in\cH_0}\iint\bigg\{\frac{1}{n(n-1)}\sum_{i=1}^n z_{iq}(u,v)-\frac{n}{n-1}\bar Z_{j_q}(u)\bar Z_{k_q}(v)\bigg\}^2\,{\rm d}u{\rm d}v \\
&~ + n\max_{q\in\cH_0}\iint\bigg|\bigg\{\frac{1}{n}\sum_{i=1}^{n}z_{iq}(u,v)\bigg\}\bigg\{\frac{1}{n(n-1)}\sum_{i=1}^n z_{iq}(u,v)-\frac{n}{n-1}\bar Z_{j_q}(u)\bar Z_{k_q}(v)\bigg\}\bigg|\,{\rm d}u{\rm d}v\\
:=&~ {\rm I}+{\rm II} \,.
\end{align*}
To specify the convergence rate of ${\rm I}$ and ${\rm II}$, we first define two events
\begin{align*}
  &~\cE_1 = \bigg\{\max_{j,k\in[p]}\bigg\|\frac{1}{n}\sum_{i=1}^n (z_{ijk}-\mathbb{E}z_{ijk})\bigg\|_{\cS}\leq \tilde C_1 \sqrt{\frac{\log p}{n}}\bigg\} \,,\\
  &~ \cE_2 = \bigg\{\max_{j\in[p]}\bigg\|\frac{1}{n}\sum_{i=1}^n Z_{ij}\bigg\|\leq \tilde C_2\sqrt{\frac{\log p}{n}}\bigg\} \,,
\end{align*}
where $\tilde C_1$ and $\tilde C_2$ are two positive constants. Since $\mathbb{E}z_{iq}(u,v)=\mathbb{E}\{Z_{ij_q}(u)Z_{ik_q}(v)\}=0$ if $q\in\cH_0$, then restricted on $\cE_1\bigcap\cE_2$, we have
\begin{align*}
 {\rm I}\lesssim &~ \frac{1}{n}\max_{q\in\cH_0}\bigg\|\frac{1}{n}\sum_{i=1}^n z_{iq}\bigg\|_{\cS}^2 + n\max_{j\in[p]}\bigg\|\frac{1}{n}\sum_{i=1}^n Z_{ij}\bigg\|^4 \lesssim \frac{(\log p)^2}{n} \,,\\
 {\rm II}\lesssim &~ n\max_{q\in\cH_0} \bigg\|\frac{1}{n}\sum_{i=1}^n z_{iq}\bigg\|_{\cS}\sqrt{\frac{{\rm I}}{n}} \lesssim \frac{(\log p)^{3/2}}{n^{1/2}}\,,
\end{align*} 
provided that $p\lesssim n^\kappa$ for some constant $\kappa>0$.
Based on Equations (S2) and (S3) in the supplementary material of {\color{blue}  Zapata et al. (2022)}, it holds that $\mathbb{P}(\cE_1)\rightarrow 1$ and $\mathbb{P}(\cE_2)\rightarrow 1$ as $n\rightarrow\infty$. Then we can conclude that $\max_{q\in\cH_0}|T_n^{(q)}-\mathring{T}_n^{(q)}| =O_{\rm p}\{n^{-1/2}(\log p)^{3/2}\}$.
Thus we complete the proof of Lemma \ref{lm_diff}. $\hfill\Box$

\subsection{Proof of Lemma \ref{partial.vb}}\label{pf.lm_vb}
Throughout this section, we use $c_0$ to denote generic positive finite constants that may be different in different uses. Recall that $\be_0=(1,0)^{\T}$, $\widetilde{\bU}_{ijt}=\{1,(U_{ijt}-u)/h_j\}^{\T}$ and
\begin{align*}
 \widehat{\bS}_{ij}(u)=\frac{1}{T_{ij}}\sum_{t=1}^{T_{ij}} \widetilde{\bU}_{ijt}\widetilde{\bU}_{ijt}^{\T}K_{h_j}(U_{ijt}-u)\,, ~~~~
 \widehat{\bR}_{ij}(u)=\frac{1}{T_{ij}}\sum_{t=1}^{T_{ij}} \widetilde{\bU}_{ij}\{W_{ijt} - Z_{ij}(u)\}K_{h_j}(U_{ijt}-u) \,.
\end{align*}
A simple calculation yields $\widehat{X}_{ij}(u) - Z_{ij}(u)=\be_0^{\T}\{\widehat{\bS}_{ij}(u)\}^{-1}\widehat{\bR}_{ij}(u)$. We can then write $\widehat{X}_{ij}(u) - Z_{ij}(u) - \widetilde{X}_{ij}(u)$ as 
\begin{align*}
    \widehat{X}_{ij}(u) - Z_{ij}(u) - \widetilde{X}_{ij}(u) = &\be_0^{\T} [\bbE\{\widehat{\bS}_{ij}(u)\}]^{-1} [\widehat \bR_{ij}(u) - \bbE\{\widehat \bR_{ij}(u)\}]  \\ &-  \be_0^{\T} \{\widehat{\bS}_{ij}(u)\}^{-1} [\widehat \bS_{ij}(u) - \bbE\{\widehat \bS_{ij}(u)\}]  [\bbE\{\widehat{\bS}_{ij}(u)\}]^{-1} \widehat \bR_{ij}(u),
\end{align*}
leading to 
\begin{equation}
    \begin{split}
        &|\widehat{X}_{ij}(u) - Z_{ij}(u) - \widetilde{X}_{ij}(u)| \\ \leq & ~\|\bbE\{\widehat{\bS}_{ij}(u)\}\|_{\min}^{-1} \|\widehat \bR_{ij}(u) - \bbE\{\widehat \bR_{ij}(u)\}\|  \\ 
        &+  \|\widehat{\bS}_{ij}(u)\|_{\min}^{-1} \|\widehat \bS_{ij}(u) - \bbE\{\widehat \bS_{ij}(u)\}\|_{\tF}  \|\bbE\{\widehat{\bS}_{ij}(u)\}\|_{\min}^{-1} \|\widehat \bR_{ij}(u)\|,
    \end{split}
\end{equation}
where
$\|\bB\|_{\min} = \{\lambda_{\min}(\bB^{\T}\bB)\}^{1/2}$ for any matrix $\bB = (B_{jk})_{p \times p}$. 

For $m_1, m_2 = 1,2$, let $\widehat S_{ijm_1m_2}(u)$ be the $(m_1,m_2)$-th element of $\widehat \bS_{ij}(u)$. By Conditions~\ref{cond.Tij}--\ref{cond.kern}, we have that for $q = 2,3,\dots$ and $s = 0,1,2$, there exits a positive constant $c_0$ such that
\begin{align*}
    \bbE\bigg\{\Big | \Big( \frac{U_{ijt} - u}{h_j}\Big)^s K_{h_j}(U_{ijt}-u)\Big|^q\bigg\} \leq \int h_j^{-q} K^q\Big(\frac{t-u}{h_j} \Big) \Big|\frac{t-u}{h_j}\Big|^{sq} f_{U}(t) {\rm d}t \leq c_0 h^{1-q}.
\end{align*}
Then we obtain that
\begin{align*}
    \sum_{t=1}^{T_{ij}}\bbE\bigg\{\Big | \Big( \frac{U_{ijt} - u}{h_j}\Big)^s K_{h_j}(U_{ijt}-u)\Big|^2\bigg\} \leq c_0 T h^{-1}.
\end{align*}
\vspace{-1cm}
\begin{align*}
    \sum_{t=1}^{T_{ij}}\bbE\bigg\{\Big | \Big( \frac{U_{ijt} - u}{h_j}\Big)^s K_{h_j}(U_{ijt}-u)\Big|^q\bigg\} \leq 2^{-1} q! c_0 Th^{-1} h^{2-q} ~~~\text{for } q\geq 3.
\end{align*}
This, together with the Bernstein inequality (see also Theorem 2.10 and Corollary 2.11 of {\color{blue}Boucheron
et al. (2014)}), implies  that 
there exists some positive constant $\tilde c_3$ (independent of $u$) such that for any $\delta>0$ and $u \in \cU$,
\begin{align}\label{eq.S}
    \bbP\big\{ \|\widehat \bS_{ij}(u) - \bbE\{\widehat \bS_{ij}(u)\}\|_{\tF} \geq \delta \big\}\leq 8 \exp\bigg(-\frac{\tilde c_3 T h \delta^2}{1+\delta}\bigg).
\end{align}

For $m = 1,2$, let $\widehat R_{ijm}(u)$ be the $m$-th element of $\widehat{\bR}_{ij}(u)$. We will then show that there exists some positive constant $\tilde c_4$ (independent of $u$) such that for any $\delta>0$ and $u \in \cU$,
\begin{align} \label{eq.R}
    \bbP\big\{ |\widehat R_{ijm}(u) - \bbE\{\widehat R_{ijm}(u)\}|\geq \delta \big\}\leq \tilde c_4 \exp\bigg(-\frac{\tilde c_4 T h \delta^2}{1+\delta}\bigg).
\end{align}
We focus on the case of $m = 1$   and the proof of  $m = 2$ follows similarly. Let
\begin{align*}
    \widehat{R}_{ij11}(u) &= \frac{1}{T_{ij}}\sum_{t=1}^{T_{ij}} K_{h_j}(U_{ijt}-u)\mu_{j}(U_{ijt}), 
    \\
    \widehat{R}_{ij12}(u) &= \frac{1}{T_{ij}}\sum_{t=1}^{T_{ij}} K_{h_j}(U_{ijt}-u)\varsigma_{ijt}, 
    \\
    \widehat{R}_{ij13}(u) &=\frac{1}{T_{ij}}\sum_{t=1}^{T_{ij}} K_{h_j}(U_{ijt}-u)\{Z_{ij}(U_{ijt}) - Z_{ij}(u)\}.
\end{align*}
Thus, $\widehat{R}_{ij1}(u) - \bbE\{\widehat{R}_{ij1}(u)\}$ can be decomposed into
$$
 \widehat{R}_{ij1}(u) - \bbE\{\widehat{R}_{ij1}(u)\}=\widehat{R}_{ij11}(u) - \bbE\{\widehat{R}_{ij11}(u)\} +\widehat{R}_{ij12}(u) +\widehat{R}_{ij13}(u)\,.
$$
Similar to \eqref{eq.S},
we obtain that there exists some positive constant $\tilde c_5$ (independent of $u$) such that for any $\delta>0$ and $u \in \cU$,
\begin{align} \label{eq.R1}
\bbP\big\{ |\widehat R_{ij11}(u) - \bbE\{\widehat R_{ij11}(u)\}|\geq \delta \big\}\leq  2 \exp\bigg(-\frac{\tilde c_5 T h \delta^2}{1+\delta}\bigg).
\end{align}
Recall that $\varsigma_{ijt}$ is a sub-Guassian random variable  with
$
    \mathbb{E}
    \{\exp(\varsigma_{ijt} \tilde x )\} \leq \exp(2^{-1}\tilde c^2 \sigma_{j}^2\tilde x^2)
$
for all $ \tilde x \in \mathbb{R}$. By Proposition 2.5.2 of {\color{blue} Vershynin
(2018)}, for any positive integer $q$, $\bbE(|\varsigma_{ijt}|^q)\leq c_0^q \sigma_j^q q^{q/2}$. Then by the fact that $q^{q/2} \leq c_0 q!/2$, we have
\begin{align*}
    \bbE\bigg\{\Big | K_{h_j}(U_{ijt}-u)\varsigma_{ijt}\Big|^q\bigg\} \leq c_0^q h^{1-q}q!/2.
\end{align*}
Thus, there exists some positive constant $\tilde c_6$ (independent of $u$) such that for any $\delta>0$ and $u \in \cU$,
\begin{align}\label{eq.R2}
\bbP\big\{ |\widehat R_{ij12}(u)\}|\geq \delta \big\}\leq  2 \exp\bigg(-\frac{\tilde c_6 T h \delta^2}{1+\delta}\bigg).
\end{align}
Define the event $V_{j} = \{U_{ijt}, t\in [T_{ij}], i \in [n]\}$. Recall that $Z_{ij}(\cdot)$ is a mean-zero sub-Gaussian process.  Similar to the previous proof, by Condition \ref{cond.cont}, we have $\bbE\{ |Z_{ij}(U_{ijt})|^q|V_j\} \leq c_0^q q!/2$.
\begin{align*}
    &\bbE\bigg\{\Big | K_{h_j}(U_{ijt}-u)\{Z_{ij}(U_{ijt}) - Z_{ij}(u)\}\Big|^q\bigg\} \\\leq &2^{q-1}  \bbE\bigg\{\Big | K_{h_j}(U_{ijt}-u)Z_{ij}(U_{ijt}) \Big|^q\bigg\} + 2^{q-1}\bbE\bigg\{\Big | K_{h_j}(U_{ijt}-u) Z_{ij}(u)\Big|^q\bigg\}
    \\\leq &2^{q-1}  \bbE\bigg\{\Big | K_{h_j}(U_{ijt}-u)\Big|^q \bbE [Z_{ij}(U_{ijt})|V_{j}] \Big|^q\bigg\} + 2^{q-1}\bbE\bigg\{\Big | K_{h_j}(U_{ijt}-u) Z_{ij}(u)\Big|^q\bigg\} \\\leq & 
    c_0^q h^{1-q}q!/2.
\end{align*}
Then we have that there exists some positive constant $\tilde c_7$ (independent of $u$) such that for any $\delta>0$ and $u \in \cU$,
\begin{align*}
\bbP\big\{ |\widehat R_{ij13}(u)\}|\geq \delta \big\}\leq  2 \exp\bigg(-\frac{\tilde c_7 T h \delta^2}{1+\delta}\bigg).
\end{align*}
This, together with \eqref{eq.R1} and \eqref{eq.R2}, implies that \eqref{eq.R} holds. By applying the same techniques used in the proofs of the local and $\ell_2$ norm concentration inequalities in Sections A.1.2 and A.1.3 of {\color{blue}Guo et al. (2023)}, we can show that there exist some positive constants $c_3, c_4$ such that for any $\delta >0$,
   $$
\bbP(\|\widehat X_{ij} - Z_{ij} - \widetilde{X}_{ij}\| \geq \delta) \leq c_3 \exp (-c_4 Th \delta^2).
$$
For $M>0$ with $1-c_4M^2<0$, we select $\delta = \log(p \vee n )^{1/2} (Th)^{-1/2} M \leq 1$. Then
   $$
\bbP\bigg\{\frac{\max_{i \in [n], j \in[p]}\|\widehat X_{ij} - Z_{ij} - \widetilde{X}_{ij}\|}{\log(p \vee n )^{1/2} (Th)^{-1/2}} \leq  M\bigg\} \leq c_3 \exp \{(1-c_4M^2)  \log(p \vee n)\}.
$$
We hence complete the proof of (i) in Lemma \ref{partial.vb}. 

Recall that $\widetilde{X}_{ij}(u)=\be_0^{\T}[\bbE\{\widehat{\bS}_{ij}(u)\}]^{-1}\bbE\{\widehat{\bR}_{ij}(u)\}$. Thus,
$$\widetilde{X}_{ij}(u)-\mu_{j}(u) =  \be_0^{\T}[\bbE\{\widehat{\bS}_{ij}(u)\}]^{-1}\bbE[\widehat{\bR}_{ij}(u) - \widehat \bS_{ij}(u)\{\mu_j(u),0\}^{\T}].$$
By applying the same techniques used in the proof of Theorem 3 in {\color{blue}Guo et al. (2023)},  we complete the proof of Lemma \ref{partial.vb}. $\hfill\Box$

\subsection{Inequality \ref{lemma_l1_inequality} and its proof}

\begin{inequality} \label{lemma_l1_inequality}
Let  $\boldsymbol{\cF}=(\cF_1,\ldots,\cF_p)^{\T}$ with each $\cF_i \in \mathbb{S}$ and $\boldsymbol{\cG}=(\cG_1,\ldots,\cG_p)^{\T}$ with each $\cG_i \in \mathbb{S}$. Then  $\iint \{\int\boldsymbol{\cF}(u,s)^{\T}\boldsymbol{\cG}(v,s)\,{\rm d}s\}^2 {\rm d}u{\rm d}v  \leq \|\boldsymbol{\cF}\|_{\cS,\max}^2\|\boldsymbol{\cG}\|_{\cS,1}^2$.
\end{inequality}

\textit{Proof.}
By Cauchy--Schwarz inequality, we obtain that
\begin{align}
    &\iint \bigg\{\int\boldsymbol{\cF}(u,s)^{\T}\boldsymbol{\cG}(v,s)\,{\rm d}s\bigg\}^2 {\rm d}u{\rm d}v  \notag
    \\ & ~~~~~~= \iint \bigg\{\sum_{i=1}^p\int\cF_i(u,s)\cG_i(v,s)\,{\rm d}s\bigg\}^2 {\rm d}u{\rm d}v  \notag
    \\ & ~~~~~~= \iint {\sum_{i,j=1}^p\iint\cF_i(u,s_1)\cG_i(v,s_1)\cF_j(u,s_2)\cG_j(v,s_2)\,{\rm d}s_1 {\rm d}s_2} {\rm d}u{\rm d}v  \notag
    \\ & ~~~~~~= \sum_{i,j=1}^p \iint \bigg\{\int\cF_i(u,s_1)\cF_j(u,s_2)\,{\rm d}u\bigg\} \bigg\{\int \cG_i(v,s_1)  \cG_j(v,s_2)\,{\rm d}v\bigg\}\,{\rm d}s_1 {\rm d}s_2 \notag
    \\ & ~~~~~~\leq \sum_{i,j=1}^p \bigg[\iint \bigg\{\int\cF_i(u,s_1)\cF_j(u,s_2)\,{\rm d}u\bigg\}^2 \,{\rm d}s_1 {\rm d}s_2 \iint\bigg\{\int \cG_i(v,s_1)  \cG_j(v,s_2)\,{\rm d}v\bigg\}^2\,{\rm d}s_1 {\rm d}s_2 \bigg]^{1/2}\notag
    \\ & ~~~~~~\leq \sum_{i,j=1}^p \|\cF_i\|_{\cS}\|\cF_j\|_{\cS}\| \cG_i\|_{\cS}\| \cG_j\|_{\cS}\notag
    = \big(\sum_{i=1}^p\|\cF_i\|_{\cS}\| \cG_i\|_{\cS}\big)^2
    \\ & ~~~~~~\leq \big(\max_{i} \|\cF_i\|_{\cS} \sum_{i=1}^p\| \cG_i\|_{\cS}\big)^2 = \|\boldsymbol{\cF}\|_{\cS,\max}^2\|\boldsymbol{\cG}\|_{\cS,1}^2. \notag
\end{align}
Hence, the proof is complete.
 $\hfill\Box$


\subsection{Lemma \ref{lm_guo_th} and its proof}
\begin{lemma} \label{lm_guo_th}
Recall that  $\bY_i(\cdot)$ follows a mean zero multivariate Gaussian process with covariance function matrix $\bXi(u,v)=\cov\{\bY_i(u),\bY_i(v)\}$ for $(u,v) \in \cU^2.$
Then there exist some positive constants $c_1$ and $c_2$ such that, (i) for  $n \gtrsim \log p $,  the estimate $\widehat{\bXi}$ satisfies $ \|\widehat{\bXi} -{\bXi}\|_{\cS,\max} \lesssim \sqrt{\log p /n}$ with probability greater than $1-c_1(pd)^{-c_2};$ 
(ii) Suppose {\rm Condition~\ref{con_eigen}} hold. Then for   $n \gtrsim \log (pd)d^{4\varpi+2}$,  the estimates $\{\hat \phi_{\ell m}(\cdot)\}$ satisfies $ \max_{\ell \in [p],\,m \in [d]}\|\hat \phi_{\ell m}-\phi_{\ell m}\|^2 \lesssim d^{\varpi+1}\sqrt{\log(pd)/n}$,  with probability greater than $1-c_1(pd)^{-c_2}$.
\end{lemma}

\textit{Proof.} The proof follows from Theorems~2 and 3 in \textcolor{blue}{Guo and Qiao (2023)} and hence  is omitted here.
$\hfill\Box$

\subsection{Lemma \ref{lm_fang_c1} and its proof}
\begin{lemma} \label{lm_fang_c1}
Suppose {\rm Conditions~{\rm \ref{cond:subG}} and \ref{con_eigen}--\ref{con_fof_eigen_min}} hold. Then there exist some positive constants $c_1$ and $c_2$ such that, for  $n \gtrsim \log (pd)d^{4\varpi+2} $, with probability greater than $1-c_1(pd)^{-c_2},$ 
$
   \max_{j,k\in[p]} \|\widehat \bB_{kj} - \bB_{kj}\|_{1}^{(d)} \lesssim s\tau_n/\underline{\mu}.
$
\end{lemma}

\textit{Proof.} The proof follows similarly to that of  Theorem 4 in \textcolor{blue}{Fang et al. (2022)} and hence  is omitted here.
$\hfill\Box$

\subsection{Lemma \ref{lm_fang_lm17} and its proof}
\begin{lemma} \label{lm_fang_lm17}
Suppose {\rm Conditions~{\rm \ref{cond:subG}} and \ref{con_eigen}--\ref{con_fof_eigen_min}} hold. Then there exist some positive constants $c_1$ and $c_2$ such that, for  $n \gtrsim \log (pd)d^{4\varpi+2} $, with probability greater than $1-c_1(pd)^{-c_2},$ 
$
n^{-1}\|\bD^{-1} \bV^{\T}\widehat\bV \widehat \bD^{-1}-\widehat\bD^{-1} \widehat\bV^{\T}\widehat\bV \widehat \bD^{-1}\|_{\max} \lesssim d^{\varpi+1}\sqrt{\log(pd)/n}.
$
\end{lemma}

\textit{Proof.} The proof follows similarly to that of  Lemma 17 in \textcolor{blue}{Fang et al. (2022)} and hence  is omitted here.
$\hfill\Box$

\subsection{Lemma \ref{lm_fang_pr7} and its proof}
\begin{lemma} \label{lm_fang_pr7}
Suppose {\rm Conditions~{\rm \ref{cond:subG}} and \ref{con_eigen}--\ref{con_fof_eigen_min}} hold. Then there exist some positive constants $c_1$ and $c_2$ such that, for  $n \gtrsim \log (pd)d^{4\varpi+2} $, with probability greater than $1-c_1(pd)^{-c_2},$ 
$\max_{j,k\in[p]} 
 n^{-1}\|\widehat{\bD}_{-k,-j}^{-1}\widehat{\bV}_{-k,-j}^{\T}(\widehat{\bV}_k-\widehat{\bV}_{-k,-j}\widehat{\bD}_{-k,-j}^{-1}  \bB_{kj})\|_{\max}^{(d)} \leq \tau_n/2$.
\end{lemma}

\textit{Proof.} The proof follows from  Proposition 7 in \textcolor{blue}{Fang et al. (2022)} and hence  is omitted here.
$\hfill\Box$

\section{Additional empirical results}
\label{supp.sec_emp}
\subsection{Simulation studies}
Table~\ref{table.sim.12.0.05} and \ref{table.sim.34.0.05} present  the empirical FDRs and powers of the proposed and BC methods over $100$ replications at target FDR level $\alpha = 0.05$ under Models 1--2, and Models 3--4, respectively. The same patterns as those from Tables~\ref{table.sim.12}--\ref{table.sim.34} can be observed. Note that the results of hard and soft thresholding methods do not depend on $\alpha$, and hence are not reported here.
\begin{table}[!htbp]

	\caption{Empirical FDRs  ($\%$) and powers ($\%$)  of  MFC and BH methods at the $5\%$ nominal level for Models 1 and 2 over 100 simulation runs. \label{table.sim.12.0.05}}
	\begin{center}
		 \vspace{-0.5cm}
		\resizebox{4.5in}{!}{
\begin{tabular}{cccccccccccc}
\hline
\multirow{3}{*}{$n$} & \multirow{3}{*}{$p$} & \multirow{3}{*}{Scenario} & \multicolumn{4}{c}{Model 1}                     &  & \multicolumn{4}{c}{Model 2}                     \\ \cline{4-7} \cline{9-12}
                     &                      &                           & \multicolumn{2}{c}{MFC} & \multicolumn{2}{c}{BC} &  & \multicolumn{2}{c}{MFC} & \multicolumn{2}{c}{BC} \\
                     &                      &                           & FDR      & Power     & FDR      & Power     &  & FDR      & Power     & FDR      & Power     \\ \hline
100                  & 30                   & Fully                     & 4.10       & 88.00         & 0.11      & 64.75      &  & 3.88      & 83.50       & 0.10       & 59.48      \\
                     &                      & $T_{ij}$ = 51             & 4.37      & 85.79      & 0.11      & 62.68      &  & 4.12      & 82.00         & 0.22      & 55.32      \\
                     &                      & $T_{ij}$ = 25             & 4.28      & 83.35      & 0.14      & 60.81      &  & 3.85      & 79.95      & 0.24      & 52.55      \\
                     & 60                   & Fully                     & 4.59      & 82.63      & 0.07      & 58.05      &  & 4.00         & 77.32      & 0.03      & 48.94      \\
                     &                      & $T_{ij}$ = 51             & 4.84      & 79.75      & 0.12      & 56.06      &  & 4.48      & 74.91      & 0.09      & 45.16      \\
                     &                      & $T_{ij}$ = 25             & 4.86      & 77.20       & 0.09      & 54.26      &  & 4.76      & 72.96      & 0.10       & 41.68      \\ \hline
200                  & 30                   & Fully                     & 4.56      & 99.67      & 0.07      & 93.32      &  & 4.77      & 98.11      & 0.12      & 90.07      \\
                     &                      & $T_{ij}$ = 51             & 4.68      & 99.4       & 0.11      & 90.95      &  & 5.00         & 97.36      & 0.15      & 88.18      \\
                     &                      & $T_{ij}$ = 25             & 4.49      & 98.93      & 0.10       & 87.81      &  & 4.57      & 96.77      & 0.11      & 86.11      \\
                     & 60                   & Fully                     & 4.59      & 99.34      & 0.06      & 87.44      &  & 4.71      & 95.04      & 0.08      & 80.85      \\
                     &                      & $T_{ij}$ = 51             & 4.65      & 98.65      & 0.05      & 83.75      &  & 4.96      & 93.98      & 0.06      & 78.33      \\
                     &                      & $T_{ij}$ = 25             & 4.50       & 97.77      & 0.03      & 80.29      &  & 4.92      & 93.09      & 0.09      & 75.92      \\ \hline
                     \end{tabular}
		}	
	\end{center}
	 \vspace{-0.5cm}
                      \caption{Empirical FDRs  ($\%$) and powers ($\%$)  of MFG and BH methods at the $5\%$ nominal level for Models 3 and 4  over 100 simulation runs. \label{table.sim.34.0.05}}
	\begin{center}
	\vspace{-0.2cm}
		\resizebox{4.5in}{!}{
\begin{tabular}{cccccccccccc}
\hline
                     \multirow{3}{*}{$n$} & \multirow{3}{*}{$p$} & \multirow{3}{*}{Scenario} & \multicolumn{4}{c}{Model 3}                     &  & \multicolumn{4}{c}{Model   4}                   \\ \cline{4-7} \cline{9-12}
                     &                      &                           & \multicolumn{2}{c}{MFG} & \multicolumn{2}{c}{BC} &  & \multicolumn{2}{c}{MFG} & \multicolumn{2}{c}{BC} \\
                     &                      &                           & FDR      & Power     & FDR      & Power     &  & FDR      & Power     & FDR      & Power     \\ \hline
200& 30                   & Fully                     & 4.58      & 60.49      & 0.48      & 26.04      &  & 4.26      & 69.56      & 0.29      & 46.31      \\
                     &                      & $T_{ij}$ = 51             & 4.55      & 59.05      & 0.38      & 24.51      &  & 4.57      & 68.61      & 0.31      & 44.86      \\
                     &                      & $T_{ij}$ = 25             & 4.46      & 57.61      & 0.17      & 24.12      &  & 4.15      & 68.00         & 0.30       & 44.42      \\
                     & 60                   & Fully                     & 4.94      & 49.09      & 0.44      & 15.27      &  & 4.86      & 65.32      & 0.27      & 37.44      \\
                     &                      & $T_{ij}$ = 51             & 4.62      & 47.83      & 0.50       & 14.48      &  & 4.74      & 64.28      & 0.19      & 36.29      \\
                     &                      & $T_{ij}$ = 25             & 5.22      & 46.14      & 0.29      & 13.72      &  & 4.91      & 63.57      & 0.28      & 35.43      \\ \hline
400& 30                   & Fully                     & 4.44      & 91.61      & 0.05      & 62.25      &  & 4.77      & 94.25      & 0.30       & 80.19      \\
                     &                      & $T_{ij}$ = 51             & 4.01      & 90.47      & 0.06      & 59.49      &  & 4.79      & 93.50       & 0.28      & 78.39      \\
                     &                      & $T_{ij}$ = 25             & 4.13      & 90.11      & 0.03      & 58.68      &  & 5.16      & 93.53      & 0.24      & 77.61      \\
                     & 60                   & Fully                     & 4.95      & 88.57      & 0.06      & 52.99      &  & 4.55      & 92.03      & 0.09      & 71.29      \\
                     &                      & $T_{ij}$ = 51             & 4.82      & 87.50       & 0.08      & 50.38      &  & 4.66      & 90.95      & 0.06      & 69.36      \\
                     &                      & $T_{ij}$ = 25             & 5.02      & 85.92      & 0.08      & 48.44      &  & 4.70       & 90.31      & 0.08      & 67.85      \\ \hline

\end{tabular}
		}	
	\end{center}
\end{table}

\subsection{Real data analysis}
\label{supp.ad}
\begin{figure}[!htbp]
	\centering
	\includegraphics[width=1\textwidth]{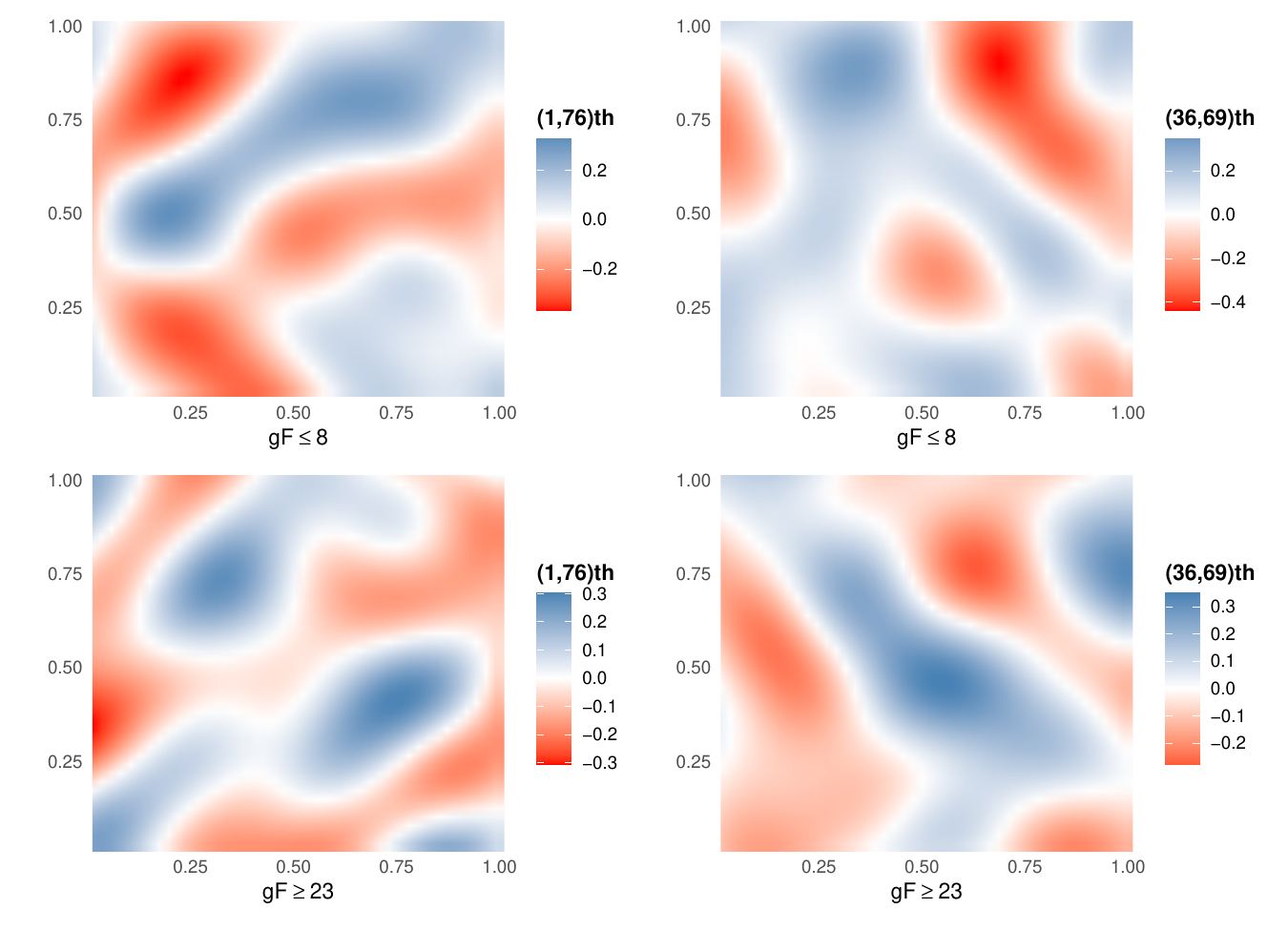}
	\caption{\label{HCP_cov.est}{ Selected significant elements of $\{\widehat\Sigma_{jk}(\cdot, \cdot): 1 \leq j < k\leq p\}$
 in Figure \ref{hcp_network}.}}
\end{figure}

\begin{figure}[!htbp]
	\centering
	\includegraphics[width=1\textwidth]{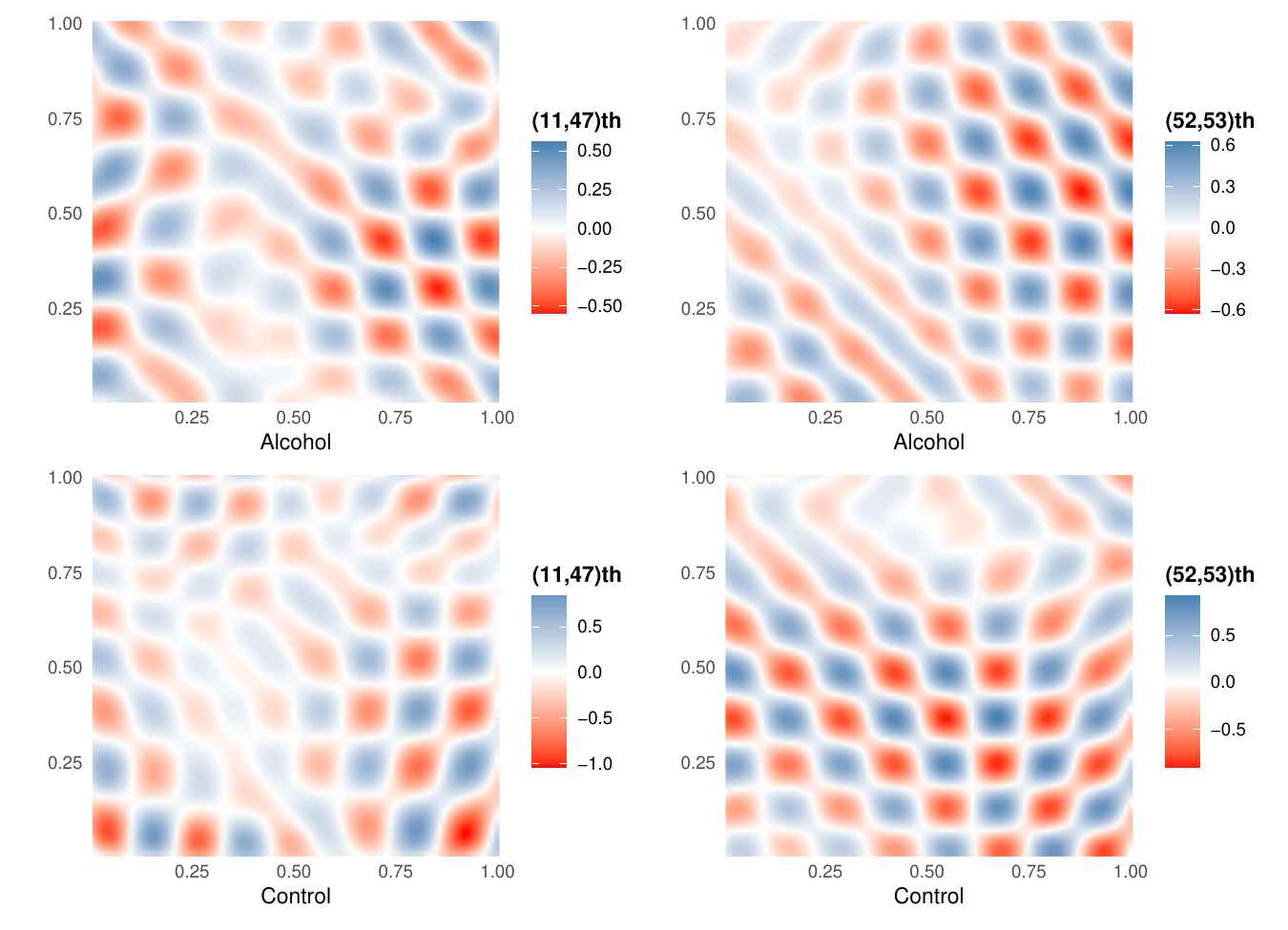}
	\caption{\label{EEG_fgm.est}{ Selected significant elements of $\{\widetilde\Sigma_{jk}^{\varepsilon}(\cdot,\cdot): 1 \leq j < k\leq p\}$
 in Figure \ref{eeg_network}.}}
\end{figure}

We present in Figure~\ref{HCP_cov.est} some randomly selected entries of $\widehat\bSigma(\cdot, \cdot)$ in Figure \ref{hcp_network} of the HCP data. The signal of the estimated covariance surfaces is observable over the entire interval $\cU$ in the sense that the values of $\widehat \Sigma_{jk}$ change from positive to negative (or from negative to positive) as a function of $(u,v)$ smoothly along some directions. Figure~\ref{EEG_fgm.est} displays some randomly selected elements of $\{\widetilde\Sigma_{jk}^{\varepsilon}(\cdot,\cdot), 1 \leq j < k\leq p\}$ in Figure \ref{eeg_network} of the EEG data, highlighting again the smoothly varying pattern of signals for the brain imaging dataset.

\newpage
\spacingset{1.2}
\section*{References}

\begin{description}

\item Berman, S. M. (1962). 
\newblock A law of large numbers for the maximum in a stationary Gaussian sequence, {\it Annals of Mathematical Statistics} {\bf 33}: 93–97.

\item Billingsley, P. (1999). \textit{Convergence of Probability Measures}, 2nd, Wiley, New York.

\item Boucheron, S., Lugosi, G. and Massart, P. (2014). \textit{Concentration inequalities: A Nonasymptotic Theory of Independence}, Oxford University Press.

\item Dauxois, J., Pousse, A. and Romain, Y. (1982). \newblock Asymptotic theory for the principal component analysis of a vector random function: Some applications to statistical inference, {\it Journal of Multivariate Analysis} {\bf 12}: 136–154.

\item Fan, J. and Gijbels, I. (1996).\textit{ Local polynomial modelling and its applications}, Chapman
and Hall, London.

    \item 
    Fang, Q., Guo, S. and Qiao, X. (2022).
	\newblock   Finite sample theory for high-dimensional functional/scalar time series with applications, {\it Electronic Journal of Statistics} {\bf 16}: 527-591.

\item Guo, S., Li, D., Qiao, X. and Wang, Y. (2023). From sparse to dense functional data in high dimensions: Revisiting phase transitions from a non-asymptotic perspective, \textit{arXiv:2306.00476v2}.

\item Guo, S. and Qiao, X. (2023). On consistency and sparsity for high-dimensional functional
time series with application to autoregressions, \textit{Bernoulli} \textbf{29}: 451–472.

\item Lin, Z. and Bai, Z. (2010). \textit{Probability Inequalities}, Science Press Beijing, Beijing.

\item Vershynin, R. (2018). \textit{High-Dimensional Probability: An Introduction with Applications in
Data Science}, Cambridge University Press, Cambridge.

\item Zaitsev, A. Y. (1987). 
\newblock On the Gaussian approximation of convolutions under multidimensional analogues of S.N. Bernstein's inequality conditions, {\it Probability Theory and Related Fields} {\bf 74}: 535–566.

\item Zapata, J., Oh, S. Y. and Petersen, A. (2022). Partial separability and functional graphical
models for multivariate Gaussian processes, \textit{Biometrika} \textbf{109}: 665–681.

\item Zhang, J.-T. (2013). \textit{Analysis of Variance for Functional Data}, CRC Press.

\item Zhang, J.-T. and Chen, J. (2007). Statistical inferences for functional data, \textit{The Annals of
Statistics} \textbf{35}: 1052–1079.

\end{description}
\end{document}